\documentclass[11pt]{amsart}
\usepackage{t1enc}
\usepackage[latin1]{inputenc}
\usepackage[english]{babel}
\usepackage{amsmath,amsthm}
\usepackage{amsfonts}
\usepackage{latexsym}
\usepackage[dvips]{graphicx}
\usepackage{float}
\usepackage[natural]{xcolor}
\usepackage{algorithm}
\usepackage{algorithmic}
\usepackage{enumerate}
\usepackage{multirow}
\usepackage{xcolor}
\usepackage[colorlinks,linkcolor=blue]{hyperref}
\usepackage{lineno}
\usepackage{setspace}
\usepackage{amssymb}
\usepackage{fullpage}
\usepackage{tikz}
\usepackage{stfloats}
\usepackage{selinput}
\usepackage{babel}
\usepackage{nicefrac}
\usepackage{lmodern}
\usetikzlibrary{decorations.pathreplacing}
\usetikzlibrary{patterns}
\usetikzlibrary{positioning} 
\usepackage{xcolor}
\usepackage{listings}
\newtheorem{theorem}{Theorem}[section]
\newtheorem{lemma}[theorem]{Lemma}

\newtheorem{conjecture}{Conjecture}[section]

\newtheorem{definition}[theorem]{Definition}
\newtheorem{proposition}[theorem]{Proposition}

\newtheorem{claim}[theorem]{Claim}
\newtheorem{corollary}{Corollary}[section]

\def\int{\textrm{int}}

\begin{document}

\onehalfspacing
\title{Rainbow Hamilton cycle in hypergraph system}

\author{Yucong Tang}
\address{YT. Department of Mathematics, Nanjing University of Aeronautics and Astronautics, Nanjing, China.
Key Laboratory of Mathematical Modelling and High Performance Computing of Air Vehicles (NUAA), Nanjing, China. \\
Email: \texttt{tangyucong@nuaa.edu.cn}.}

\author{Bin Wang}
\address{BW and GW. School of Mathematics, Shandong University, Jinan, China \\
Email: \texttt{(BW) binwang@mail.sdu.edu.cn, (GW) ghwang@sdu.edu.cn}.}
\author{Guanghui Wang}
\author{Guiying Yan}
\address{GY. Academy of Mathematics and Systems Science, Chinese
Academy of Sciences, Beijing, China \\
Email: \texttt{yangy@amss.ac.cn}.}

\begin{abstract}

In this paper,  we develop a new rainbow Hamilton framework, which is of independent interest, settling the problem proposed by Gupta, Hamann, M\"{u}yesser, Parczyk, and Sgueglia when $k=3$, and draw the general conclusion for any $k\geq3$ as follows.
A $k$-graph system $\emph{\textbf{H}}=\{H_i\}_{i\in[n]}$ is a family of not necessarily
distinct $k$-graphs on the same $n$-vertex set $V$, moreover, a $k$-graph $H$ on $V$ is rainbow if $E(H)\subseteq \bigcup_{i\in[n]}E(H_i)$
and $|E(H)\cap E(H_i)|\leq1$ for $i\in[n]$.
We show that given $\gamma> 0$, sufficiently large $n$ and
an $n$-vertex $k$-graph system $\emph{\textbf{H}}=\{H_i\}_{i\in[n]}$
, if $\delta_{k-2}(H_i)\geq(5/9+\gamma)\binom{n}{2}$ for $i\in[n]$ where $k\geq3$, then there exists a rainbow tight Hamilton cycle.
This result implies the conclusion in a single graph, which was proved by Lang and Sanhueza-Matamala [\emph{J. Lond. Math. Soc., 2022}], Polcyn, Reiher, R\"{o}dl and Sch\"{u}lke [\emph{J. Combin. Theory Ser. B, 2021}] independently.
\bigskip

\end{abstract}

\maketitle
\section{Introduction}
Finding Hamilton cycles in graphs is one of the key areas in graph theory and extremal combinatorics with a profound history.
The classical Dirac's theorem \cite{Dirac} states that every $n$-vertex graph with minimum degree at least $n/2, n\geq3$, contains a Hamilton cycle.
There are also many extensions of Dirac's theorem in hypergraphs.

\subsection{Hamilton cycles in hypergraphs}
Let $[a,b]$, $a,b\in \mathbb{Z}$, denote the set $\{a,a+1,\ldots,b\}$ and the set $[1,n]$ is denoted by $[n]$ in short.
Given a $k$-graph $H$ with a set $S$ of $d$ vertices ($d\in[k-1]$), we define $\deg_H(S)$ to be the number of edges containing $S$ (the subscript $H$ is omitted if it is clear from the context), the relative degree $\overline{\deg}(S)$ to be $\deg(S)/\binom{n-d}{k-d}$.
The \emph{minimum relative} $d$-\emph{degree} of a $k$-graph $H$, written by $\overline{\delta}_d(H)$, is the minimum of $\overline{\deg}(S)$ over all sets $S$ of $d$ vertices.

Katona and Kierstead \cite{MR1671170} defined a type of cycle in hypergraphs, which has been studied extensively.
A $k$-graph is called an $\ell$-cycle if its vertices can be ordered cyclically such that each of its edges consists of $k$ consecutive vertices and every two consecutive edges (in the natural order of the edges) share exactly $\ell$ vertices.
In $k$-graphs, a $(k-1)$-cycle is often called a tight cycle.
We say that a $k$-graph contains a Hamilton $\ell$-cycle if it contains an $\ell$-cycle as a spanning subhypergraph.
Without special instruction, the tight cycle is referred as cycle for short.

Katona and Kierstead \cite{MR1671170} gave a sufficient condition for finding a Hamilton cycle in a $k$-graph with minimum $(k-1)$-degree:
every $n$-vertex $k$-graph $H$ with $\delta_{k-1}(H)>(1-1/(2 k))n+4-k-5/(2k)$ admits a Hamilton cycle.
They conjectured that the bound on the minimum $(k-1)$-degree can be reduced to roughly $n/2$, which was confirmed asymptotically by R\"{o}dl, Ruci\'{n}ski and Szemer\'{e}di in \cite{3uniform,approximate}.
The same authors gave the exact version for $k=3$ in \cite{EXACT}.

\begin{theorem}[\cite{approximate,EXACT}]
\label{main1}
Let $k\geq3, \gamma>0$ and $H$ be an $n$-vertex $k$-graph, where $n$ is sufficiently large. If $\delta_{k-1}(H)\geq(1/2+\gamma)n$, then $H$ contains a Hamilton cycle.
Furthermore, when $k=3$ it is enough to have $\delta_2(H)\ge \lfloor n/2\rfloor$.
\end{theorem}

More generally, K\"{u}hn and Osthus~\cite{D2014Hamilton} and Zhao~\cite{2015Recent} noted that it is much more difficult to determine the minimum $d$-degree condition for tight Hamilton cycle for $d\in[k-2]$.
Based on the results of Cooley and Mycroft~\cite{Oliver2017The}, Glebov, Person and Weps~\cite{2011On}, R\"{o}dl and Ruci\'{n}ski~\cite{2014Families} and R\"{o}dl, Ruci\'{n}ski, Schacht and Szemer\'{e}di~\cite{Rodl2017}, Reiher, R\"{o}dl, Ruci\'{n}ski, ~Schacht, and Szemer\'{e}di~\cite{2019Minimum} gave the asymptotic version when $d=k-2$ and $k=3$, while Polcyn, Reiher, R\"{o}dl, Ruci\'{n}ski, Schacht, and Sch\"{u}lke~\cite{2020Minimum} gave the asymptotic version when $d=k-2$ and $k=4$.
Glebov, Person and Weps~\cite{2011On} proved the minimum relative $d$-degree condition for a tight Hamilton cycle is a function of $k$.
The best general bound was given by Lang and Sanhueza-Matamala~\cite{Lang}, Polcyn, Reiher, R\"{o}dl and Sch\"{u}lke~\cite{2021On} independently.
They proved the following theorem.
\begin{theorem}[\cite{Lang,2021On}]
Let $k\geq3$, $\gamma>0$ and $H$ be an $n$-vertex $k$-graph, where $n$ is sufficiently large.
If $\delta_{k-2}(H)\geq(5/9+\gamma)\binom{n}{2}$, then $H$ contains a Hamilton cycle.
\end{theorem}

A construction due to Han and Zhao~\cite{2015Forbidding} showed that the constant $5/9$ appearing in the above theorem is optimal.
For more background, we refer the readers to the recent surveys of K\"{u}hn and Osthus~\cite{D2014Hamilton}, R\"{o}dl and Ruci\'{n}ski \cite{30}, Simonovits and Szemer\'{e}di \cite{39} and Zhao \cite{2015Recent}.

\subsection{Rainbow settings in hypergraph systems}
A $k$-graph system $\emph{\textbf{H}}=\{H_i\}_{i\in[n]}$ is a family of not necessarily
distinct $k$-graphs on the same $n$-vertex set $V$, moreover, a $k$-graph $H$ on $V$ is rainbow if $E(H)\subseteq \bigcup_{i\in[n]}E(H_i)$
and $|E(H)\cap E(H_i)|\leq1$ for $i\in[n]$.
Let $|H|$ denote the size of the vertex set of $H$.

The study of rainbow structures in graph systems has attracted much more attention.
Aharoni, DeVos, Maza, Montejano, and \v{S}\'{a}mal.~\cite{MR4125343} conjectured that: for $|V|=n\geq 3$ and $n$-vertex graph system $\emph{\textbf{G}}=\{G_i\}_{i\in[n]}$ on $V$, if $\delta(G_i)\geq n/2$ for each $i\in[n]$, then there exists a rainbow Hamilton cycle with edge set $\{e_1,\ldots,e_n\}$ such that $e_i\in E(G_i)$ for $i\in[n]$.
This was recently verified asymptotically by Cheng, Wang and Zhao~\cite{2019Rainbow}, and completely by Joos and Kim~\cite{MR4171383}.
In~\cite{from}, Bradshaw, Halasz, and Stacho strengthened the Joos-Kim result by showing that given an $n$-vertex graph system $\emph{\textbf{G}}=\{G_i\}_{i\in[n]}$ with $\delta(G_i)\geq n/2$ for $i\in[n]$, then $\emph{\textbf{G}}$ has exponentially many rainbow Hamilton cycles.
Similarly, a degree condition of Moon and Moser \cite{Moon} for Hamiltonicity in bipartite graphs has been generalized to the rainbow setting by Bradshaw in \cite{Bradshaw}.
Generally, for each graph $F$, let $\delta_F$ be the smallest real number $\delta\geq 0$ such that, for each
$\varepsilon>0$ there exists some $n_0$ such that, for every $n\geq n_0$ with $|F|$ dividing $n$, if an $n$-vertex graph $G$ has minimum
degree at least $(\delta+\varepsilon)n$, then $G$ contains an $F$-factor.
Cheng, Han, Wang, and Wang~\cite{2020B} proved that the minimum degree bound $\delta_{K_r}$ is asymptotically sufficient for the existence of rainbow $K_r$-factor in graph systems.
Montgomery, M\"{u}yesser and  Pehova~\cite{factor} generalized the above conclusion for some $F$ satisfying $\delta_F \geq1/2$ or $F$ has a bridge.

In hypergraph systems, Cheng, Han, Wang, Wang and Yang~\cite{CHWWY} proved that given $k\geq3, \gamma>0$, sufficiently large $n$ and an $n$-vertex $k$-graph system $\emph{\textbf{H}}=\{H_i\}_{i\in[n]}$, if $\delta_{k-1}(H_i)\geq(1/2+\gamma)n$ for $i\in[n]$, then there exists a rainbow tight Hamilton cycle.
There are also some works on rainbow subgraphs, see~\cite{Preprint,MR4157094,2012The,2021G,New,2020A,LU,Hongliang2018ON,2020C,2020B,MR4055023}.
Recently, Gupta, Hamann, M\"{u}yesser, Parczyk and Sgueglia~\cite{2022new} gave a unified approach to this problem.
However, they mentioned that ``there is a well-known (uncoloured) Dirac-type result whose rainbow version is missing'' (Given a $3$-graph system $\emph{\textbf{H}}=\{H_i\}_{i\in[n]}$, with minimum vertex degree condition of each $H_i$, does $\emph{\textbf{H}}$ admit a rainbow Hamilton cycle?) and  ``it would be an interesting challenge to obtain this result''.
It hits a technical barrier to tackle.
In this paper, we develop the new rainbow Hamilton framework, whose uncolored version is first established in ~\cite{Lang}, and give a general result as follows.

\begin{theorem}
\label{main}
For every $k\geq3, \gamma>0$, there exists $n_0$ such that the following holds for $n\geq n_0$.
Given a $k$-graph system $\textbf{H}=\{H_i\}_{i\in[n]}$, if $\delta_{k-2}(H_i)\geq(5/9+\gamma)\binom{n}{2}$ for $i\in[n]$, then $\textbf{H}$ admits a rainbow Hamilton cycle.
 \end{theorem}

\subsection{Notation and preliminary}
We call a hypergraph $H$ a $(1,k)$-graph if $V(H)$ can be partitioned into $V_1$ and $V_2$ such that every edge contains exactly one vertex of $V_1$ and $k$ vertices of $V_2$.
Given a partition $V(H)=V_1\cup V_2$, a $(1,d)$-subset $S$ of $V(H)$ contains one vertex in $V_1$ and $d$ vertices in $V_2$.
Let $\delta_{1,d}(H):=\min\{\deg_H(S):S {\rm\ is\ a }\ (1,d)$-${\rm subset\ of}\ V(H)\}$ for $d\in[k-1]$.

A $k$-partite graph is a graph whose vertices are (or can be) partitioned into $k$ different independent sets.
Given a $(k+1)$-partite $(k+1)$-graph $H$ with $V(H)=V_0\cup V_1\cup\cdots\cup V_k$.
A $(k+1)$-uniform sequentially path $P$ of \emph{length} $t$ in $H$ is a $(k+1)$-graph with vertex set $V(P)=C(P)\cup I(P)$ where $C(P)=\{c_1,\ldots,c_{t-k+1}\}\subseteq V_0$, $I(P)=\{v_1,\ldots,v_{t}\}\subseteq V_1\cup\cdots \cup V_k$ and edge set $\{e_1,\ldots,e_{t-k+1}\}$ such that $e_i=\{c_i,v_i,\ldots,v_{i+k-1}\}$ for $i\in[t-k+1]$.
Denote the length of $P$ by $\ell(P)$.
We call $c_1,\ldots,c_{t-k+1}$ the \emph{colors} of $P$ and $v_1,\ldots,v_t$ the \emph{points} of $P$.
For convenience, we use $(C(P),I(P))$ to denote the above sequentially tight path.
Furthermore, if $(v_1,\ldots,v_t)$ is a cyclically ordered set, then we call this sequentially path a \emph{sequentially cycle}.
A $(k+1)$-uniform sequentially walk is an ordered set of points with an ordered set of colors such that the $i_{th}$ $k$ consecutive points along with the $i_{th}$ color forms an edge.
Note that the points, edges and colors in a sequentially walk are allowed to be repeated.
The length of a sequentially walk is its number of points.

Before we give the proof of Theorem \ref{main}, we use the following similar definitions with \cite{Lang}.
\begin{definition}[Sequentially Hamilton cycle threshold]
The minimum $(1,k-2)$-degree threshold for sequentially Hamilton cycles, denoted by $thc_{k-2}(k)$, is the smallest number $\delta>0$ such that,
for every $\varepsilon>0$, there exists an $n_0\in \mathbb{N}$ such that every $(1,k)$-graph $H$ on $[n]\cup V$ with minimum degree $\delta_{1,k-2}(H)\geq(\delta+\varepsilon)\binom{n}{2}$ contains a sequentially Hamilton cycle where $|V|=n\geq n_0$.
\end{definition}

\begin{definition}[Sequentially tight connectivity]
A subgraph $H'$ of a $(1,k)$-graph $H$ is sequentially tightly connected, if any two edges of $H'$ can be connected by a sequentially walk.
A sequentially tight component of $H$ is an edge maximal sequentially tightly connected subgraph.
\end{definition}

Given $\textbf{b}$: $V(H)\rightarrow[0,1]$, we define the  $\textbf{b}$-\emph{fractional matching} to be a function $\textbf{w}$: $E(H)\rightarrow[0,1]$ such that $\sum_{e:v\in e}\textbf{w}(e)\leq\textbf{b}(v)$ for every vertex $v\in V(H)$.
Moreover, if the equality holds, then we call $\textbf{w}$ perfect.
Denote the maximum size of a $\textbf{b}$-fractional matching by $\nu(H,\textbf{b})=\max_{\textbf{w}}\sum_{e\in E(H)}\textbf{w}(e)$ where $\textbf{w}$ is a $\textbf{b}$-fractional matching.
It is well-known that perfect matchings are closely related to its fractional counterpart.
In particular, when $\textbf{b}\equiv1$, the $\textbf{b}$-\emph{fractional matching} is called \emph{fractional matching}.
The \emph{density} of a $\textbf{b}$-fractional matching is $\sum_{e\in E(H)}\textbf{w}(e)/|V(H)|$.
Besides, we require the following characterization.
Given a $k$-graph $H$, we say that $H$ is $\gamma$-\emph{robustly matchable} if the following holds.
For every vertex weight $\textbf{b}$: $V(H)\rightarrow[1-\gamma,1]$, there is an edge weight $\textbf{w}$: $E(H)\rightarrow[0,1]$ with $\sum_{e:v\in e}\textbf{w}(e)=\textbf{b}(v)/(k-1)$ for every vertex $v\in V(H)$.
Note that a $\gamma$-robustly matchable $k$-graph $H$ admits a $\textbf{b}$-fractional matching of size $\sum_{v\in V(H)}\textbf{b}(v)/k(k-1)$ for every vertex weighting $\textbf{b}$: $V(H)\rightarrow[1-\gamma,1]$.
The following definition plays an important role in our proof.
\begin{definition}[Link graph]
Consider a $(1,k)$-graph $H$ on $V(H)=[n]\cup V$ where $|V|=n$ and a set $S$ of $(1,\ell)$-subset of $V(H)$.
We define the link $(k-\ell)$-graph of $S$ in $H$ as the graph $L_H(S)$ with vertex set $V$ and edge set $\{X: X\cup S\in E(H)\}$ for $\ell\in[0,k-1]$.
If $H$ is clear, then we simply write $L(S)$.
\end{definition}

Let $H=(V,E)$ be a $k$-graph, $V'\subseteq V$, an \emph{induced subgraph} $H[V']$ of a $k$-graph $H$ is a $k$-graph with vertex set $V'$ and edge set $E'$ where each edge is precisely the edge of $H$ consisting of $k$ vertices in $V'$.
We usually denote $H'$ by $H[V']$.

\begin{definition}[Rainbow Hamilton framework]
\label{1.7}
Let $\alpha,\gamma,\delta$ be positive constants.
Suppose $R$ is a $(1,k)$-graph on $[t]\cup V$ where $|V|=t$, we call a subgraph $H$ of $R$ an $(\alpha,\gamma,\delta)$-rainbow Hamilton framework, if $H$ has the following properties.
\begin{itemize}
\item[(F1)] $H_i:=H[\{i\}\cup V]$ is sequentially tightly connected for $i\in[t]$,
\item[(F2)] $H_i$ contains a sequentially closed walk of length 1 mod $k$ for $i\in[t]$,
\item[(F3)] $H_{W_i}:=H[[t(i-1)/k+1,ti/k]\cup V]$ is $\gamma$-robustly matchable for $i\in[k]$,
\item[(F4)] For every color $i\in[t]$, there are at least $(1-\alpha)t$ points $v\in V$ such that $\{i,v\}$ has relative  $(1,1)$-degree at least $1-\delta+\gamma$,
\item[(F5)] $L_H(\{i\})$ and $L_H(\{j\})$ intersect in an edge for each $i,j\in[t]$.
\end{itemize}
\end{definition}
We write $x\ll y$ to mean that for any $y\in(0,1]$, there exists an $x_0\in(0,1)$ such that for all $x\leq x_0$, the subsequent statements hold.
Hierarchies with more constants are defined similarly to be read from right to left.
\begin{definition}[Rainbow Hamilton framework threshold]
\label{1.8}
The minimum $(1,k-2)$-degree threshold for $(1,k)$-uniform rainbow Hamilton framework, denoted by $rhf_{k-2}(k)$, is the smallest value of $\delta$ such that the following holds.

Suppose $\varepsilon,\alpha,\gamma,\mu>0$ and $t\in \mathbb{N}$ with $1/t\ll\varepsilon\ll\alpha\ll\gamma\ll\mu$.
If $R$ is a $(1,k)$-graph on $[t]\cup V$ where $|V|=t$, with minimum relative $(1,k-2)$-degree at least $\delta+\mu$ and a set $I\subseteq E(R)$ of at most $\varepsilon t\binom{t}{k}$ perturbed edges, then $R$ contains an $(\alpha,\gamma,\delta)$-rainbow Hamilton framework $H$ that avoids the edges of $I$.
\end{definition}

We transform the problem of bounding the sequentially Hamilton cycle threshold to bound the rainbow Hamilton framework threshold.
\begin{theorem}[Framework Theorem]
\label{1.9}
For $k\geq3$, we have $thc_{k-2}(k)\leq rhf_{k-2}(k)$.
\end{theorem}

Let the \emph{shadow graph} $\partial_j(H)$ of $(1,k)$-graph $H$ at level $j$ be the $(1,j)$-graph on $[n]\cup V$ whose edges are $(1,j)$-sets contained in the edges of $H$ for $j\in[k]$.
\begin{definition}[Vicinity]
\label{defvicinity}
Given a $(1,k)$-graph $R$ on $[t]\cup V$, we say that $\mathcal{C}_i=\{C_S\subseteq L(S):S\in\partial_{k-2}(R){\rm\ and\  }i\in S\}$ for each $i\in[t]$ is a $(k-2)$-vicinity.
We define the $(1,k)$-graph $H$ generated by $\mathcal{C}_i$ as the subgraph of $R$ with vertex set $V(H)=\{i\}\cup V$ and edge set
\[
E(H)=\bigcup_{i\in S, S\in\partial_{k-2}(R)}\{A\cup S:A\in C_S\}.
\]
\end{definition}
Besides, we need the following structures.
\begin{definition}[Switcher]
\label{switcher}
A switcher in a graph $G$ is an edge $ab$ such that $a$ and $b$ shares a common neighbor in $G$.
\end{definition}

Note that a switcher together with its common neighbor generates a triangle.
\begin{definition}[Arc]
\label{arc}
Let $R_i$ be a $(1,k)$-graph on $\{i\}\cup V$ with $(k-2)$-vicinity $\mathcal{C}_i=\{C_S:S\in\partial_{k-2}(R_i)\}$.
We say that a $(1,k+1)$-tuple $(i,v_1,\ldots,v_{k+1})$ is an arc for $\mathcal{C}_i$ if the following holds.
\begin{itemize}
  \item $\{i,v_1,\ldots,v_{k-2}\}\in\partial_{k-2}(R_i)$ with $\{v_{k-1},v_k\}\in C_{\{i,v_1,\ldots,v_{k-2}\}}$.
  \item $\{i,v_2,\ldots,v_{k-1}\}\in\partial_{k-2}(R_i)$ with $\{v_k,v_{k+1}\}\in C_{\{i,v_2,\ldots,v_{k-1}\}}$.
\end{itemize}
\end{definition}

\begin{definition}[Rainbow Hamilton vicinity]
\label{1.15}
Let $\gamma,\delta>0$.
Suppose that $R$ is a $(1,k)$-graph on $[t]\cup V$, let $R_i:=R[\{i\}\cup V]$.
We say that a family $\mathcal{C}=\{\mathcal{C}_i:i\in[t]\}$ of $(k-2)$-vicinities where $\mathcal{C}_i=\{C_S:S\in\partial_{k-2}(R_i)\}$ is $(\gamma,\delta)$-rainbow Hamilton if for any $S,S'\in \partial_{k-2}(R_i)$ and $T\in \partial_{k-2}(R_j)$ where $i\neq j$, the followings hold,
\begin{itemize}
  \item[(V1)] $C_S$ is tightly connected,
  \item[(V2)] $C_S$ and $C_{S'}$ intersect in an edge,
  \item[(V3)] $C_S$ has a switcher and the vicinity $\mathcal{C}_i$ has an arc for $i\in[t]$,
  \item[(V4)] $C_S$ has a fractional matching of density $ (1+1/k)(1/(k+1)+\gamma)$,
  \item[(V5)] $C_S$ has edge density at least $1-\delta+\gamma$,
  \item[(V6)] $C_{S}$ and $C_T$ intersect in an edge.
\end{itemize}
\end{definition}
\begin{definition}[Perturbed degree]
\label{1.16}
Let $\alpha,\delta>0$.
We say that a $(1,k)$-graph $R$ has $\alpha$-perturbed minimum relative $(1,k-2)$-degree at least $\delta$ if the followings hold for $j\in[k-2]$.
\begin{itemize}
  \item[(P1)] every edge of $\partial_{j}(R)$ has relative degree at least $\delta$ in $R$,
  \item[(P2)] $\overline{\partial_{j}(R)}$ has edge density at most $\alpha$, where $\overline{\partial_{j}(R)}$ denotes the complement of $\partial_{j}(R)$,
  \item[(P3)] each $(1,j-1)$-tuple of $\partial_{j-1}(R)$ has relative degree less than $\alpha$ in $\overline{\partial_{j}(R)}$.
\end{itemize}
\end{definition}
\begin{definition}[Rainbow Hamilton vicinity threshold]
\label{1.17}
The minimum $(1,k-2)$-degree threshold for $(1,k)$-uniform rainbow Hamilton vicinities, denoted by $rhv_{k-2}(k)$, is the smallest value $\delta>0$ such that the following holds.
Let $\alpha,\gamma,\mu>0$, $t\in \mathbb{N}$ with $1/t\ll\alpha\ll\gamma\ll\mu$ and $R$ be a $(1,k)$-graph on $[t]\cup V$.
If each $R_i:=R[\{i\}\cup V]$ has $\alpha$-perturbed minimum relative $(1,k-2)$-degree at least $\delta+\mu$ for $i\in[t]$, then $R$ admits a family of $(\gamma,\delta)$-rainbow Hamilton $(k-2)$-vicinities.
\end{definition}
\begin{theorem}[Vicinity Theorem]
\label{vicinity}
For $k\geq3$, $rhf_{k-2}(k)\leq rhv_{k-2}(k)$.
\end{theorem}

Combining Theorem \ref{vicinity} with Theorem \ref{1.9}, we just need to prove the following theorem, and we can obtain Theorem \ref{main}.
\begin{theorem}
\label{1.20}
For $k\geq3$, $rhv_{k-2}(k)\leq5/9$.
\end{theorem}
We use the following concentration inequalities.
\begin{proposition}[Chernoff's inequality~\cite{Probabilistic}]
\label{chernoff}
Suppose that $X$ has the binomial distribution and $0<a<3/2$,
then $\Pr(|X-\mathbb{E}X|\geq a\mathbb{E}X)\leq2e^{-a^2\mathbb{E}X/3}$.
\end{proposition}
\begin{proposition}[McDiarmid's inequality~\cite{1989On}]
\label{Mc}
Suppose $X_1,\ldots,X_m$ are independent Bernoulli random variables and $b_i\in[0,B]$ for $i\in[m]$.
Suppose that $X$ is a real-valued random variable determined by $X_1,\ldots,X_m$ such that altering the value of $X_i$ changes $X$ by at most $b_i$ for $i\in[m]$.
For all $\lambda>0$, we have
\[
\Pr(|X-\mathbb{E}X|>\lambda)\leq2\exp\left(\frac{-2\lambda^2}{B\Sigma_{i=1}^mb_i}\right).
\]

\end{proposition}
\subsection{Organisation of the paper}
The paper is organised as follows.
In Section 2, we show that how a rainbow Hamilton vicinity deduce a rainbow Hamilton framework.
In Section 3, we show the minimum degree condition guarantees a rainbow Hamilton vicinity.
We review the hypergraph regularity method in Section 4.
In Section 5, the proof of Theorem \ref{1.9} is obtained by the absorption method and almost cover lemma, whose details can be seen in Section 6 and Section 7 respectively.
We conclude the paper with a discussion in Section 8.
For the proof of absorption lemma and almost cover lemma, we develop a new rainbow Hamilton framework.
It is of great interest to tackle the rainbow Hamilton cycle embedding problem with other conditions.
In the proof of absorption lemma, which was widely popularised by R\"{o}dl, Ruci\'{n}ski and Sz\'{e}mer\'{e}di~\cite{3uniform}, we have an innovation point where we provides a mentality for absorbing a color set and a point set.
An absorber can be divided into two parts, one for the color set and the other for the point set.
The almost cover lemma is obtained by tools with regularity.
However, connecting the end-pairs of paths arising in the proof requires more involved changes.
The traditional connecting lemma asserts that for every pair of disjoint pairs of vertices there exists a relatively short tight
path, but there might be pairs of vertices that are not contained in any hyperedge at all as the following example shows.
Consider a 3-graph system $\textbf{\emph{H}}=\{H_i\}_{i\in[n]}$, $V(H_i)=V=X\cup Y$ where $|X|<\frac{1}{3}n$, each $H_i$ has edge set $E=\{e\in V^{(3)}:|X\cap e|\neq2\}$, it is easy to confirm that this 3-graph system satisfies the degree condition in Theorem~\ref{main}, but every tight path starting with a pair of vertices
in $X$ is bound to stay in $X$.
We overcome the obstacle in Section 6.

\section{From vicinity to framework}
Our goal is to prove Theorem \ref{vicinity} in this part.
We need the followings lemmas.
\begin{lemma}
\label{2.1}
Let $R_i$ be a $(1,k)$-graph on $\{i\}\cup V$ with a $(k-2)$-vicinity $\mathcal{C}_i=\{C_S:S\in\partial_{k-2}(R_i)\}$ for $i\in[t]$.
For every $S,S'\in\partial_{k-2}(R_i)$, if the vicinity $\mathcal{C}_i$ has an arc for $i\in[t]$, $C_S$ and $C_{S'}$ intersect, $C_S$ is tightly connected and has a switcher,
then the vertex spanning subgraph $H$ of $R_i$ generated by $\mathcal{C}_i$ is sequentially tightly connected and contains a sequentially closed walk of length 1 mod $k$.
\end{lemma}
\begin{lemma}
\label{2.2}
Let $\gamma,\alpha,\delta>0$ such that $1/t\ll\alpha,\gamma\ll1/k$.
Let $R$ be a $(1,k)$-graph on $[t]\cup V$ where $|V|=t$ and each $R_i$ has $\alpha$-perturbed minimum relative $(1,k-2)$-degree at least $\delta$.
Let $\mathcal{C}=\{\mathcal{C}_i:i\in[t]\}$ be a family of $(k-2)$-vicinities where $\mathcal{C}_i=\{C_S:S\in\partial_{k-2}(R_i)\}$.
If for every $S\in \partial_{k-2}(R)$, $C_S$ has a fractional matching of density $(1+1/k)(1/(k+1)+\gamma)$, then the graph $H\subseteq R$ generated by $\mathcal{C}_{W_i}:=\{\mathcal{C}_j:j\in[t(i-1)/k+1,ti/k]\}$ is $\gamma$-robustly matchable for each $i\in[k]$.
\end{lemma}
\begin{lemma}
\label{2.3}
Let $t,k\in \mathbb{N},i\in[t]$ and $\delta,\alpha,\varepsilon>0$ with $1/t\ll\varepsilon\ll\alpha\ll\delta,1/k$.
Let $R_i$ be a $(1,k)$-graph on $\{i\}\cup V$ with minimum relative $(1,k-2)$-degree at least $\delta$ where $|V|=t$.
Let $I$ be a subgraph of $R_i$ with edge density at most $\varepsilon$, there exists a vertex spanning subgraph $R_i'\subseteq R_i-I$ of $\alpha$-perturbed minimum relative $(1,k-2)$-degree at least $\delta-\alpha$.
\end{lemma}
\begin{proof}[Proof of Theorem \ref{vicinity}]
Let $\delta=rhv_{k-2}(k)$ and $\varepsilon,\alpha,\gamma>0$ with $t_0\in \mathbb{N}$ such that
\[
1/t\ll\varepsilon\ll\alpha\ll\alpha'\ll\gamma\ll\mu\ll\delta,1/k.
\]
Moreover,the constants $t,\varepsilon,\alpha,\mu$ are compatible with the constant hierarchy given by Definition \ref{1.17}, $t,\varepsilon,2\alpha,\mu$ satisfy the conditions of Lemma \ref{2.2} and $t,\varepsilon,\alpha,\delta$ satisfy the conditions of Lemma \ref{2.3}.

Given a $(1,k)$-graph $R_i$ on $\{i\}\cup V$ with minimum relative $(1,k-2)$-degree at least $\delta+2\mu$ and a set $I$ of at most $\varepsilon \binom{t}{k}$ perturbed edges.
We start by selecting a subgraph of $R_i$.
By Lemma \ref{2.3}, we obtain a vertex spanning subgraph $R_i'\subseteq R_i-I$ of $\alpha$-perturbed minimum relative $(1,k-2)$-degree at least $\delta+\mu$.

By the definition of \ref{1.17}, $R':=\bigcup_{i\in[t]}R_i'$ has a family of $(2\gamma,\delta)$-rainbow Hamilton $(k-2)$-vicinities $\mathcal{C}=\{\mathcal{C}_i:i\in[t]\}$ where $\mathcal{C}_i=\{C_S:S\in\partial_{k-2}(R_i)\}$.
Each $\mathcal{C}_i$ generates a $(1,k)$-graph $G_i$.
Let $H=\bigcup_{i\in[t]}G_i$.
Note that $G_i$ does not contain the edges of $I$ and $V(G_i)=V(R_i')$.
By Lemma \ref{2.1} and \ref{2.2}, $H$ also satisfies (F1)-(F3).
For $k\geq4$, by repeatedly applying Definition \ref{1.16}, we deduce that for all but at most $\alpha t$ $(1,1)$-sets of $V(R_i')$ is contained in at least $(1-2\alpha)^{k-3}\binom{|V'|-1}{k-3}\geq(1-2(k-3)\alpha)\binom{|V'|-1}{k-3}$
many $(1,k-2)$-sets in $\partial_{k-2}(R_i')$.
Note that $\partial_{k-2}(R_i')=\partial_{k-2}(G_i)$.
 This implies that for all but at most $\alpha t$ $(1,1)$-tuples of $V(G_i)$ has relative degree at least $1-2(k-3)\alpha$ in $\partial_{k-2}(G_i)$.
Moreover, every $(1,k-2)$-set in $\partial_{k-2}(G_i)$ has relative degree at least $1-\delta+2\gamma$ in $G_i$, since $G_i$ is generated from $(2\gamma,\delta)$-rainbow Hamilton $(k-2)$-vicinity and Definition \ref{1.15}.
Thus, we obtain that for each color $i\in[t]$, there are at least $(1-\alpha)t$ points $v\in V$ such that $\{i,v\}$ has relative $(1,1)$-degree at least $1-\delta+\gamma$, which implies (F4) for $k\geq 4$.
While for $k=3$, by Definition \ref{1.15}, we have every $(1,1)$-set has relative degree at least $1-\delta+2\gamma$ in $G_i$, which implies (F4) for $k=3$.
Besides, it is obvious that (V6) implies (F5), we obtain an $(\alpha,\gamma,\delta)$-framework, as desired.
\end{proof}
\subsection{The Proof of Lemma \ref{2.1}}
We define a \emph{directed edge} in a $k$-graph to be a $k$-tuple whose vertices correspond to an underlying edge.
Note that  the directed edges $(a,b,c), (b,c,a)$ corresponds to the same underlying edge $\{a,b,c\}$.
Given a $k$-graph system $\emph{\textbf{H}}=\{H_i\}_{i\in[n]}$ on vertex set $V$, we consider the hypergraph $H$ with vertex set $[n]\cup V$ and edge set $\{\{i\}\cup e: e\in E(H_i), i\in[n]\}$.
Define a directed edge to be a $(1,k)$-tuple $(i,v_1,\ldots,v_k)$ with $k$ points corresponding to an underlying edge $\{v_1,\ldots,v_k\}$ in $H_i$.
Given a $k$-tuple $\overrightarrow{S}:=(v_1,\ldots,v_k)$, abbreviated as $v_1\cdots v_k$,  we use $\overrightarrow{S}\subseteq V$ to mean that the corresponding $k$-set of $\overrightarrow{S}$ is a subset of $V$.
Similarly, given a family $F$ of $k$-sets and a $k$-tuple $\overrightarrow{S}$, we use $\overrightarrow{S}\in F$ to denote that the corresponding $k$-set of $\overrightarrow{S}$ is an element of $F$.
Let $\overrightarrow{S}=(v_1,\ldots,v_k)$, $\overrightarrow{S}\setminus\{v_i\}$ is the $(k-1)$-tuple $(v_1,\ldots,v_{i-1},v_{i+1},\ldots,v_k)$ for $i\in[k]$,
$\{v_i'\}\cup\overrightarrow{S}\setminus\{v_i\}$ is the $k$-tuple $(v_1,\ldots,v_{i-1},v_i',v_{i+1},\ldots,v_k)$.
\begin{definition}[Strong Connectivity]
A hypergraph is called strongly connected, if every two directed edges lie on a sequentially walk.
\end{definition}
\begin{claim}
\label{2.5}
If $G$ is a tightly connected graph, then $G$ is strongly connected.
\end{claim}
\begin{proof}
Let $ab$ be a switcher in $G$, by Definition \ref{switcher}, we obtain that $a$ and $b$ share a neighbor $c$.
If we can prove that $(a,b)$ and $(b,a)$ are on a walk $W$, then we can obtain that $G$ is strongly connected.
Since we consider any two directed edges $D_1$ and $D_2$ of $G$, there are walks $W_1$ and $W_2$ starting from $D_1$ and $D_2$ respectively and ending with $\{a,b\}$, $W_1WW_2$ is a tight walk starting from $D_1$ and ending with $D_2$.
While it is easy to see that $aba$ is a tight walk from $(a,b)$ to $(b,a)$, as desired.
\end{proof}
Next, we want to show that switchers can control the length of sequentially walks.
\begin{proposition}
\label{2.6}
If $G$ is a tightly connected graph containing a switcher, then $G$ has a closed tight walk of odd length.
\end{proposition}
\begin{proposition}
\label{2.7}
Let $R$ be a $(1,k)$-graph with a subgraph $H$ which is generated by $\mathcal{C}_i$.
Suppose that $\mathcal{C}_i$ satisfies the conditions of Lemma \ref{2.1}, for any $(1,k-2)$-tuple $\overrightarrow{S}\in \partial_{k-2}(H)$ and two directed edges $D_1,D_2\in C_{\overrightarrow{S}}$, there exists a sequentially walk $W$ of length 0 mod $k$ in $H$ starting from $\overrightarrow{S}D_1$ and ending with $\overrightarrow{S}D_2$.
\end{proposition}
\begin{proof}
Let $\mathcal{C}_i=\{C_{\overrightarrow{S}}:\overrightarrow{S}\in\partial_{k-2}(R){\rm\ and}\ i\in \overrightarrow{S}\}$ and $\overrightarrow{S}=\{i\}\cup \overrightarrow{S}'$ where $\overrightarrow{S}'$ is a $(k-2)$-tuple. By Proposition \ref{2.6}, there is a closed tight walk $W_1$ of odd length in $C_{\overrightarrow{S}}$. By Proposition \ref{2.5}, there is a tight walk $W_2$ starting from $D_1$, ending with $D_2$ and containing $W_1$ as a subwalk.
Let $\ell(W_2)=p$.
We obtain $W_3$ from $W_2$ by replacing $W_1$ with the concatenation of $p+1$ mod 2 copies of $W_1$.
Hence, $W_3$ is a tight walk of even length in $C_{\overrightarrow{S}}$ starting from $D_1$ and ending with $D_2$.

Suppose that $W_3=(a_1,a_2\ldots,a_{2m})$, we have $D_1=(a_1,a_2)$ and $D_2=(a_{2m-1},a_{2m})$.
Note that $(i\ldots i,\overrightarrow{S}'a_1a_2\overrightarrow{S}'a_3a_4\cdots \overrightarrow{S}'a_{2m-1}a_{2m})$ is a sequentially walk in $H$. Moreover, it has length 0 mod $k$, as desired.
\end{proof}
\begin{proposition}
\label{2.8}
Let $R$ be a $(1,k)$-graph with a subgraph $H$ that is generated by $\mathcal{C}_i$.
Suppose $\mathcal{C}_i$ satisfies the conditions of Lemma \ref{2.1}, we consider directed edges $\overrightarrow{S},\overrightarrow{T}\in\partial_{k-2}(H)$ and $D_1\in C_{\overrightarrow{S}}$, $D_2\in C_{\overrightarrow{T}}$.
If $\overrightarrow{S}$ and $\overrightarrow{T}$ differ in exactly one coordinate, then there is sequentially walk of length 0 mod $k$ in $H$ starting from $\overrightarrow{S}D_1$ and ending with $\overrightarrow{T}D_2$.
\end{proposition}
\begin{proof}
Let $\overrightarrow{S}=(i, v_1\ldots v_i\ldots v_{k-2})$ and $\overrightarrow{T}=(i,v_1\ldots u_i\ldots v_{k-2})$ where $u_i\neq v_i$.
By Definition \ref{1.15}, there is a directed edge $D_3$ in $C_{\overrightarrow{S}}\cap C_{\overrightarrow{T}}$, thus $(ii,\overrightarrow{S}\setminus \{i\}D_3\overrightarrow{T}\setminus \{i\})$ is a sequentially walk in $H$.
By Proposition \ref{2.7}, there is a sequentially walk $W_1$ of length 0 mod $k$ starting from $\overrightarrow{S}D_1$ and ending with $\overrightarrow{S}D_3$, $W_2$ of length 0 mod $k$ starting from $\overrightarrow{T}D_3$ and ending with $\overrightarrow{T}D_2$, $(C(W_1)C(W_2),I(W_1)I(W_2))$ is the desired walk.

\end{proof}
\begin{proposition}
\label{2.9}
Let $R$ be a $(1,k)$-graph with a subgraph $H$ that is generated by $\mathcal{C}_i$.
Suppose $\mathcal{C}_i$ satisfies the conditions of Lemma \ref{2.1}, we consider directed edges $\overrightarrow{S},\overrightarrow{T}\in \partial_{k-2}(H)$ and $D_1\in C_{\overrightarrow{S}}$, $D_2\in C_{\overrightarrow{T}}$.
There is a sequentially walk of length 0 mod $k$ in $H$ starting from $\overrightarrow{S}D_1$ and ending with $\overrightarrow{T}D_2$.
\end{proposition}
\begin{proof}
Let $r\in[k-2]$ be the number of indices where $\overrightarrow{S}$ and $\overrightarrow{T}$ differ.
If $r=1$, the result follows from Proposition \ref{2.8}.
Suppose the result is known for $r-1$.
By Definition \ref{1.15}, there exists an edge $pq$ in $C_{\overrightarrow{S}}\cap C_{\overrightarrow{T}}$.

Suppose that the $i$th coordinate vertex of $\overrightarrow{S}$ and $\overrightarrow{T}$ are different, which are replaced with $p$, we obtain $\overrightarrow{S}'$ and $\overrightarrow{T}'$.
Note that $\overrightarrow{S}',\overrightarrow{T}'\in \partial_{k-2}(H)$.
Choose $D_1'\in C_{\overrightarrow{S}'}$.
By Proposition \ref{2.8}, there is a sequentially walk $W_1$ of length 0 mod $k$ from $\overrightarrow{S}D_1$ to $\overrightarrow{S}'D_1'$, similarly, there is a sequentially walk $W_3$ of length 0 mod $k$ from $\overrightarrow{T}'D_2'$ to $\overrightarrow{T}D_2$ where $D_2'\in C_{\overrightarrow{T}'}$.
By induction, there is a sequentially walk $W_2$ from $\overrightarrow{S}'D_1'$ to $\overrightarrow{T}'D_2'$ of length 0 mod $k$.
Thus, $(C(W_1)C(W_2)C(W_3),I(W_1)I(W_2)I(W_3))$ is the desired walk.
\end{proof}
\begin{proof}[The proof of Lemma \ref{2.1}]
Consider any two edges $X$ and $Y$ of $H$.
Since $H$ is generated by $\mathcal{C}_i$, let $X=S\cup A$ and $Y=T\cup B$ where $A\in C_S$ and $B\in C_T$.
The desired walk can be obtained from Proposition \ref{2.9}.

Next, we need to show that $H$ contains a closed walk of length 1 mod $k$.
Since $\mathcal{C}_i$ admits an arc $\{i,v_1,\ldots,v_{k+1}\}$, by Proposition \ref{2.9}, there is a sequentially walk $W$ of length 0 mod $k$ from $\{i,v_2,\ldots,v_{k+1}\}$ to $\{i,v_1,\ldots,v_{k}\}$.
Thus, $(C(W)i, I(W)v_{k+1})$ is a closed walk of length 1 mod $k$.
\end{proof}
\subsection{The proof of Lemma \ref{2.2}}
In this section, we show the details of proving Lemma \ref{2.2}.
The following claim can be seen in \cite{Lang}, we use a corollary of the claim in this paper.
\begin{claim}\cite{Lang}
\label{bf}
Let $H$ be a $k$-graph and $\textbf{b}:V(H)\rightarrow[0,1]$.
Suppose that there exists $m\leq\sum_{v\in V(H)}\textbf{b}(v)/k$ such that for every $v\in V(H)$, the link graph $L_H(\{v\})$ has a $\textbf{b}$-fractional matching of size $m$, then $H$ has a $\textbf{b}$-fractional matching of size $m$.
\end{claim}
\begin{corollary}\label{isolated}
Let $H$ be a $k$-graph, $\alpha\in[0,1)$ and $\textbf{b}:V(H)\rightarrow[0,1]$.
Suppose that there exists $m\leq\sum_{v\in V(H)}\textbf{b}(v)/k$ such that for all but at most $\alpha|V(H)|$ isolated vertices $v$, the link graph $L_H(\{v\})$ has a $\textbf{b}$-fractional matching of size $m$, then $H$ has a $\textbf{b}$-fractional matching of size $m$.
\end{corollary}
\begin{proof}
We first delete the isolated vertices of $H$ and obtain a subgraph $H'$ of $H$.
Thus, $L_{H'}(\{v\})$ has a $\textbf{b}$-fractional matching of size $m$.
By Claim \ref{isolated}, we obtain that $H'$ has a $\textbf{b}$-fractional matching $\textbf{w}$ of size $m$.
Assign a weight $\textbf{b}'(u)\in[0,1]$ to each isolated vertex $u$ of $H$, and $\textbf{b}'(v)=\textbf{b}(v)$ for each non-isolated vertex $v$ of $H$, it is obvious that $H$ has a $\textbf{b}'$-fractional matching $\textbf{w}$ of size $m$ since $\sum_{e\ni u}\textbf{w}(e)=0$ for any isolated vertex $u$ and $E(H')=E(H)$.
\end{proof}
\begin{proposition}
\label{2.11}
Let $R$ be a $(1,k)$-graph on $[n/k]\cup V$ where $|V|=n$, $\gamma>0$, $\alpha\in[0,1)$, $\textbf{b}:[n/k]\cup V\rightarrow [1-\gamma,1]$.
Suppose that there exists $m\leq\sum_{v\in V(R)}\textbf{b}(v)/(k+1)$ such that given $c\in[n/k]$, for all but at most $\alpha n$ vertices $v\in V$, the link graph $L_R(\{c,v\})$ has a $\textbf{b}$-fractional matching of size $m$, then $R$ has a $\textbf{b}$-fractional matching of size $m/k$.
\end{proposition}
\begin{proof}
By Corollary \ref{isolated} with $H$ being $L_R(\{c\})$ for $c\in[n/k]$, we obtain that $L_R(\{c\})$ has a $\textbf{b}$-fractional matching of size $m$ for $c\in[n/k]$.

Next, we want to construct a $\textbf{b}$-fractional matching of size $m/k$ for $R$.
Let $\textbf{w}_c:E(L_R(\{c\}))\rightarrow[0,1]$ such that $\sum_{v\in e,e\in L_R(\{c\})}\textbf{w}_c(e)\leq \textbf{b}(v)$ where $\sum_{e\in L_R(\{c\})}\textbf{w}_c(e)=m$.
Let $\textbf{w}(f)=\frac{1}{n}\textbf{w}_c(e)$ for $e\in L_R(\{c\})$ and $f=e\cup \{c\}$, $c\in[n/k]$.
Thus, we have $\sum_{f\in E(R)}\textbf{w}(f)=\sum_{c\in[n/k]}\sum_{e\in L_R(\{c\})}\frac{1}{n}\textbf{w}_c(e)=\frac{m}{k}$.
It is easy to see that
$\sum_{c\in f}\textbf{w}(f)=\sum_{e\in L_R(\{c\})}\frac{1}{n}\textbf{w}_c(e)=\frac{m}{n}\leq\frac{1}{k}\leq\textbf{b}(c)$.
And $\sum_{v\in f}\textbf{w}(f)=\sum_{c\in[n/k]}\sum_{v\in e, e\in L_R(\{c\})}\frac{k}{n}\textbf{w}_c(e)\leq\sum_{c\in[n/k]}\frac{k}{n}\textbf{b}(v)=\textbf{b}(v)$ for $v\in V$.
As desired.
\end{proof}
We use the following results of \cite{Lang} directly.
\begin{proposition}\cite{Lang}
\label{2.12}
Let $H$ be a $k$-graph and $m\leq v(H)/k$.
If for every vertex $v$ of $V(H)$, $L_H(\{v\})$ has a fractional matching of size $m$, then $H$ has a fractional matching of size $m$.
\end{proposition}
\begin{proposition}
\label{2.13}
Let $d\in[k-2]$ and $\alpha,\gamma,\delta>0,k\geq3$ such that $\alpha,\gamma\ll1/k$.
Let $R$ be a $(1,k)$-graph on $[t]\cup V$ with $\alpha$-perturbed minimum $(1,k-2)$-degree $\delta$ where $|V|=t$.
If for every $S\in\partial_{d}(R)$, the link graph $L(S)$ contains a fractional matching of size at least $(1+1/k)(1/(k+1)+\gamma)t$, then for every edge $S'\in \partial_{1}(R)$, the link graph $L(S')$ contains a fractional matching of size at least $(1+1/k)(1/(k+1)+\gamma)t$.
\end{proposition}
\begin{proof}
We prove it by induction on $d$. Note that the base case when $d=1$ is obvious.
Suppose that given $d\in[2,k-2]$, we obtain the conclusion for $d'<d$.
Let $S\subseteq V(R)$ be a $(1,d-1)$-set in $\partial_{d-1}(R)$.
Consider any vertex $s'$ in $\partial_1(L_R(S))$,  $S\cup \{s'\}$ is an edge in $\partial_{d}(R)$.
By assumption, $L_R(S\cup \{s'\})$ has a fractional matching of size at least $(1+1/k)(1/(k+1)+\gamma)t$,
thus, we have $L_{R'}(\{s'\})$ contains a fractional matching of size at least $(1+1/k)(1/(k+1)+\gamma)t$ for any vertex $s'$ of $V$ where $R'$ is the subgraph of $L_R(S)$ induced on the non-isolated vertices of $L_R(S)$.

By Definition \ref{1.16}, $S$ has at most $\alpha t$ neighbors in $\overline{\partial_{d}(R)}$.
It follows that $v(R')=\partial_1(L_R(S))\geq(1-\alpha)t$ and $(1+1/k)(1/(k+1)+\gamma)t\leq v(R')/(k-d+1)$ since $\alpha,\gamma\ll1/k$.
By Proposition \ref{2.12} with the condition that $L_{R'}(\{s'\})$ contains a fractional matching of size at least $(1+1/k)(1/(k+1)+\gamma)t$ for any vertex $s'$ of $V$, we obtain $R'$(and thus $L_R(S)$) contains a fractional matching of size $(1+1/k)(1/(k+1)+\gamma)t$.
Since $S$ is arbitrary, for any $S\in \partial_{d-1}(R)$, $L_R(S)$ contains a fractional matching of size $(1+1/k)(1/(k+1)+\gamma)t$.
Hence, we are done by the induction hypothesis.

\end{proof}
\begin{proof}[The proof of Lemma \ref{2.2}]
Suppose that $V(H)=[t/k]\cup V'$ where $|V'|=t$.
By assumption, $C_S$ contains a fractional matching of size $(1+1/k)(1/(k+1)+\gamma)t$ for every $S\in \partial_{k-2}(H)$ and $C_S$ is a subgraph of $L_H(S)$.
By Proposition \ref{2.13}, we have $L_H(\{i,v\})$ contains a fractional matching of size $(1+1/k)(1/(k+1)+\gamma)t$ for every $\{i,v\}\in\partial_1(H)$.

We want to show that $H$ is $\gamma$-robustly matchable.
Given a vertex weight $\textbf{b}$: $[t/k]\cup V'\rightarrow[1-\gamma,1]$, we have to find a $\textbf{b}$-fractional matching $\textbf{w}$ such that $\sum_{e\ni v}\textbf{w}(e)=\textbf{b}(v)/k$ for any vertex $v\in V(H)$.
That is, we need to find a $\textbf{b}$-fractional matching with size $\sum_{v\in V(H)}\textbf{b}(v)/k(k+1)$.
Given $i\in[t/k]$, there are at most $\alpha t$ isolated $(1,1)$-tuples by Definition \ref{1.16}.
For any non-isolated $(1,1)$-tuple $(i,v)$ of $V(H)$, let $\textbf{x}$ be a fractional matching in $L_H(\{i,v\})$ of size at least $(1+1/k)(1/(k+1)+\gamma)t$ and let $\textbf{w}'=(1-\gamma)\textbf{x}$, since $1-\gamma\leq \textbf{b}(v)$ for any $v\in V(H)$, thus $\textbf{w}'$ is a $\textbf{b}$-fractional matching in $L_H(\{i,v\})$.
Moreover $\textbf{w}'$ has size at least $(1-\gamma)(1+1/k)(1/(k+1)+\gamma)t\geq (1+1/k)t/(k+1)\geq\sum_{v\in V(H)}\textbf{b}(v)/(k+1)$ since $1/t\ll\gamma\ll1/k$.
We can assume that $\textbf{w}'$ has size exactly $\sum_{v\in V(H)}\textbf{b}(v)/(k+1)$.
By Proposition \ref{2.11}, we obtain that $H$ has a $\textbf{b}$-fractional matching of size $\sum_{v\in V(H)}\textbf{b}(v)/k(k+1)$, as desired.
\end{proof}
\subsection{The Proof of Lemma \ref{2.3}}
We use the following claim directly, which can be seen in \cite{Lang}.
\begin{claim}\label{lem2.3}\cite{Lang}
Let $t,d,k$ be integers with $d\in[k-1]$ and $\delta,\varepsilon,\alpha>0$ with $1/t\ll\varepsilon\ll\alpha\leq\delta,1/k$.
Let $R$ be a $k$-graph on $t$ vertices with minimum relative $d$-degree $\overline{\delta_d(R)}\geq\delta$.
Let $I$ be a subgraph of $R$ of edge density at most $\varepsilon$.
Then there exists a vertex spanning subgraph $R'\subseteq R-I$ od $\alpha$-perturbed minimum relative $d$-degree at least $\delta-\alpha$.
\end{claim}

The $(1,k)$-graph $R_i$ on $\{i\}\cup V$ with minimum relative $(1,k-2)$-degree at least $\delta$ is equivalent to a $k$-graph $R_i'$ on $V$ with minimum relative $(k-2)$-degree at least $\delta$.
Thus, by Claim~\ref{lem2.3}, we obtain Lemma~\ref{2.3}.
\section{Obtaining vicinity}
In this section, we determine the $(k-2)$-vicinity threshold of $(1,k)$-graphs.
Lov\'{a}sz's formulation of the Kruskal-Katona theorem states
that, for any $x>0$, if $G$ is a $k$-graph with $e(G)\geq\binom{x}{k}$ edges, then $e_j(G)\geq\binom{x}{j}$ for every $j\in[k]$ (Theorem 2.14 in \cite{frankl2018extremal}).
By approximating the binomial coefficients, they~\cite{Lang} deduce the following variant.
\begin{lemma}[Kruskal-Katona theorem]
\label{3.1}\cite{Lang}
Let $1/t\ll\varepsilon\ll1/k$ and $G$ be a graph on $t$ vertices and edge density $\delta$, then $\partial(G)$ has at least $(\delta^{1/2}-\varepsilon)t$ vertices.
\end{lemma}
\begin{proposition}
\label{3.6}
Let $t\in N$ and $\gamma,\delta',\delta>0$ with $1/t\ll\varepsilon\ll\delta$ and $\delta+\delta^{1/2}>1+\varepsilon$.
Let $R_i$ be a $(1,k)$-graph on $\{i\}\cup V$ where $|V|=t$ with a subgraph that is generated by a $(k-2)$-vicinity $\mathcal{C}_i$.
Suppose that each $C_S\in \mathcal{C}_i$ has edge density at least $\delta+\mu$, then $\mathcal{C}_i$ admits an arc.
\end{proposition}
\begin{proof}
Consider an arbitrary set $S=\{i,v_1,\ldots,v_{k-2}\}\in\partial_{k-2}(R_i)$.
By averaging, there is a vertex $v_{k-1}$ with relative vertex degree at least $\delta$ in $C_S$.
Set $S'=\{i,v_2,\ldots,v_{k-1}\}$, we have $S'\in\partial_{k-2}(R_i)$.
Thus, $C_{S'-\{v_1\}}$ has edge density at least $\delta+\mu/2$.
By Lemma \ref{3.1}, $\partial(C_{S'}-\{v_1\})$ has at least $(\delta^{1/2}-\varepsilon)t$ vertices.

By the choice of $v_{k-1}$ and the pigeonhole principle, $\partial(C_{S'}-\{v_1\})$ and $L(\{i,v_1,\ldots,v_{k-1}\})$ must share a common vertex $v_k$.
Since $v_k\in\partial(C_{S'}-\{v_1\})$, there is another vertex $v_{k+1}$ such that $\{v_k,v_{k+1}\}\in C_{S'}-\{v_1\}$.
Thus, $\{i,v_1,\ldots,v_{k+1}\}$ is an arc.
\end{proof}

We use the following result of \cite{Lang}.
\begin{lemma}\cite{Lang}
\label{3.7}
Let $1/t\ll\gamma\ll\mu$, suppose that $L_1$ and $L_2$ are graphs on a common vertex set of size $t$ such that $L_1$, $L_2$ has edge density at least $5/9+\mu$.
For $i\in[2]$, let $C_i$ be a tight component of $L_i$ with a maximum number of edges.
We have
\begin{itemize}
  \item[(i)] $C_1$ and $C_2$ has an edge in common,
  \item[(ii)] $C_i$ has a switcher for $i\in[2]$,
  \item[(iii)] $C_i$ has a fractional matching of density $1/3+\gamma$ for $i\in[2]$,
  \item[(iv)] $C_i$ has edge density at least $4/9+\gamma$ for $i\in[2]$.
\end{itemize}
\end{lemma}
\begin{proof}[The proof of Lemma \ref{1.20}]
Let $\alpha,\gamma,\mu>0$ with
\[
1/t\ll\alpha\ll\delta\ll\mu\ll5/9.
\]
Consider a $(1,k)$-graph $R$ on $[t]\cup V$ where $|V|=t$ and each $R_i:=R[\{i\}\cup V]$ has $\alpha$-perturbed minimum relative $(1,k-2)$-degree at least $5/9+\mu$.
For every $S\in\partial_{k-2}(R)$, let $C_S$ be a tight component of $L(S)$ with a maximum number of edges and $\mathcal{C}_i=\{C_S:S\in \partial_{k-2}(R){\rm\ and}\ i\in S\}$.
By the choice of $C_S$, (V1) holds obviously.
By Lemma \ref{3.7}, $\mathcal{C}_i$ satisfies (V2), (V4), (V5) and (V6).
Every $C_S\in\mathcal{C}_i$ contains a switcher.
By Proposition \ref{3.6}, $\mathcal{C}_i$ contains an arc since $4/9+(4/9)^{1-1/2}=1+1/9$, thus $\mathcal{C}=\{\mathcal{C}_i:i\in[t]\}$ satisfies (V3), as desired.
\end{proof}

\section{Tools}
\subsection{Regular Complexes}
A hypergraph $H=(V,E)$ is a \emph{complex} if its edge set is down-closed, meaning that whenever $e\in E$ and $e'\subseteq e$, we have $e'\in E$.
A $k$-complex is a complex where all edges have size at most $k$.
Given a complex $H$, we use $H^{(i)}$ to denote the $i$-graph obtained by taking all vertices of $H$ and edges of size $i$.
Denote the number of edges of size $i$ in $H$ by $e_i(H)$.

Let $\mathcal{P}$ partition a vertex set $V$ into parts $V_1,\ldots,V_s$.
Then we say that a subset $S\subseteq V$ is $\mathcal{P}$-\emph{partite} if $|S\cap V_i|\leq1$ for every $i\in[s]$.
Similarly, we say that hypergraph $\mathcal{H}$ is $\mathcal{P}$-\emph{partite} if all of its edges are $\mathcal{P}$-partite.
In this case we refer to the parts of $\mathcal{P}$ as the \emph{vertex class} of $\mathcal{H}$.
We say that a hypergraph $\mathcal{H}$ is $s$-\emph{partite} if there is some partition $\mathcal{P}$ of $V(\mathcal{H})$ into $s$ parts for which $\mathcal{H}$ is $\mathcal{P}$-partite.

Let $\mathcal{H}$ be a $\mathcal{P}$-partite complex.
Then for any $A\subseteq[s]$ we write $V_A$ for $\bigcup_{i\in A}V_i$.
The \emph{index} of a $\mathcal{P}$-partite set $S\subseteq V$ is $i(S):=\{i\in[s]:|S\cap V_i|=1\}$.
We write $\mathcal{H}_A$ to denote the collection of edges in $\mathcal{H}$ with index $A$, that is, $\mathcal{H}_A$ can be regarded as an $|A|$-partite $|A|$-graph on vertex set $V_A$.\
Similarly, if $X$ is a $j$-set of indexes of vertex classes of $\mathcal{H}$ we write $\mathcal{H}_X$ for the $j$-partite $j$-uniform subgraph of $\mathcal{H}^{(j)}$ induced by $\bigcup_{i\in X} V_i$.
We write $\mathcal{H}_{X<}$ for the $j$-partite hypergraph with vertex set $\bigcup_{i\in V_X} V_i$ and edge set $\bigcup_{X'\subset X}\mathcal{H}_{X'}$.

Let $H_i$ be any $i$-partite $i$-graph and $H_{i-1}$ be any $i$-partite $(i-1)$-graph on a common vertex set $V$ partitioned into $i$ common vertex classes.
Denote $K_{i}(H_{i-1})$ by the $i$-partite $i$-graph on $V$ whose edges are all $i$-sets which are supported on $H_{i-1}$(i.e. induce a copy of complete $(i-1)$-graph $K_i^{i-1}$ on $i$ vertices in $H_{i-1}$).
The \emph{density of} $H_i$ \emph{with respect to} $H_{i-1}$ is defined to be
\[
d(H_i|H_{i-1}):=\frac{|K_i(H_{i-1})\cap H_i|}{|K_i(H_{i-1})|}
\]
if $|K_i(H_{i-1})|>0$.
For convenience, we take $d(H_i|H_{i-1}):=0$ if $|K_i(H_{i-1})|=0$.
When $H_{i-1}$ is clear from the context, we simply refer $d(H_i|H_{i-1})$ as the \emph{relative density of} $H_i$.
More generally, if $\textbf{Q}:=(Q_1,\ldots,Q_r)$ is a collection of $r$ not necessarily disjoint subgraphs of $H_{i-1}$, we define
\[
K_i(\textbf{Q}):=\bigcup_{j=1}^rK_i(Q_j)
\]
and
\[
d(H_i|\textbf{Q}):=\frac{|K_i(\textbf{Q})\cap H_i|}{|K_i(\textbf{Q})|}
\]
if $|K_i(\textbf{Q})|>0$.
Similarly, we take $d(H_i|\textbf{Q}):=0$ if $|K_i(\textbf{Q})|=0$.
We say that $H_i$ is $(d_i,\varepsilon,r)$-\emph{regular with respect to} $H_{i-1}$ if we have $d(H_i|\textbf{Q})=d_i\pm\varepsilon$ for every $r$-set $\textbf{Q}$ of subgraphs of $H_{i-1}$ such that $|K_i(\textbf{Q})|>\varepsilon|K_i(H_{i-1})|$.
We refer to $(d_i,\varepsilon,1)$-regularity simply as $(d_i,\varepsilon)$-\emph{regularity}.
We say that $H_i$ is $(\varepsilon,r)$-regular with respect to $H_{i-1}$ to mean that there exists some $d_i$ for which $H_i$ is $(d_i,\varepsilon,r)$-regular with respect to $H_{i-1}$.
Given an $i$-graph $G$ whose vertex set contains that of $H_{i-1}$, we say that $G$ is $(d_i,\varepsilon,r)$-\emph{regular with respect to} $H_{i-1}$ if the $i$-partite subgraph of $G$ induced by the vertex classes of $H_{i-1}$ is $(d_i,\varepsilon,r)$-regular with respect to $H_{i-1}$.
Similarly, when $H_{i-1}$ is clear from the context, we refer to the relative density of this $i$-partite subgraph of $G$ with respect to $H_{i-1}$ as the \emph{relative density of }$G$.

Now let $\mathcal{H}$ be an $s$-partite $k$-complex on vertex classes $V_1,\ldots,V_s$, where $s\geq k\geq3$.
Since $\mathcal{H}$ is a complex, if $e\in \mathcal{H}^{(i)}$ for some $i\in[2,k]$, then the vertices of $e$ induce a copy of $K_i^{i-1}$ in $\mathcal{H}^{(i-1)}$.
This means that for any index $A\in\binom{[s]}{i}$, the density $d(\mathcal{H}^{(i)}[V_A]|\mathcal{H}^{(i-1)}[V_A])$ can be regarded as the proportion of `possible edges' of $\mathcal{H}^{(i)}[V_A]$ which are indeed edges.
We say that $\mathcal{H}$ is $(d_2,\ldots,d_k,\varepsilon_k,\varepsilon,r)$-\emph{regular} if
\begin{enumerate}
  \item for $i\in[2,k-1]$ and $A\in\binom{[s]}{i}$, the induced subgraph $\mathcal{H}^{(i)}[V_A]$ is $(d_i,\varepsilon)$-regular with respect to $\mathcal{H}^{(i-1)}[V_A]$ and
  \item for any $A\in\binom{[s]}{k}$, the induced subgraph $\mathcal{H}^{(k)}[V_A]$ is $(d_k,\varepsilon_k,r)$-regular with respect to $\mathcal{H}^{(k-1)}[V_A]$.
\end{enumerate}

\subsection{Regular Slices}
The Regular Slice Lemma says that any $k$-graph $G$ admits a regular slice.
Informally speaking, a regular slice of $G$ is a partite $(k-1)$-complex $\mathcal{J}$ whose vertex classes have equal size, whose subgraphs $\mathcal{J}^{(2)},\ldots,\mathcal{J}^{(k-1)}$ satisfy certain regularity properties and which moreover has the property that $G$ is regular with respect to $\mathcal{J}^{(k-1)}$.
The first two of these conditions are formalised in the following definition: we say that a $(k-1)$-complex $\mathcal{J}$ is $(t_0,t_1,\varepsilon)$-\emph{equitable}, if it has the following properties.
\begin{enumerate}
  \item $\mathcal{J}$ is $\mathcal{P}$-partite for a $\mathcal{P}$ which partitions $V(\mathcal{J})$ into $t$ parts of equal size, where $t_0\leq t\leq t_1$.
      We refer to $\mathcal{P}$ as the \emph{ground partition} of $\mathcal{J}$, and to the parts of $\mathcal{P}$ as the \emph{clusters} of $\mathcal{J}$.
  \item There exists a \emph{density vector} $\textbf{d}=(d_2,\ldots,d_{k-1})$ such that for $i\in[2,k-1]$ we have $d_i\geq1/t_1$ and $1/d_i\in \mathbb{N}$ and for each $A\subseteq\mathcal{P}$ of size $i$, the $i$-graph $\mathcal{J}^{(i)}[V_A]$ induced on $V_A$ is $(d_i,\varepsilon)$-regular with respect to $\mathcal{J}^{(i-1)}[V_A]$.
\end{enumerate}
If $\mathcal{J}$ has density vector $\textbf{d}=(d_2,\ldots,d_{k-1})$, then we will say that $\mathcal{J}$ is $(d_2,\ldots,d_{k-1},\varepsilon)$-regular, or $(\textbf{d},\varepsilon)$-\emph{regular}, for short.
For any $k$-set $X$ of clusters of $\mathcal{J}$, we write $\hat{\mathcal{J}}_X$ for the $k$-partite $(k-1)$-graph $\mathcal{J}_{X<}^{(k-1)}$.
Given a $(t_0,t_1,\varepsilon)$-equitable $(k-1)$-complex $\mathcal{J}$, a $k$-set $X$ of clusters of $\mathcal{J}$ and a $k$-graph $G$ on $V(\mathcal{J})$, we say that $G$ is $(d,\varepsilon_k,r)$-\emph{regular with respect to} $X$ if $G$ is $(d,\varepsilon_k,r)$-regular with respect to $\hat{\mathcal{J}}_X$.
We will also say that $G$ is $(\varepsilon_k,r)$-\emph{regular with respect to} $X$ if there exists a $d$ such that $G$ is $(d,\varepsilon_k,r)$-regular with respect to $X$.
We write $d_{\mathcal{J},G}^*(X)$ for the relative density of $G$ with respect to $\hat{\mathcal{J}}_X$, or simply $d^*(X)$ if $\mathcal{J}$ and $G$ are clear from the context, which will always be the case in applications.

We now give the key definition of the Regular Slice Lemma.
\begin{definition}[Regular Slice]
Given $\varepsilon,\varepsilon_k>0$, $r,t_0,t_1\in \mathbb{N}$, a $k$-graph $G$ and a $(k-1)$-complex $\mathcal{J}$ on $V(G)$, we call $\mathcal{J}$ a $(t_0,t_1,\varepsilon,\varepsilon_k,r)$-regular slice for $G$ if $\mathcal{J}$ is $(t_0,t_1,\varepsilon)$-equitable and $G$ is $(\varepsilon_k,r)$-regular with respect to all but at most $\varepsilon_k\binom{t}{k}$ of the $k$-sets of clusters of $\mathcal{J}$, where $t$ is the number of clusters of $\mathcal{J}$.
\end{definition}

It will sometimes be convenient not to specify all parameters, we may write that $\mathcal{J}$ is $(\cdot,\cdot,\varepsilon)$-equitable or is a $(\cdot,\cdot,\varepsilon,\varepsilon_k,r)$-slice for $G$, if we do not wish to specify $t_0$ and $t_1$.

Given a regular slice $\mathcal{J}$ for a $k$-graph $G$, it will be important to know the relative densities $d^*(X)$ for $k$-sets $X$ of clusters of $\mathcal{J}$.
To keep track of these we make the following definition.
\begin{definition}[Weighted reduced $k$-graph]
Let $G$ be a $(1,k)$-graph and let $\mathcal{J}$ be a $(t_0,t_1,\varepsilon,\varepsilon_{k+1},r)$-regular slice for $G$.
We define the \emph{weighted reduced} $(1,k)$-\emph{graph} of $G$, denoted by $R(G)$, to be the complete weighted $(1,k)$-graph whose vertices are the clusters of $\mathcal{J}$ and where each edge $X$ is given weight $d^*(X)$.

Similarly, for $d_{k+1}>0$, we define the $d_{k+1}$-reduced $(1,k)$-graph $R_{d_{k+1}}(G)$ to be the (unweighted) $(1,k)$-graph whose vertices are the clusters of $\mathcal{J}$ and whose edges are all $(1,k)$-sets $X$ of clusters of $\mathcal{J}$ such that $G$ is $(\varepsilon_{k+1},r)$-regular with respect to $X$ and $d^*(X)\geq d_{k+1}$.
\end{definition}

Given a $(1,k)$-graph $G$ on $[n]\cup V$, a vertex $v\in V$ and a color $c\in[n]$, recall that $\deg_G(c,v)$ is the number of edges of $G$ containing $c$ and $v$ and $\overline{\deg}_G(c,v)=\deg_G(c,v)/\binom{n-1}{k-1}$ is the relative degree of $v$ in $G$.
Given a $(t_0,t_1,\varepsilon)$-equitable $(k-1)$-complex $\mathcal{J}$ with $V(\mathcal{J})\subseteq V(G)$, the \emph{rooted degree} of $(c,v)$ \emph{supported by} $\mathcal{J}$, written by $\deg_G((c,v),\mathcal{J})$, is defined as the number of $(k-1)$-sets $T$ in $\mathcal{J}^{(k-1)}$ such that $T\cup \{c,v\}$ forms an edge in $G$.
Then the relative degree $\overline{\deg}_G((c,v);\mathcal{J})$ of $(c,v)$ in $G$ supported by $\mathcal{J}$ is defined as $\overline{\deg}_G((c,v);\mathcal{J})=\deg_G((c,v);\mathcal{J})/e(\mathcal{J}^{(k-1)})$.
\begin{definition}[Representative rooted degree]
Let $\eta>0$, $G$ be a $(1,k)$-graph on $[n]\cup V$ and $\mathcal{J}$ be a $(t_0,t_1,\varepsilon,\varepsilon_{k+1})$-regular slice for $G$.
We say that $\mathcal{J}$ is $\eta$-rooted-degree-representative if for any vertex $v\in V$ and any color $c\in[n]$, we have
\[
|\overline{\deg}_G((c,v);\mathcal{J})-\overline{\deg}_G(c,v)|<\eta.
\]
\end{definition}
\begin{definition}[Regular Setup]
Let $k,m,r,t\in \mathbb{N}$ and $\varepsilon,\varepsilon_{k+1},d_2,\ldots,d_{k+1}>0$.
We say that $(G,G_{\mathcal{J}},\mathcal{J},\mathcal{P},R)$ is a $(k,m,t,\varepsilon,\varepsilon_{k+1},r,d_2,\ldots,d_{k+1})$-regular setup, if
\begin{itemize}
  \item[(RS1)] $G$ is a $(1,k)$-graph on $[n]\cup V$ where $|V|=n$ and $G_{\mathcal{J}}\subseteq G$,
  \item[(RS2)] $\mathcal{J}$ is a $(\cdot,\cdot,\varepsilon,\varepsilon_{k+1},r)$-regular slice for $G$ with density vector $\textbf{d}=(d_2,\ldots,d_{k})$,
  \item[(RS3)] $\mathcal{P}$ is the ground partition of $\mathcal{J}$ with initial partition of $[n]\cup V$ and $2t$ clusters, each of size $m$,
  \item[(RS4)] $R$ is a subgraph of $R_{d_{k+1}}(G)$,
  \item[(RS5)] for each $X\in E(R)$, $G_{\mathcal{J}}$ is $(d_{k+1},\varepsilon_{k+1},r)$-regular with respect to $X$.
\end{itemize}
We further say that $(G,G_{\mathcal{J}},\mathcal{J},\mathcal{P},R)$ is representative if
\begin{itemize}
  \item[(RS6)] $\mathcal{J}$ is $\varepsilon_{k+1}$-rooted-degree-representative.
\end{itemize}
\end{definition}

The Regular Slice Lemma of ~\cite{Peter2017Tight} ensures that every sufficiently large $k$-graph has a representative regular slice.
Given the existence of a regular slice, it is easy to derive the existence of a regular setup.
In ~\cite{Lang}, it is stated directly in terms of regular setups.
And it is an easy corollary of giving a sufficiently large $(1,k)$-graph.
\begin{lemma}[Regular Setup Lemma~\cite{Peter2017Tight}]
\label{5.5}
Let $k,t_0$ be positive integers, $\delta,\mu,\alpha,\varepsilon_{k+1},d_{k+1}$ be positive and $r:\mathbb{N}\rightarrow\mathbb{N}$ and $\varepsilon: \mathbb{N}\rightarrow(0,1]$ be functions.
Suppose that
\[
k\geq3,\varepsilon_{k+1}\ll\alpha,d_{k+1}
\ll\mu.
\]
Then there exists $t_1$ and $m_0$ such that the following holds for all $n\geq2t_1m_0$.
Let $G$ be a $(1,k)$-graph on $[n]\cup V$ where $|V|=n$ and suppose that $G$ has minimum relative $(1,k-2)$-degree $\overline{\delta}_{1,k-2}(G)\geq\delta+\mu$.
There exists $\textbf{d}=(d_2,\ldots,d_{k+1})$ and a representative $(k,m,2t,\varepsilon(t_1),\varepsilon_{k+1},r(t_1),\textbf{d})$-regular setup $(G,G_{\mathcal{J}},\mathcal{J},\mathcal{P},R_{d_{k+1}})$ with $t\in[t_0,t_1]$, $m_0\leq m$ and $n\leq(1+\alpha)mt$.
Moreover, there is a $(1,k)$-graph $I$ on $\mathcal{P}$ of edge density at most $\varepsilon_{k+1}$ such that $R=R_{d_{k+1}}\cup I$ has minimum relative $(1,k-2)$-degree at least $\delta+\mu/2$.
\end{lemma}
\subsection{Tools for working with regularity}
Let $\mathcal{G}$ be a $\mathcal{P}$-partite $k$-complex and $X_1,\ldots,X_s\in \mathcal{P}$(possibly with repetition), and let $\mathcal{H}$ be a $k$-complex on vertices $[s]$.
We say that an embedding of $\mathcal{H}$ in $\mathcal{G}$ is \emph{partition}-\emph{respecting}, if $i$ is embedded in $X_i$ for $i\in[s]$.
Note that this notion depends on the labeling of $V(\mathcal{H})$ and the clusters $X_1,\ldots,X_s$, but these will be clear in the paper.
Denote the set of labelled partition-respecting copies of $\mathcal{H}$ in $\mathcal{G}$ by $\mathcal{H}_{\mathcal{G}}[\bigcup_{i\in S} X_i]$.
When $X_1,\ldots,X_s$ are clear, we denote it by $\mathcal{H}_{\mathcal{G}}$ for short.
Recall that $e_i(\mathcal{H})$ denotes the number of edges of size $i$ in $\mathcal{H}$.

The following lemma states that the number of copies of a given small $k$-graph inside a regular slice is roughly what we expect if the edges inside a regular slice were chosen randomly.
There are many different versions in \cite{Peter2017Tight,2006Embeddings,17,40} and we use the following version in \cite{2006Embeddings}.
\begin{lemma}[Counting Lemma~\cite{2006Embeddings}]
\label{counting lemma}
Let $k,s,r,m$ be positive integers and let $\beta,d_2,\ldots,d_k,\varepsilon,\varepsilon_k$ be positive constants such that $1/d_i\in \mathbb{N}$ for $i\in[2,k-1]$ and such that
\[
1/m\ll1/r,\varepsilon\ll\varepsilon_k,d_2,\ldots,d_{k-1},
\]
\[
\varepsilon_k\ll\beta,d_k,1/s.
\]
Let $H$ be a $k$-graph on $[s]$ and let $\mathcal{H}$ be the $k$-complex generated by the down-closure of $H$.
Let $\textbf{d}=(d_2,\cdots,d_k)$, let $(G,G_{\mathcal{J}},\mathcal{J},\mathcal{P},R)$ be a $(k,m,\cdot,\varepsilon,\varepsilon_k,r,\textbf{d})$-regular setup and $\mathcal{G}=\mathcal{J}\cup G_{\mathcal{J}}$.
Suppose $X_1,\ldots,X_s$ are such that $i\mapsto X_i$ is a homomorphism from $H$ into $R$, then the number of labelled partition-respecting copies of $\mathcal{H}$ in $\mathcal{G}$ satisfies
\[
|\mathcal{H}_{\mathcal{G}}|=(1\pm\beta)\left(\prod_{i=2}^kd_i^{e_i(\mathcal{H})}\right)m^s.
\]
\end{lemma}

The following tool allows us to extend small subgraphs into a regular slice.
It was given by Cooley, Fountoulakis, K\"{u}hn and Osthus \cite{2006Embeddings}.
\begin{lemma}[Extension Lemma~\cite{2006Embeddings}]
\label{extension}
Let $k,s,s',r,m$ be positive integers, where $s'<s$ and let $\beta,d_2,\ldots,d_k,\varepsilon,\varepsilon_k$ be positive constants such that $1/d_i\in \mathbb{N}$ for $i\in[2,k-1]$ and such that
\[
1/m\ll1/r,\varepsilon\ll\varepsilon_k,d_2,\ldots,d_{k-1},
\]
\[
\varepsilon_k\ll\beta,d_k,1/s.
\]
Suppose $H$ is a $k$-graph on $[s]$.
Let $\mathcal{H}$ be the $k$-complex generated by the down-closure of $H$ and $\mathcal{H}'$ be an induced subcomplex of $\mathcal{H}$ on $s'$ vertices.
Let $\textbf{d}=(d_2,\ldots,d_k)$ and $(G,G_{\mathcal{J}},\mathcal{J},\mathcal{P},R)$ be a $(k,m,\cdot,\varepsilon,\varepsilon_k,r,\textbf{d})$-regular setup and $\mathcal{G}=\mathcal{J}\cup G_{\mathcal{J}}$.
Suppose $X_1,\ldots,X_s$ are such that $i\mapsto X_i$ is a homomorphism from $H$ into $R$.
Then all but at most $\beta|\mathcal{H}_{\mathcal{G}}'|$ labelled partition-respecting copies of $\mathcal{H}'$ in $\mathcal{G}$ extend to
\[
(1\pm\beta)\left(\prod_{i=2}^kd_i^{e_i(\mathcal{H})-e_i(\mathcal{H}')}\right)m^{s-s'}
\]
labelled partition-respecting copies of $\mathcal{H}$ in $\mathcal{G}$.
\end{lemma}

In some certain situation, we look for structures whose edges lie entirely in the $(k-1)$-complex $\mathcal{J}$ of a regular setup.
We can no longer use the above lemmas whose input is a regular setup rather than an equitable complex.
Also, the above lemmas requires $r$ to be large enough with respect to $\varepsilon_k$ while the $(k-1)$-th level of  $\mathcal{J}$ will only need to be $(d_{k-1},\varepsilon)$-regular with respect to the lower level.
We can use a Dense Counting Lemma as proved by Kohayakawa, R\"{o}dl and Skokan~\cite{2002Hypergraphs}.
We state the following version given by Cooley, Fountoulakis, K\"{u}hn and Osthus~\cite{2006Embeddings}.
\begin{lemma}[Dense Counting Lemma~\cite{2006Embeddings}]
\label{dense counting}
Let $k,s,m$ be positive integers and $\varepsilon,d_2,\ldots,d_{k-1},\beta$ be positive constants such that
\[
1/m\ll\varepsilon\ll\beta\leq d_2,\ldots,d_{k-1},1/s.
\]
Suppose $H$ is a $(k-1)$-graph on $[s]$ and $\mathcal{H}$ is the $(k-1)$-complex generated by the down-closure of $H$.
Let $\textbf{d}=(d_2,\ldots,d_{k-1})$ and $\mathcal{J}$ be a $(\textbf{d},\varepsilon)$-regular $(k-1)$-complex with ground partition $\mathcal{P}$, each size of whose vertex class is $m$.
If $X_1,\ldots,X_s\in \mathcal{P}$, then
\[
|\mathcal{H}_{\mathcal{J}}|=(1\pm\beta)\prod_{i=2}^{k-1}d_i^{e_i(\mathcal{H})}m^s.
\]
\end{lemma}
The following lemma gives the number of edges in each layer of a regular slice.
\begin{lemma}~\cite{Peter2017Tight}
\label{7.7}
Suppose that $1/m\ll\varepsilon\ll\beta\ll d_2,\ldots,d_{k-1},1/k$ and that $\mathcal{J}$ is a $(\cdot,\cdot,\varepsilon)$-equitable $(k-1)$-complex with density vector $(d_2,\ldots,d_{k-1})$ and clusters of size $m$.
Let $X$ be a set of at most $k-1$ clusters of $\mathcal{J}$.
Then
\[
|\mathcal{J}_X|=(1\pm\beta)\left(\prod_{i=2}^{|X|}d_i^{\binom{|X|}{i}}\right)m^{|X|}.
\]
\end{lemma}

Analogously, we have a dense version of Extension Lemma~\cite{2006Embeddings}.
\begin{lemma}[Dense Extension Lemma~\cite{2006Embeddings}]
\label{dense extension}
Let $k,s,s',m$ be positive integers, where $s'<s$ and $\varepsilon,\beta, d_2,\ldots,d_{k-1}$ be positive constants such that $1/m\ll\varepsilon\ll\beta\ll d_2,\ldots,d_{k-1},1/s$.
Let $H$ be a $(k-1)$-graph on $[s]$.
Let $\mathcal{H}$ be the $(k-1)$-complex generated by the down-closure of $H$ and $\mathcal{H}'$ be an induced subcomplex of $\mathcal{H}$ on $s'$ vertices.
Let $\textbf{d}=(d_2,\ldots,d_{k-1})$ and let $\mathcal{J}$ be a $(\textbf{d},\varepsilon)$-regular $(k-1)$-complex, with ground partition $\mathcal{P}$ with vertex classes of size $m$ each.
If $X_1,\ldots,X_s\in \mathcal{P}$, then all but at most $\beta|\mathcal{H}_{\mathcal{J}}'|$ labelled partition-respecting copies of $\mathcal{H}'$ in $\mathcal{J}$ extend to
\[
(1\pm\beta)\left(\prod_{i=2}^{k-1}d_i^{e_i(\mathcal{H})-e_i(\mathcal{H}')}\right)m^{s-s'}
\]
labelled partition-respecting copies of $\mathcal{H}$ in $\mathcal{J}$.
\end{lemma}

The restriction of a regular complex to a large subset of its vertex is also a regular complex, with slightly altered constants.
\begin{lemma}[Regular Restriction Lemma~\cite{Peter2017Tight}]
\label{restriction}
Let $k,r,m,s$ be integers and $\alpha,\varepsilon,\varepsilon_k,d_2,\ldots,d_{k}$ be positive constants such that $1/d_i\in \mathbb{N}$ for $\in[2,k]$ and
\[
1/m\ll\varepsilon\ll\varepsilon_k,d_2,\ldots,d_{k-1},
\]
and
\[
\varepsilon_k\ll\alpha.
\]
Let $\mathcal{G}$ be an $s$-partite $k$-complex on vertex classes $V_1,\ldots,V_s$, each of size $m$ and which is $(\textbf{d},\varepsilon_k,\varepsilon,r)$-regular where $\textbf{d}=(d_2,\ldots,d_k)$.
Choose any $V_i'\subseteq V_i$ with $|V_i'|\geq\alpha m$ for $i\in[s]$.
Then the induced subcomplex $\mathcal{G}[V_1'\cup\cdots\cup V_s']$ is $(\textbf{d},\sqrt{\varepsilon_k},\sqrt{\varepsilon},r)$-regular.
\end{lemma}
\section{Framework lemma}
In this section, we use the following Absorption Lemma and Almost Cover Lemma to prove Theorem \ref{1.9}.
The proof of these two lemmas will be found in Section 8 and 9.
Before we give these two lemmas, we need some definition.

\begin{definition}[Extensible paths]
Let $(G,G_{\mathcal{J}},\mathcal{J},\mathcal{P},R)$ be a $(k,m,2t,\varepsilon,\varepsilon_{k+1},r,\textbf{d})$-regular setup, $G$ be a $(1,k)$-graph on $[n]\cup V$ where $|V|=n$ and $c,\nu>0$.
A $(k-1)$-tuple $A$ in $V^{k-1}$ is said to be $(c,\nu)$-extensible rightwards to an ordered edge $Y=(Y_0,Y_1,\ldots,Y_k)$ in $R$ if there exists a connection $S\subseteq [n]\cup V$ and a target set $T\subseteq \mathcal{J}_{(Y_2,\ldots,Y_k)}$  with the following properties.
\begin{itemize}
  \item $|T|\geq\nu|\mathcal{J}_{(Y_2,\ldots,Y_k)}|$,
  \item for every $(v_2,\ldots,v_k)\in T$, there are at least $cm^{3k+1}$ many $(3k+1)$-tuples $(c_1,\ldots,c_{2k},w_1,\ldots w_k,v_1)$ with $v_1\in S \cap Y_1$, $w_i\in S \cap Y_i$ and $c_j\in Y_0$ for $i\in[k]$ and $j\in[2k]$ such that $(c_1\ldots c_{2k},Aw_1\ldots w_kv_1\ldots v_k)$ is a sequentially path in $G$.
\end{itemize}
\end{definition}

Given a sequentially path $P$ in a $(1,k)$-graph $G$ and an ordered edge $X$ in $R$, we say that $P$ is $(c,\nu)$-\emph{extensible rightwards} to $X$ if the $(k-1)$-tuple corresponding $P$'s last $k-1$ vertices is $(c,\nu)$-extensible rightwards to $X$.
We call $X$ as the right extension.
We can define leftwards path extensions for $(k-1)$-tuples and
for tight paths in an analogous way (this time corresponding to the first $k-1$ vertices of $P$).
A \emph{connection set} of a sequentially path is the union of the connection set of the initial $(k-1)$-tuple and the connection set of the end $(k-1)$-tuple.

Given that $X=(a,b,c)$ and $Y=(a,c,b)$, there is no guarantee that $H$ contains a walk from $X$ to $Y$.
While if $Y$ is a cyclic shift of $X$, that is, $(b,c,a)$ or $(c,a,b)$, then a walk from $X$ to $Y$ does exist.
More generally, a \emph{cyclic shift} of a tuple $(v_1,\ldots,v_k)$ is any $k$-tuple of the form $(v_i,\ldots,v_k,v_1,\ldots,v_{i-1})$ for $i\in[k]$.

An orientation of a $(1,k)$-graph $G$ on $[n]\cup V$ is a family of ordered $(1,k)$-tuples $\{\overrightarrow{e}\in [n]\times V^k:e\in E(G)\}$.
We say that a family $\overrightarrow{G}$ of ordered $(1,k)$-tuples is an \emph{oriented} $(1,k)$-\emph{graph} if there exists a $(1,k)$-graph $G$ such that $\overrightarrow{G}=\{\overrightarrow{e}\in [n]\times V^k:e\in E(G)\}$.
Given an oriented $(1,k)$-graph $\overrightarrow{R}$, we say that $(G,G_{\mathcal{J}},\mathcal{J},\mathcal{P},\overrightarrow{R})$ is an \emph{oriented} $(k,m,2t,\varepsilon,\varepsilon_{k+1},r,\textbf{d})$-\emph{regular setup} if $\overrightarrow{R}$ is an orientation of $R$ and $(G,G_{\mathcal{J}},\mathcal{J},\mathcal{P},R)$ is a $(k,m,2t,\varepsilon,\varepsilon_{k+1},r,\textbf{d})$-regular setup.
Consider a $(1,k)$-graph $G$ with an orientation $\overrightarrow{G}$ and vertex set $[n]\cup V$.
Given an ordered $k$-tuple $Y$ of distinct vertices in $V$ and $c\in[n]$, we say that $\{c\}\cup Y$ is \emph{consistent with} $\overrightarrow{G}$ if there exists an oriented edge $\{c\}\cup \overrightarrow{e}\in\overrightarrow{G}$ such that $\overrightarrow{e}$ is a cyclic shift of $Y$.
We say that an extensible path is \emph{consistent with} $\overrightarrow{G}$ if its left and right extensions are consistent with $\overrightarrow{G}$.
Finally, when considering multiple paths, we refer to the union of their connection sets as their \emph{joint connection set}.

Let $\overrightarrow{G}$ be an orientation of a $(1,k)$-graph $G$.
A sequentially walk $W$ in $G$ is said to be \emph{compatible} with $\overrightarrow{G}$ if each oriented edge of $\overrightarrow{G}$ appears at least once in $W$ as a sequence of $k$ consecutive vertices.

Let $G$ be a $(1,k)$-graph on $[n]\cup V$ where $|V|=n$, and $S\subseteq V, O\subseteq[n], |O|=|S|=k$, $P$ be a sequentially path.
Recall that $(C(P),I(P))$ is used to denote a sequentially path where $C(P)$ is the color set of $P$ and $I(P)$ is the point set of $P$.
We say that $P$ is $(S,O)$-\emph{absorbing} in $G$ if there exits a sequentially path $P'$ in $G$ with the same initial $(k-1)$-tuple and the same terminal $(k-1)$-tuple with $P$, $I(P')=I(P)\cup S$ and $C(P')=C(P)\cup O$.
We say that $P$ is $\eta$-\emph{absorbing} in $G$ if it is $(S,O)$-absorbing in $G$ for every $S$ of size at most $\eta n$ divisible by $k$, any $O$ of size $|S|$, and $S\cap I(P)=\emptyset, O\cap C(P)=\emptyset$.
\begin{lemma}[Absorption Lemma]
\label{absorption lemma}
Let $k,r,m,t\in \mathbb{N}$ and $d_2,\ldots,d_{k+1},\varepsilon,\varepsilon_{k+1},\eta,\mu,\delta,\alpha,c,\nu,\lambda$ be such that
\begin{equation}
\nonumber
\begin{split}
1/m&\ll1/r,\varepsilon\ll1/t,c,\varepsilon_{k+1},d_2,\ldots,d_k,\\
c&\ll d_2,\ldots,d_k,\\
1/t&\ll\varepsilon_{k+1}\ll d_{k+1},\nu\leq1/k,\\
c&\ll\varepsilon_{k+1}\ll\alpha\ll\eta\ll\lambda\ll\nu\ll\mu\ll\delta,1/k.
\end{split}
\end{equation}
Let $\textbf{d}=(d_2,\ldots,d_{k+1})$ and let $\mathfrak{S}=(G,G_{\mathcal{J}},\mathcal{J},\mathcal{P},\overrightarrow{H})$ be an oriented representative $(k,m,2t,\varepsilon,\varepsilon_{k+1},r,\textbf{d})$-regular setup.
Let $G$ be $(1,k)$-graph on $[n]\cup V$ with minimum relative $(1,1)$-degree being at least $\delta+\mu$ where $|V|=n$, $n\leq(1+\alpha)mt$.
Suppose that there exists a closed sequentially walk which is compatible with the orientation $\overrightarrow{H}$ of $H$ and
\begin{itemize}
  \item[(F1)] $H_i$ is sequentially tightly connected,
  \item[(F4)] For every color $i\in[t]$, there are at least $(1-\alpha)t$ points $v\in V$ such that $\{i,v\}$ has relative  $(1,1)$-degree at least $1-\delta+\gamma$.
\end{itemize}
Then there exists a sequentially path $P$ in $G$ such that the following holds.
\begin{itemize}
 \item[(1)] $P$ is $(c,\nu)$-extensible and consistent with $\overrightarrow{H}$,
 \item[(2)] $V(P)$ is $\lambda$-sparse in $\mathcal{P}$ and $V(P)\cap S=\emptyset$, where $S$ denotes the connection set of $P$,
 \item[(3)] $P$ is $\eta$-absorbing in $G$.
 \end{itemize}
\end{lemma}
\begin{lemma}[Almost Cover Lemma]
\label{almost cover}
Let $k,r,m,t\in \mathbb{N}$, $d_2,\ldots,d_{k+1},\varepsilon,\varepsilon_{k+1},\alpha,\gamma,c,\nu,\lambda$ be such that
\[
1/m\ll1/r,\varepsilon\ll1/t,c,\varepsilon_{k+1},d_2,\ldots,d_k,
\]
\[
c\ll d_2,\ldots,d_k,
\]
\[
1/t\ll \varepsilon_{k+1}\ll d_{k+1},\nu,\alpha\leq1/k,
\]
\[
\alpha\ll\eta\ll\lambda\ll\nu\ll\gamma.
\]
Let $\textbf{d}=(d_2,\ldots,d_{k+1})$ and let $\mathfrak{S}=(G,G_{\mathcal{J}},\mathcal{J},\mathcal{P},\overrightarrow{H})$ be an oriented $(k,m,2t,\varepsilon,\varepsilon_{k+1},r,\textbf{d})$-regular setup.
Suppose that $G$ is a $(1,k)$-graph on $[n]\cup V$ where $|V|=n$ and $n\leq(1+\alpha)mt$, $H$ is a $(1,k)$-graph on $[t]\cup V'$ where $|V'|=t$ and
\begin{itemize}
  \item[(F1)] $H_i$ is sequentially tightly connected,
  \item[(F2)] $H_i$ contains a sequentially closed walk $W$ compatible with $\overrightarrow{H}$ whose length is 1 mod $k$,
  \item[(F3)] $H_{W_i}$ is $\gamma$-robustly matchable for $i\in[k]$,
  \item[(F5)] $L_H(\{i\})$ and $L_H(\{j\})$ intersect in an edge for each $i,j\in[t]$.
\end{itemize}
Suppose that $P$ is a sequentially path in $G$ such that
\begin{itemize}
  \item[(1)] $P$ is $(c,\nu)$-extensible and consistent with $\overrightarrow{H}$,
  \item[(2)] $V(P)$ is $\lambda$-sparse in $\mathcal{P}$ and $V(P)\cap S=\emptyset$ where $S$ is the connection set of $P$,
\end{itemize}
then there exists a sequentially cycle $C$ of length at least $(1-\eta)n$ which contains $P$ as a subpath.
Moreover, the number of uncovered points of $V$ is divisible by $k$ and the number of uncovered colors of $[n]$ has the same size with the number of uncovered points.
\end{lemma}
\begin{proof}[The proof of Theorem \ref{1.9}]
Let $\delta=rhf_{k-2}(k)$, $\mu>0$ and
\[
\varepsilon_{k+1}\ll\alpha\ll\eta\ll\lambda\ll\nu\ll\gamma\ll\mu,
\]
\[
1/t_0\ll\varepsilon_{k+1}\ll d_{k+1}\ll\mu.
\]
We apply Lemma \ref{5.5} with input $\varepsilon_{k+1},1/t_0,r,\varepsilon$ to obtain $t_1,m_0$.
Choose $c\ll1/t_1$ and $1/n_0\ll1/t_1,1/m_0,c,1/r,\varepsilon$.
Let $G$ be a $(1,k)$-graph on $[n]\cup V$ where $|V|=n$ and $2n\geq n_0$ vertices with $\overline{\delta}_{1,k-2}(G)\geq\delta+\mu$.
Our goal is to prove that $G$ contains a sequentially Hamilton cycle.
By Lemma \ref{5.5}, there exists a representative $(k,m,2t,\varepsilon,\varepsilon_{k+1},r,d_2,\ldots,d_{k+1})$-regular setup $(G,G_{\mathcal{J}},\mathcal{J},\mathcal{P},R_{d_{k+1}})$ with $t_0\leq t\leq t_1$ and $n\leq(1+\alpha)mt$.
Moreover, there is a $(1,k)$-graph $I$ of edge density at most $\varepsilon_{k+1}$ such that $R=R_{d_{k+1}}\cup I$ has minimum relative $(1,k-2)$-degree at least $\delta+\mu/2$.
By Definition \ref{1.8} and $\delta=rhf_{k-2}(k)$, we obtain that $R$ contains an $(\alpha,\gamma,\delta)$-rainbow Hamilton framework $H$ that avoids edges of $I$.
Thus, $H\subseteq R_{d_{k+1}}$.

Next, we want to fix an orientation $\overrightarrow{H}$ and a compatible walk $W$.
Since $H$ is an $(\alpha,\gamma,\delta)$-rainbow Hamilton framework, $H_i$ is sequentially tightly connected and has a sequentially closed walk of length 1 mod $k$, $L_H(\{i\})$ and $L_H(\{j\})$ intersect in an edge for each $i,j\in[t]$.
We obtain a sequentially closed walk of length 1 mod $k$ visiting all edges of $H$.
Define an orientation $\overrightarrow{H}=\{\overrightarrow{e}\in V(H)^k:e\in H\}$ by choosing for every edge $e$ of $H$ a $k$-tuple (or subpath) $\overrightarrow{e}$ in $W$ which contains the vertices of $e$.
Note that $W$ is compatible with $\overrightarrow{H}$.

Firstly, we select a sequentially absorbing path $P$.
Note that $1/t_1\leq d_2,\ldots,d_k$, since $\mathcal{J}$ is a $(t_0,t_1)$-equitable complex.
Since $H$ is an $(\alpha,\gamma,\delta)$-rainbow Hamilton framework, it follows that there exists a sequentially path $P$ in $G$ by Lemma \ref{absorption lemma} such that
\begin{enumerate}
 \item $P$ is $(c,\nu)$-extensible and consistent with $\overrightarrow{H}$,
 \item $V(P)$ is $\lambda$-sparse in $\mathcal{P}$ and $V(P)\cap T=\emptyset$, where $T$ denotes the connection set of $P$,
 \item $P$ is $\eta$-absorbing in $G$.
 \end{enumerate}

Next, by Lemma \ref{almost cover}, there is a sequentially cycle $A$ of length at least $(1-\eta)n$ which contains $P$ as a subpath.
Moreover, the number of uncovered points $|V\setminus I(A)|$ is divisible by $k$ and the number of uncovered colors is of size $|[n]\setminus C(A)|=|V\setminus I(A)|$.

Finally, we absorb the uncovered points and colors  into $A$.
Note that $|V\setminus I(A)|\leq\eta n$.
Thus, there is a sequentially path $P'$ with point set $I(P)\cup (V\setminus I(A))$ and color set $C(P)\cup ([n]\setminus C(A))$, which has the same endpoints as $P$, as desired.
\end{proof}
\section{Almost Covering}

\subsection{Embedding sequentially paths}
Given sequentially walks $W$ and $W'$ with the property that the terminal $(k-1)$-tuple of $W$ is identical to the initial $(k-1)$-tuple of $W'$, we may \emph{concatenate} $W$ and $W'$ to form a new sequentially walk with color set $C(W)+ C(W')$, which we denote $W+W'$.
Note that a rainbow path in $k$-graph system is a sequentially path in the auxiliary $(1,k)$-graph $G$.


\begin{lemma}
\label{7.9}
Let $k,r,n_0,t,B$ be positive integers and $\psi,d_2,\ldots,d_{k+1},\varepsilon,\varepsilon_{k+1},\nu$ be positive constants such that $1/d_i\in \mathbb{N}$ for $i\in[2,k]$ and such that $1/n_0\ll1/t$,
\[
\frac{1}{n_0},\frac{1}{B}\ll\frac{1}{r},\varepsilon\ll\varepsilon_{k+1},d_2,\ldots,d_k,
\]
\[
\varepsilon_{k+1}\ll\psi,d_{k+1},\nu,\frac{1}{k}.
\]
Then the following holds for all integers $n\geq n_0$.

Let $G$ be a $(1,k)$-graph on $[n]\cup V$ where $|V|=n$, $\mathcal{J}$ be a $(\cdot,\cdot,\varepsilon,\varepsilon_{k+1},r)$-regular slice for $G$ on $[t]\cup V'$ where $|V'|=t$ with density vector $\textbf{d}=(d_2,\ldots,d_k)$.
Let
$\mathcal{J}_{W_i}$ be the induced subcomplex of $\mathcal{J}$ on $[t(i-1)/k+1,ti/k]\cup V'$ for $i\in[k]$.
We call $[t]$ the family of color clusters and $V'$ the family of point clusters.
Let $R_{W_i}:=R\left[[t(i-1)/k+1,ti/k]\cup V'\right]$ be the induced subgraph of $R:=R_{d_{k+1}}(G)$.
Let $R_{W_i}$ be sequentially tightly connected for $i\in[k]$ and $\textbf{w}_i$ be a fractional matching of size $\mu_i=\sum_{e\in E(R_{W_i})}\textbf{w}_i(e)$ for $i\in[k]$ and $\mu_i(Z)=\sum_{Z\in e, e\in E(R_{W_i})}\textbf{w}_i(e)\leq1/k$ for each cluster $Z$.
Also, let $X$ and $Y$ be $(k-1)$-tuples of point clusters, $S_X$ and $S_Y$ be the subsets of $\mathcal{J}_X$ and $\mathcal{J}_Y$ of sizes at least $\nu|\mathcal{J}_X|$ and $\nu|\mathcal{J}_Y|$ respectively.
Finally, let $W$ be a sequentially walk from $X$ to $Y$ of length at most $t^{2k+1}$ in $R_{W_i}$ and denote $\ell(W)$ by $p$.
For $i\in[k]$, we have
\begin{enumerate}
  \item[(i)] for any $\ell$ divisible by $k$ with $4k\leq\ell\leq(1-\psi)\mu_i kn/t$, there is a sequentially path $P$ in $G$ of length $\ell-1+\ell(W)(k+1)$ whose initial $(k-1)$-tuple belongs to $S_X$ and whose terminal $(k-1)$-tuple belongs to $S_Y$,
  \item[(ii)] $P$ uses at most $\mu_i(Z)n/t+B$ vertices from any point cluster $Z\in V'$ and at most $k\mu_i(C)n/t+B$ vertices from any color cluster $C\in[t(i-1)/k+1,ti/k]$ where $\mu_i(Z')=\sum_{Z'\in e,e\in R_{W_i}}\textbf{w}_i(e)$ for any cluster $Z'$.
\end{enumerate}
\end{lemma}

\begin{proof}
Let $\alpha=\psi/5$ and $\beta=1/200$.
When using Lemma \ref{extension}, we require that $\varepsilon\ll c^2$ and choose $m_0$ to be large enough so that $m\geq\alpha m_0$ is acceptable for all these applications.
Given $t$, let
\begin{equation}
\label{parameter}
n_0=t\cdot\max(m_0,\frac{200k^2}{\varepsilon},\frac{8k^2}{\alpha\sqrt{\varepsilon}},\frac{10k(k+1)t^{2k+1}}{\alpha}).
\end{equation}
We write $\mathcal{G}$ for the $(k+1)$-complex obtained from $\mathcal{J}_{W_i}$ by adding all edges of $G$ supported on $\mathcal{J}^{(k)}_{W_i}$ as the `$(k+1)$th level' of $\mathcal{G}$.
So for any edge $X=(X_0,X_1,\ldots,X_k)\in R_{W_i}$, $\mathcal{G}[\bigcup_{i\in [0,k]}X_i]$ is a $(d_2,\ldots,d_k,d^*(X),\varepsilon,\varepsilon_{k+1},r)$-regular $(k+1)$-partite $(k+1)$-complex with $d^*(X)\geq d_{k+1}$.

Since $\mathcal{J}$ is a regular slice for $G$, for any $(1,k)$-set of clusters $X=\{X_0,X_1,\ldots,X_k\}$ in $\mathcal{J}_{W_i}$, the $(k+1)$-partite $k$-complex $\mathcal{J}_{W_i}[\bigcup_{j\in[0,k]}X_j]$ is $(\textbf{d},\varepsilon)$-regular.
By adding all $(k+1)$-sets supported on $\hat{\mathcal{J}_{W_i}}_X$ as the `$(k+1)$th level', we may obtain a $(d_2,\ldots,d_k,1,\varepsilon,\varepsilon_{k+1},r)$-regular $(k+1)$-partite $(k+1)$-complex, whose vertex clusters are subsets $Y_j\subseteq X_j$ for $j\in[0,k]$ of size $|Y_1|=\cdots=|Y_k|=\alpha m/k$ and $|Y_0|=\alpha m$.
$Y_0$ can be seen as $\bigcup_{i\in[k]}{Y_{0,i}}$ where $|Y_{0,i}|=\alpha m/k$ for $i\in[k]$ and we obtain a $(d_2,\ldots,d_k,1,\sqrt{\varepsilon},\sqrt{\varepsilon_{k+1}},r)$-regular by Lemma \ref{restriction}.
We conclude by Lemma \ref{7.7} that for any subset $Y_i$, $i\in[k-1]$ of distinct clusters of $\mathcal{J}$, each of size $\alpha m$, we have
\begin{equation}
\label{31}
|\mathcal{G}(Y_1,\ldots,Y_{k-1})|\geq\varepsilon m^{k-1}.
\end{equation}

The following claim plays an important role in Lemma \ref{7.9}.
\begin{claim}
\label{37}
Let $\{X_0,X_1,\ldots,X_k\}$ be an edge of $R$ and choose any $Y_j\subseteq X_j$ for each $j\in[0,k]$ so that $|Y_0|=k|Y_1|=\cdots=k|Y_k|=\alpha m$.
Let $\mathcal{P}$ be a collection of at least $\frac{1}{2}|\mathcal{G}(Y_1,\ldots,Y_{k-1})|$ sequentially paths in $G$(not necessarily contained in $\bigcup_{j\in[k]}Y_j$) each of length at most $3k$ and whose terminal $(k-1)$-tuples are distinct members of $\mathcal{G}(Y_1,\ldots,Y_{k-1})$.
Then for each $\sigma\in\{0,1\}$ there is a path $P\in\mathcal{P}$ and a collection $\mathcal{P}'$ of $\frac{9}{10}e(\mathcal{G}(Y_{\sigma+1},\ldots,Y_{\sigma+k-1}))$
sequentially paths in $G$, each of length $2k-1+\sigma$, all of whose initial $(k-1)$-tuples are the same (terminal $(k-1)$-tuple of $P$).
Furthermore, the terminal $(k-1)$-tuples of paths in $\mathcal{P}'$ are distinct members of $\mathcal{G}(Y_{\sigma+1},\ldots,Y_{\sigma+k-1})$.
If $j\leq k-1$, then the $j$th vertex $x$ of each path in $\mathcal{P}'$ lies in $Y_j$, if $j\geq k$, then $x$ is not contained in $P$, and $k$ new colors are not contained in $P$.
\end{claim}
\begin{proof}
Let $\sigma\in\{0,1\}$ be fixed, we take $\mathcal{H}$ to be the $(k+1)$-complex generated by the down-closure of a sequentially path of length $2k-1+\sigma$ with vertex set $\{c_1,\ldots,c_{k+\sigma}\}\cup\{v_1,\ldots,v_{2k-1+\sigma}\}$ and consider its $(k+1)$-partition $V_0\cup V_1\cup\cdots\cup V_k$ where $\{c_1,\ldots,c_{k+\sigma}\}\subseteq V_0$ and the $i$th vertex of the path lies in the vertex class $V_j$ with $j=i$ mod $k$.
We take $\mathcal{H}'$ to be the subcomplex of $\mathcal{H}$ induced by $\{v_1,\ldots,v_{k-1},v_{k+1+\sigma},\ldots,v_{2k-1+\sigma}\}$.
Consider the pair $(e,f)$, where $e$ is an ordered $(k-1)$-tuple of $\mathcal{G}(Y_1,\ldots,Y_{k-1})$ and $f$ is an ordered $(k-1)$-tuple of $\mathcal{G}(Y_{\sigma+1},\ldots,Y_{\sigma+k-1})$.
For any such ordered $(k-1)$-tuple $e$, there are at most $km^{k-2}$ such ordered $(k-1)$-tuples $f$ which intersect $e$, thus there are at most $1/200$-proportion of the pairs $(e,f)$ are not disjoint.
On the other hand, if $e$ and $f$ are disjoint, then the down-closure of the pair $(e,f)$ forms a labelled copy of $\mathcal{H}'$ in $\mathcal{G}[\bigcup_{j\in[0,k]}Y_j]$, so by Lemma \ref{extension} with $s=3k+2\sigma-1$ and $s'=2k-2$, for all but at most $1/200$-proportion of the disjoint pairs $(e,f)$, there are at least $c(\alpha m/k)^{k+2\sigma+1}\geq \sqrt{\varepsilon}(\alpha m/k)^{k+2\sigma+1}$ extensions to copies of $\mathcal{H}$ in $\mathcal{G}[\bigcup_{j\in[0,k]}Y_j]$.
Each such copy of $\mathcal{H}$ corresponds to a sequentially path in $G$ of length $2k-1+\sigma$ with all vertices in the desired clusters.
We conclude that at least $99/100$-proportion of all pairs $(e,f)$ of ordered $(k-1)$-tuples are disjoint and are linked by at least $\sqrt{\varepsilon}(\alpha m/k)^{k+2\sigma+1}$ sequentially paths in $G$ of length $2k-1+\sigma$, where $c_i\in V_0$ for $i\in[k+\sigma]$ and $v_{\ell}\in V_j$ with $j=\ell$ mod $k$.
We call these pairs \emph{extensible}.

We call an ordered $(k-1)$-tuple $e\in \mathcal{G}(Y_1,\ldots,Y_{k-1})$ \emph{good} if at most $1/20$ of the ordered edges $f\in\mathcal{G}(Y_{\sigma+1},\ldots,Y_{\sigma+k-1})$ do not make an extensible pair with $e$.
Then at most $1/5$ of the ordered $(k-1)$-tuples in $\mathcal{G}(Y_1,\ldots,Y_{k-1})$ are not good.
Thus, there exists a path $P\in \mathcal{P}$ whose terminal $(k-1)$-tuple is a good ordered $(k-1)$-tuple $e$.
Fix such a $P$ and $e$, and any ordered $(k-1)$-tuple $f$ in $\mathcal{G}(Y_{\sigma+1},\ldots,Y_{\sigma+k-1})$ which is disjoint from $P$, suppose that $(e,f)$ is an extensible pair, there are at least $\sqrt{\varepsilon}(\alpha m/k)^{k+2\sigma+1}$ sequentially paths in $G$ from $e$ to $f$.
We claim that at least one of these paths has the further property that if $j\geq k$, then the $j$th vertex is not contained in $P$ and the $k+\sigma$ new colors are not contained in $P$, we can therefore put it in $\mathcal{P}'$.
Indeed as $f$ is disjoint from $P$, if $\sigma=0$, then it suffices to show that one of these paths has the property that $v_k\in Y_k\setminus V(P)$ and $c_i\in Y_0\setminus V(P)$ for $i\in[k]$.
This is true because there are only at most $(2k+1)(\alpha m)^k+k(2k+1)(\alpha m/k)^k<\sqrt{\varepsilon}(\alpha m/k)^{k+1}$ paths which do not have this property  by (\ref{parameter}).
If $\sigma=1$, then we need a path whose $k$th and $(k+1)$st vertices are not in $V(P)$ and $c_i\in Y_0\setminus V(P)$ for $i\in[k+1]$, which is possible since $2(2k+1)(\alpha m/k)^{k+2}+(k+1)(2k+1)(\alpha m/k)^{k+2}<\sqrt{\varepsilon}(\alpha m/k)^{k+3}$ by (\ref{parameter}).

Finally, considering the ordered $(k-1)$-tuple $f\in\mathcal{G}(Y_{\sigma+1},\ldots,Y_{\sigma+k-1})$, we have $20|V(P)|(k-1)(\alpha m/k)^{k-2}$
$\leq\varepsilon m^{k-1}\leq e(\mathcal{G}(Y_{\sigma+1},\ldots,Y_{\sigma+k-1}))$ by (\ref{parameter}) and (\ref{31}), at most $1/20$ of these $(k-1)$-tuples $f$ intersect $P$ and by the choice of $e$, at most $1/20$ of these $(k-1)$-tuples $f$ are such that $(e,f)$ is not extensible.
This leaves at least $9/10$ of $(k-1)$-tuples $f$ remaining, and choose a sequentially path for each such $f$ as described above gives the desired set $\mathcal{P}'$.
\end{proof}

Let $X=(X_1,\ldots,X_{k-1})$, $Y=(Y_1,\ldots,Y_{k-1})$, $X_k$ be the cluster following $X$ in $W$ and $Y_k$ be the cluster preceding $Y$ in $W$.
Without loss of generality, we may assume that $\{X_0,X_1,\ldots,X_k\}$ is an edge of $R_1$ and $\{Y_0,Y_1,\ldots,Y_k\}$ is an edge of $R_k$.
By the condition, we have $S_X$ constitutes at least a $\nu$ proportion of $\mathcal{G}(X_1,\ldots,X_{k-1})$ and $S_Y$ constitutes at least a $\nu$ proportion of $\mathcal{G}(Y_1,\ldots,Y_{k-1})$.
Given any subsets $X_j'\subseteq X_j$ of size $\alpha m/k$ for $j\in[k]$ and $X_0'\subseteq X_0$ of size $\alpha m$, we say that a $(k-1)$-tuple $e\in\mathcal{G}(X_1,\ldots,X_{k-1})$ is \emph{well}-\emph{connected to} $(X_1',\ldots,X_{k-1}')$ via $X_k'$ and $X_0'$ if for at least $9/10$ of the $(k-1)$-tuples $f$ in $\mathcal{G}(X_1',\ldots,X_{k-1}')$, there exist distinct $k$-subsets $\{c_1,\ldots,c_k\}$, $\{f_1,\ldots,f_k\}$ of $X_0'$ and distinct $u,v\in X_k'$ such that $(c_1\cdots c_k,e(u)f)$ and $(f_1\cdots f_k,e(v)f)$ are sequentially paths in $G$ of length $2k-1$.
\begin{claim}
\label{38}
For any subsets $X_j'\subseteq X_j$ of size $\alpha m/k$, $Z_j\subseteq X_j$ of size $\alpha m/k$ for $j\in[k]$ and $X_0'\subseteq X_0$, $Z_0\subseteq X_0$ of size $\alpha m$   such that each $X_j'$ is disjoint from $Z_j$, the following statements hold.
\begin{itemize}
  \item[(1)] At least $9/10$ of the $(k-1)$-tuples $e$ in $\mathcal{G}(Z_1,\ldots,Z_{k-1})$ are well-connected to $(Z_1,\ldots,Z_{k-1})$ via $Z_k$ and $Z_0$.
  \item[(2)] At least $9/10$ of the $(k-1)$-tuples $e$ in $\mathcal{G}(Z_1,\ldots,Z_{k-1})$ are well-connected to $(X_1',\ldots,X_{k-1}')$ via $X_k'$ and $X_0'$.
  \item[(3)] At least $9/10$ of the $(k-1)$-tuples $e$ in $\mathcal{G}(X_1',\ldots,X_{k-1}')$ are well-connected to $(Z_1,\ldots,Z_{k-1})$ via $X_k'$ and $X_0'$.
\end{itemize}
\begin{proof}
From the proof of Claim \ref{37}, we know that all but at most $1/100$-proportion of pairs $(e,f)$, where $e,f\in\mathcal{G}(Z_1,\ldots,Z_{k-1})$, are disjoint and are linked by at least $\sqrt{\varepsilon}(\alpha m/k)^{k+1}$ sequentially tight paths in $G$ of length $2k-1$.
It is obvious that at least $9/10$-proportion $(k-1)$-tuples of $\mathcal{G}(Z_1,\ldots,Z_{k-1})$ can be extended to at least $9/10$-proportion $(k-1)$-tuples of $\mathcal{G}(Z_1,\ldots,Z_{k-1})$ by at least $\sqrt{\varepsilon}(\alpha m/k)^{k+1}$ sequentially paths.
To prove (2), we apply Lemma \ref{dense extension} with $\mathcal{H}$ being the $(k+1)$-complex generated by the down-closure of a sequentially path of length $2k-1$ and $\mathcal{H}'$ being the subcomplex induced by its initial and terminal $(k-1)$-tuples.
We regard $\mathcal{H}$ as a $(2k)$-partite $(k+1)$-complex with $k$ colors in the color cluster and one vertex in each point cluster.
The role of $\mathcal{G}$ in Lemma \ref{extension} is the $(2k)$-partite subcomplex of $\mathcal{G}$ with vertex classes $X_0',Z_1,\ldots,Z_{k-1},X_k',X_1',\ldots,X_{k-1}'$, the colors of $\mathcal{H}$ are embedded in $X_0'$, the first vertex of $\mathcal{H}$ is to be embedded in $Z_1$, the second one in $Z_2$, and so forth.
By Lemmas \ref{restriction} and \ref{extension}, the proportion of pairs $(e,f)$ for which there is no path as in (2) is at most $1/200$, and the remainder of the argument can be followed in (1).
(3) can be proved similarly.
\end{proof}
\end{claim}

We are ready to construct our path.
Arbitrarily choose a subset $X_0^{(0)}\subseteq X_0$, $Z_0\subseteq Y_0$ of size $\alpha m$ and $X_j^{(0)}\subseteq X_j$, $Z_j\subseteq Y_j$ of size $\alpha m/k$ for $j\in[k]$.
By Theorem~\ref{counting lemma}, Theorem~\ref{extension}, Theorem~\ref{7.7}, there are at least $|S_X||\mathcal{G}(X_1^{(0)},\ldots,X_{k-1}^{(0)})|/2$ pairs $(e,f)$, where $e\in S_X$ and $f\in\mathcal{G}(X_1^{(0)},\ldots,X_{k-1}^{(0)})$, can be extended to $\sqrt{\varepsilon}(\alpha m/k)^{k+1}$ sequentially paths whose remaining point lies in $X_k^{(0)}$ and colors lie in $X_0^{(0)}$.
Thus, we choose a $(k-1)$-tuple $P^{(0)}$ of $S_X$ such that the following holds,
there is a set $\mathcal{P}^{(0)}$ of sequentially paths of the form $(c_1\cdots c_k,P^{(0)}(v)f)$ for $v\in X_k^{(0)}$, $c_1,\ldots,c_k\in X_0^{(0)}$ and $f\in \mathcal{G}(X_1^{(0)},\ldots,X_{k-1}^{(0)})$ for which the terminal $(k-1)$-tuples of paths in $\mathcal{P}^{(0)}$ are all distinct and constitute at least half of the ordered $(k-1)$-tuples of $\mathcal{G}(X_1^{(0)},\ldots,X_{k-1}^{(0)})$.
Similarly, we can choose $e\in S_Y$ such that for at least half the members $e'$ of $\mathcal{G}(Z_1,\ldots,Z_{k-1})$, there is a sequentially path of length $2k-1$ in $G$ from $e'$ to $e$ whose remaining point lies in $Z_k$ and colors lie in $Z_0$.

We now construct the desired path.
Since $H_{W_i}$ is sequentially tightly connected, we can obtain $W=e_1\cdots e_s$ passing all edges of $H_{W_i}$.
For each $i\in[s]$, let $n_i$ be any integer with $0\leq n_i\leq(1-3\alpha)\textbf{w}({e_i})m$.
Set the initial state to be `filling the edge $e_1$', we proceed for $j\geq1$ as follows,
\begin{description}
  \item[$\bigstar$] The terminal $(k-1)$-tuple of the path family $\mathcal{P}^{(j)}$ constitute at least half of the ordered $(k-1)$-tuples $\mathcal{G}(X_1^{(j)},\ldots,X_{k-1}^{(j)})$.
\end{description}

Suppose that our current state is `filling the edge $e_i$' for some $i$, if we have previously completed $n_i$ steps in this state, then we do nothing and change the state to `position 1 in traversing the walk $W$'.
Otherwise, since $\bigstar$ holds for $j-1$, we apply Claim \ref{37} with $\sigma=0$ to obtain a path $P\in\mathcal{P}^{(j-1)}$ and a collection $\mathcal{P}^{(j)}$ of $\frac{9}{10}e(\mathcal{G}(X_1^{(j-1)},\ldots,X_{k-1}^{(j-1)}))$ sequentially paths of length $2k-1$, all of whose initial $(k-1)$-tuples are the same (the terminal $(k-1)$-tuple of $P$) and whose terminal $(k-1)$-tuples are distinct numbers of $\mathcal{G}(X_1^{(j-1)},\ldots,X_{k-1}^{(j-1)})$ and are disjoint from $V(P)$,
whose colors lie in $X_0^{(j-1)}\setminus C(P)$,
and whose remaining vertex lies in $X_k^{(j-1)}\setminus V(P)$.
We define $P^{(j)}$ to be the concatenation $P^{(j-1)}+P$ with color classes $C(P^{(j-1)})\cup C(P)$.
For $p\in[0,k]$, we generate $X_p^{(j)}$ from $X_p^{(j-1)}$ by removing the vertices of $P^{(j)}$ in $X_p^{(j)}$ and replacing them by vertices from the same cluster which do not lie in $Z$ or in $P^{(j)}$.
We will prove that this is possible in Claim \ref{39}.

Now suppose that our current state is `position $q$ in traversing the walk $W$'.
Since $\bigstar$ holds for $j-1$, applying Claim \ref{37} with $\sigma=1$ to obtain a path $P\in\mathcal{P}^{(j-1)}$ and a collection $\mathcal{P}^{(j)}$ of $\frac{9}{10}e(\mathcal{G}(X_1^{(j-1)},\ldots,X_{k-1}^{(j-1)}))$ sequentially paths of length $2k$, all of whose initial $(k-1)$-tuples are the same (the terminal $(k-1)$-tuple of $P$) and whose terminal $(k-1)$-tuples are distinct numbers of $\mathcal{G}(X_2^{(j-1)},\ldots,X_{k}^{(j-1)})$ and are disjoint from $V(P)$, and whose two remaining vertices lie in $X_k^{(j-1)}\setminus V(P)$ and $X_1^{(j-1)}\setminus V(P)$ respectively with colors in $X_0^{(j-1)}\setminus C(P)$.
Exactly as before we define $P^{(j)}$ to be the concatenation $P^{(j-1)}+P$.
We generate $X_p^{(j)}$ from $X_{p+1}^{(j-1)}$ for $p\in[0,k-1]$ by removing the vertices of $P^{(j-1)}$ in $X_{p+1}^{(j-1)}$ and replacing them by vertices from the same cluster do not lie in $Z$ or $P^{(j)}$.
If we have not reached the end of $W$, we choose $X_k^{(j)}$ to be a subset of the cluster at position $q+k$ in the sequence of $W$ such that $X_k^{(j)}$ is disjoint from $P^{(j)}\cup Z$.
In this case, we change our state to `position $q+1$ in traversing $W$'.
Alternatively, if we have reached the end of $W$, meaning that the $(k-1)$-tuple of clusters containing $X_1^{(j)},\ldots,X_{k-1}^{(j)}$ is $(Y_1\ldots,Y_{k-1})$, then we choose $X_k^{(j)}$ to be a subset of $Y_k$ which has size $\alpha m/k$ and is disjoint from $P^{(j)}\cup Z$.
We may choose a path $P\in \mathcal{P}^{(j-1)}$ such that the terminal $(k-1)$-tuple $f\in G(X_1^{(j)},\ldots,X_{k-1}^{(j)})$ of $P$ is well-connected to $(Z_1,\ldots,Z_{k-1})$ via $Z_k$ and $Z_0$.
This implies that we may choose a $(k-1)$-tuple $e'$ in $\mathcal{G}(Z_1,\ldots,Z_{k-1})$, $v,v'$ in $Z_k$ with new colors $C^*, C^{**}$ in $Z_0$ with $|C^*|=|C^{**}|=k$ such that $(C^*,f(v')e')$ is a sequentially path $Q'$ and $(C^{**},e'(v)e)$ is a sequentially path $Q$.
Return $P^{(j)}+Q'+Q$ as the output sequentially path in $G$.
Note that an edge may appear multiple times.
When it first appears in the walk, the process executes `filling the edge'.
When it appears later, `filling the edge' is no longer needed.
Again we prove Claim \ref{39} that these choices are all possible.


\begin{claim}
\label{39}
The algorithm described above is well-defined(that is, it is always possible to construct the sets $X_p^{(j)}$), maintains $\bigstar$ and returns a sequentially path of length
\[
4k-1+\left(\sum_{i\in[s]}n_i\right)\cdot k+\ell(W)\cdot(k+1).
\]
\end{claim}
\begin{proof}

We prove that $\bigstar$ is maintained, recall that $e(\mathcal{G}(X_1^{(j)},\ldots,X_{k-1}^{(j)}))\geq\varepsilon m^{k-1}$ for each $j$.
Fixing some $j$, for either $A_p:=X_{p}^{(j-1)}$ or $A_p:=X_{p+1}^{(j-1)}$, we obtain sets $A_1,\ldots,A_{k-1}$, each  with size $\alpha m$ such that the terminal $(k-1)$-tuples of $\mathcal{P}^{(j)}$ constitute at least $9/10$ of the ordered edges of $\mathcal{G}(A_1,\ldots,A_{k-1})$ and for each $i\in[k-1]$, $X_i^{(j)}$ is formed from $A_i$ by removing at most two vertices and replacing them with the same number of vertices.
Since each vertex is in at most $m^{k-2}$ ordered $(k-1)$-tuples of either $\mathcal{G}(A_1,\ldots,A_{k-1})$ or $\mathcal{G}(X_1^{(j)},\ldots,X_{k-1}^{(j)})$, we conclude that the fraction of ordered $(k-1)$-tuples of $\mathcal{G}(X_1^{(j)},\ldots,X_{k-1}^{(j)})$ which are the terminal $(k-1)$-tuples of paths in $\mathcal{P}^{(j)}$ is at least
\begin{equation}
\begin{split}
&\frac{\frac{9}{10}e(\mathcal{G}(A_1,\ldots,A_{k-1}))-2(k-1)m^{k-2}}{e(\mathcal{G}(X_1^{(j)},\ldots,X_{k-1}^{(j)}))}\\
&\geq\frac{\frac{9}{10}(e(\mathcal{G}(X_1^{(j)},\ldots,X_{k-1}^{(j)}))-2(k-1)m^{k-2})-2(k-1)m^{k-2}}{e(\mathcal{G}(X_1^{(j)},\ldots,X_{k-1}^{(j)}))}\\
&\geq\frac{9}{10}-\frac{4(k-1)m^{k-2}}{\varepsilon m^{k-1}}\geq\frac{1}{2},
\end{split}
\end{equation}
where the last equality holds since $m\geq m_0\geq16(k-1)/\varepsilon$.
Thus, we obtain $\bigstar$.

To prove that we can always construct the set $X_p^{(j)}$, observe that it is enough to check that at termination every cluster still have at least $2\alpha m$ vertices not in $P^{(j)}$, as then there are at least $\alpha m$ vertices outside $Z$.
In each walk-traversing step, each path in $\mathcal{P}^{(j)}$ contains precisely $k+1$ new vertices and $k+1$ new colors and the total number of walk-traversing steps is precisely $\ell(W)$.
Recall that this number is at most $t^{2k+1}$, we have $(k+1)t^{2k+1}<\frac{\alpha m}{2k}$ and $(k+1)^2t^{2k+1}<\frac{\alpha m}{2}$ by (\ref{parameter}).
When we are in the state `filling the edge $e_i$', we have $n_i$ steps and in each step, each path in $\mathcal{P}^{(j)}$ contains $k$ new vertices, one from each cluster of $e_i\setminus C(e_i)$ and $k$ new colors from $C(e_i)$.
So for any color cluster $C$, the number of whose vertices which are added to $P^{(j)}$ is at most $\sum_{i:C\in e_i}kn_i\leq\sum_{i:C\in e_i}(1-3\alpha)k\textbf{w}(e_i)m\leq(1-3\alpha)m$.
And for any point cluster $X$, the number of whose vertices which are added to $P^{(j)}$ is at most $\sum_{i:X\in e_i}n_i\leq\sum_{i:X\in e_i}(1-3\alpha)\textbf{w}(e_i)m\leq(1-3\alpha)m/k$.
Together with $e$ and the $k$ vertices of the chosen path in $\mathcal{P}^{(0)}$, we conclude that there are at most $(1-2\alpha)m$ vertices of any color cluster and at most $(1-2\alpha)m/k$ vertices of any point cluster contained in $P^{(j)}$ at termination.

Finally, the length of the path is equal to the number of vertices.
Recall that $P^{(0)}$ contains $k-1$ vertices.
Next, $k$ vertices and $k$ colors are added from $P^{(0)}$ to form $P^{(1)}$.
Each of the $\sum_{i\in[s]}n_i$ edge-filling steps resulted in $k$ new vertices and $k$ new colors being added to $P^{(j)}$ and each of the $\ell(W)$ walk-traversing steps resulted in $k+1$ new vertices and $k+1$ new colors being added to $P^{(j)}$.
When completing the path, we need $2k$ vertices which are not in the final paths $P^{(j)}$ ($v,v',e$ and $e'$).
Thus, the final path has length
\[
(k-1)+k+\left(\sum_{i\in[s]}n_i\right)\cdot k+\ell(W)\cdot(k+1)+2k.
\]
\end{proof}

We obtain the shortest sequentially path by never entering the state `filling an edge', in which case we can obtain a sequentially path of length $4k-1+\ell(W)(k+1)$.
On the other hand, by extending $W$ to include all edges of $R_{W_i}$, we take $n_i$ to be $(1-\psi)\textbf{w}(e_i)m$ for each $i\in[s]$.
We can obtain a sequentially path of length at least $(1-\psi)\mu_i kn/t$, with using at most $k\mu_i(C)n/t+B$ vertices from any color cluster $C$ in $R_{W_i}$ and at most $\mu_i(X)n/t+B$ where $\mu_i(Z)=\sum_{Z\in e,e\in R_{W_i}}\textbf{w}_i(e)$ for $i\in[k]$ and $B=B(t,k)$.
By choosing $n_i$ appropriately, we can obtain tight cycles of certain length between two extremes.
\end{proof}

Similarly with Lemma \ref{7.9}, we can obtain the following lemma.
\begin{lemma}
\label{coro}
Let $k,r,n_0,t,B$ be positive integers and $\psi,d_2,\ldots,d_{k+1},\varepsilon,\varepsilon_{k+1},\nu$ be positive constants such that $1/d_i\in \mathbb{N}$ for $i\in[2,k]$ and such that $1/n_0\ll1/t$,
\[
\frac{1}{n_0}\ll\frac{1}{t}\ll\frac{1}{B}\ll\frac{1}{r},\varepsilon\ll\varepsilon_{k+1},d_2,\ldots,d_k,
\]
\[
\varepsilon_{k+1}\ll\psi,d_{k+1},\nu,\frac{1}{k}.
\]
Then the following holds for all integers $n\geq n_0$.

Let $G$ be a $(1,k)$-graph on $[n]\cup V$ where $|V|=n$, $\mathcal{J}$ be a $(\cdot,\cdot,\varepsilon,\varepsilon_{k+1},r)$-regular slice for $G$ on $[t]\cup V'$ where $|V'|=t$ with density vector $\textbf{d}=(d_2,\ldots,d_k)$.
Let
$\mathcal{J}_{W_i}$ be the induced subcomplex of $\mathcal{J}$ on $[t(i-1)/k+1,ti/k]\cup V'$ for $i\in[k]$.
Let $R_{W_i}:=R\left[[t(i-1)/k+1,ti/k]\cup V'\right]$ be the induced subgraph of $R:=R_{d_{k+1}}(G)$.
Let $R_{W_i}$ be sequentially tightly connected for $i\in[k]$ and $\textbf{w}_i$ be a fractional matching of size $\mu_i=\sum_{e\in E(R_{W_i})}\textbf{w}_i(e)$ for $i\in[k]$ with $\mu_i(Z)\leq1/k$ for each cluster $Z$ and $i\in[k]$.
Also, let $X$ and $Y$ be $(k-1)$-tuples of point clusters, $S_X$ and $S_Y$ be the subsets of $\mathcal{J}_X$ and $\mathcal{J}_Y$ of sizes at least $\nu|\mathcal{J}_X|$ and $\nu|\mathcal{J}_Y|$ respectively.
Finally, let $W$ be a sequentially walk traversing all edges of each $H_{W_i}$ from $X$ to $Y$ of length at most $t^{2k+1}$ and denote $\ell(W)$ by $p$.
For $i\in[k]$, we have
\begin{enumerate}
  \item for any $\ell$ divisible by $k$ with $4k\leq\ell\leq(1-\psi)\sum_{i\in[k]}\mu_i kn/t$, there is a sequentially path $P$ in $G$ of length $\ell-1+\ell(W)(k+1)$ whose initial $(k-1)$-tuple belongs to $S_X$ and whose terminal $(k-1)$-tuple belongs to $S_Y$,
  \item $P$ uses at most $\sum_{i\in[k]}\mu_i(Z)n/t+B$ vertices from any point cluster $Z\in V'$ and at most $k\mu_i(C)n/t+B$ vertices from any color cluster $C\in[t]$ where $\mu_i(Z')=\sum_{Z'\in e,e\in R_{W_i}}\textbf{w}_i(e)$ for any cluster $Z'$.
\end{enumerate}
\end{lemma}
\subsection{Connecting}
Let us begin with the existence of extensible paths. The following proposition states that
most tuples in the complex induced by an edge of the reduced graph of a regular slice also extend
to that edge.
\begin{proposition}
\label{8.1}
Let $k,m,t,r\in \mathbb{N}$ and $\varepsilon,\varepsilon_{k+1},d_2,\ldots,d_{k+1},\beta,c,\nu$ be such that
\[
1/m\ll1/r,\varepsilon\ll c\ll\varepsilon_{k+1},d_2,\ldots,d_k,
\]
\[
\varepsilon_{k+1}\ll\beta\ll d_{k+1},\nu.
\]
Let $\textbf{d}=(d_2,\ldots,d_{k+1})$ and let $(G,G_{\mathcal{J}},\mathcal{J},\mathcal{P},R)$ be a $(k,m,2t,\varepsilon,\varepsilon_{k+1},r,\textbf{d})$-regular setup.
Let $Y=(Y_0,Y_1,\ldots,Y_k)$ be an ordered edge in $R$, then all but at most $\beta|\mathcal{J}_{(Y_1,\ldots,Y_{k-1})}|$ many tuples $(v_1,\ldots,v_{k-1})$
$\in\mathcal{J}_{(Y_1,\ldots,Y_{k-1})}$ are $(c,\nu)$-extensible both left and rightwards to $Y$.
\end{proposition}
\begin{proof}
Let $P=(c_1,\ldots,c_{2k},v_1,\ldots,v_{3k-1})$ be a sequentially path. Partition its vertex set in $k+1$ clusters $X_0,X_1,\ldots,X_k$ such that $X_0=\{c_1,\ldots,c_{2k}\}$, and $X_i=\{v_j:j=i\ {\rm mod}\ k\}$ for $i\in[k]$.
Thus, $P$ is a $(k+1)$-partite $(k+1)$-graph.

Let $\mathcal{H}$ be the down-closure of the path $P$, which is a $(k+1)$-partite $(k+1)$-complex.
Let $V_1=\{v_1,\ldots,v_{k-1}\}$ and $V_2=\{v_{2k+1},\ldots,v_{3k-1}\}$.
Let $\mathcal{H}'$ be the induced subcomplex of $\mathcal{H}$ on $V_1\cup V_2$.
Thus, $\mathcal{H}'$ is a $k$-partite $(k-1)$-complex on $2k-2$ vertices.
Let $\mathcal{G}=\mathcal{J}\cup G_{\mathcal{J}}$.

Let $\mathcal{H}_{\mathcal{G}}'$ be the set of labelled partition-respecting copies of $\mathcal{H}'$ in $\mathcal{G}$.
It follows that
\begin{equation}
|\mathcal{H}_{\mathcal{G}}'|=(1\pm\varepsilon_{k+1})|\mathcal{J}_{(Y_1,\ldots,Y_{k-1})}|^2,
\end{equation}
where the error term accounts for the fact that we do not count the intersecting pairs of $(k-1)$-tuples in $\mathcal{J}_{(Y_1,\ldots,Y_{k-1})}$.
Since $Y$ is an edge of $R$, any function $\phi: V(P)\rightarrow V(R)$ such that $\phi(X_i)\subseteq Y_i$ is a homomorphism.
By Lemma \ref{extension} with $\beta^2$ playing the role of $\beta$, we deduce that all but at most $\beta^2|\mathcal{H}_{\mathcal{G}}'|$ of labelled partition-respecting copies of $\mathcal{H}'$ in $\mathcal{G}$ extend to at least $cm^{3k+1}$ labelled partition-respecting copies of $\mathcal{H}$ in $\mathcal{G}$, since $c\ll d_2,\ldots,d_{k-1}$.
For each $e\in\mathcal{J}_{(Y_1,\ldots,Y_{k-1})}$, let $T(e)$ be the number of tuples $e'$ in $\mathcal{J}_{(Y_1,\ldots,Y_{k-1})}$ such that $e\cup e'$ can be extended to at least $cm^{3k+1}$ copies of $\mathcal{H}$ in $\mathcal{G}$,
We have
\begin{equation}
\label{8.2}
\sum_{e\in\mathcal{J}_{(Y_1,\ldots,Y_{k-1})}}T(e)\geq(1-2\beta^2)|\mathcal{J}_{(Y_1,\ldots,Y_{k-1})}|^2.
\end{equation}

Let $S\subseteq \mathcal{J}_{(Y_1,\ldots,Y_{k-1})}$ be the set of $(k-1)$-tuples $e$ which is not $(c,\nu)$-extensible leftwards to $Y$, that is $T(e)<\nu|\mathcal{J}_{(Y_1,\ldots,Y_{k-1})}|$.
Combining with (\ref{8.2}) and $\beta\ll\nu$, we have
\begin{equation}
\nonumber
\sum_{e\in\mathcal{J}_{(Y_1,\ldots,Y_{k-1})}}T(e)\leq|S|\cdot\nu|\mathcal{J}_{(Y_1,\ldots,Y_{k-1})}|+(|\mathcal{J}_{(Y_1,\ldots,Y_{k-1})}|-|S|)|\mathcal{J}_{(Y_1,\ldots,Y_{k-1})}|,
\end{equation}
furthermore, we have
\begin{equation}
\nonumber
|S|\leq\frac{2\beta^2}{1-\nu}|\mathcal{J}_{(Y_1,\ldots,Y_{k-1})}|
\leq\frac{\beta}{2}|\mathcal{J}_{(Y_1,\ldots,Y_{k-1})}|.
\end{equation}
A symmetric fact shows that all but at most $\frac{\beta}{2}|\mathcal{J}_{(Y_1,\ldots,Y_{k-1})}|$ $(k-1)$-tuples in $\mathcal{J}_{(Y_1,\ldots,Y_{k-1})}$ are not $(c,\nu)$-extensible rightwards to $Y$.
Thus, all but at most $\beta|\mathcal{J}_{(Y_1,\ldots,Y_{k-1})}|$ pairs in $\mathcal{J}_{(Y_1,\ldots,Y_{k-1})}$ are not $(c,\nu)$-extensible both left and rightwards to $Y$.
\end{proof}

In Proposition \ref{8.1}, we know that most tuples in the complex induced by an edge of the reduced graph of a regular slice also extend to that edge.
The following lemma allows us to connect up two extensible paths using either very few or quite a lot of vertices.
\begin{lemma}
\label{lem8.1}
Let $k,r,m,t\in \mathbb{N}$, and $d_2,\ldots,d_{k+1},\varepsilon,\varepsilon_{k+1},c,\nu,\lambda$ be such that
\[
1/m\ll1/r,\varepsilon\ll c\ll\varepsilon_{k+1},d_2,\ldots,d_k,
\]
\[
\lambda\ll\nu\ll1/k,
\]
\[
\varepsilon_{k+1}\ll d_{k+1}.
\]
Let $\textbf{d}=(d_2,\ldots,d_{k+1})$ and let $\mathfrak{S}=(G,G_{\mathcal{J}},\mathcal{J},\mathcal{P},H)$ be a $(k,m,2t,\varepsilon,\varepsilon_{k+1},r,\textbf{d})$-regular setup where $\mathcal{P}$ has an initial partition of $[n]\cup V$ and $H$ is a $(1,k)$-graph on $[t]\cup V'$.
Suppose that $H_{W_i}=H[[t(i-1)/k,ti/k]\cup V']$ and $H_{W_i}$ is sequentially tightly connected for $i\in[k]$.
Let $P_1$, $P_2\subseteq G$ be $(c,\nu)$-extensible paths such that $P_1$ extends rightwards to $X$ and $P_2$ extends leftwards to $Y$.
Suppose that $P_1$ and $P_2$ are either identical or disjoint, let $W$ be a sequentially walk traversing each $H_{W_i}$ of length at most $t^{2k+1}$ that starts from $X$ and ends with $Y$.
Let $T$ be the joint connection set of $P_1$ and $P_2$.
Suppose that $T$ and $S\subseteq V(G)$ are $\lambda$-sparse in $\mathcal{P}$, $V(P_1)\cup V(P_2)\subseteq S$ and $T\cap S=\emptyset$, then

(1) there is a sequentially path $Q$ of length $4k-1+(\ell(W)+2)(k+1)$ in $G[V(\mathcal{P})]$ such that $P_1QP_2$ is a sequentially path, containing no vertices of $S$ and exactly $6k+2$ vertices of $T$,

(2) consider $\psi$ with $\varepsilon_{k+1}\ll\psi$, let $\textbf{w}$ be a fractional matching of size $\mu=\sum_{i\in[k]}\sum_{e\in E(H_{W_i})}\textbf{w}_i(e)$ $\geq5/m$ such that $\sum_{Z\in e,e\in H_{W_i}}\textbf{w}_i(e)\leq(1-2\lambda)/k$ for each $Z\in \mathcal{P}$.
There is a sequentially path $Q$ of length $\ell(W)+1$ mod $k$ in $G[V(\mathcal{P})]$ such that $P_1QP_2$ is a sequentially path, containing no vertices of $S$ and exactly $6k+2$ vertices of $T$.
Moreover, there is a set $U\subseteq V(\mathcal{P})$ of size at most $\psi mt$ such that $U\cup V(Q)$ has exactly 
$\lceil\sum_{i\in[k]}\sum_{Z\in e,e\in H_{W_i}}\textbf{w}_i(e)m\rceil+B$  vertices in each point cluster $Z$.
\end{lemma}
\begin{proof}
Let $X=(X_0,X_1,\ldots,X_k)$, since $P_1$ extends rightwards to $X$, thus there exists a target set $T_1\subseteq \mathcal{J}_{(X_2,\ldots,X_k)}$ of size $|T_1|\geq\nu|\mathcal{J}_{(X_2,\ldots,X_k)}|$ such that for every $(v_2,\ldots,v_k)\in T_1$, there are at least $cm^{3k+1}$ many $(3k+1)$-tuples $(c_1,\ldots,c_{2k},w_1,\ldots,w_k,v_1)$ with $c_i\in T\cap X_0$ for $i\in[2k]$, $w_i\in T\cap X_i$ for $i\in[k]$ and $v_1\in T\cap X_1$ such that $((c_1,\ldots,c_{2k}), P_1(w_1,\ldots,w_k,v_1,\ldots,v_k))$ is a sequentially path.
Let $Y=(Y_0,Y_1,\ldots,Y_k)$, $P_2$ extends leftwards to $Y$ with target set $T_2\subseteq \mathcal{J}_{(Y_2,\ldots,Y_k)}$.

For each $Z\in \mathcal{P}$, let $Z'\subseteq Z\setminus (S\cup T)$ of size $m'=(1-2\lambda)m$ since $S$ and $T$ are $\lambda$-sparse.
Let $\mathcal{P}'=\{Z'\}_{Z\in \mathcal{P}}$, $G'=G[V(\mathcal{P}')]$ and $\mathcal{J}'=\mathcal{J}[V(\mathcal{P}')]$.
By lemma \ref{restriction}, $\mathfrak{S}':=(G',G_{\mathcal{J}}',\mathcal{J}',\mathcal{P}',H)$ is a $(k,m',2t,\sqrt{\varepsilon},\sqrt{\varepsilon_{k+1}},r,\textbf{d})$-regular setup.

For (2), let $\mu'=\mu/(1-2\lambda)$ be the scaled size of $\textbf{w}$ and $B\in \mathbb{N}$ such that $1/B\ll1/r,\varepsilon$.
Let $\ell$ be the largest integer divisible by $k$ with $4k\leq\ell\leq(1-\psi/4)\mu'm'k$.
Note that such an $\ell$ exists since $(1-\psi/4)\mu'm'\geq4$, where the latter inequality follows from $\mu\geq5/m$.
Applying Lemma \ref{coro} with $G', \mathcal{J}', W, \ell, \textbf{w}, \mu'$ and $T_1,T_2$, we obtain a sequentially path $Q'$ whose initial $(k-1)$-tuple belongs to $T_1$ and whose terminal $(k-1)$-tuple belongs to $T_2$.
Furthermore, $Q'$ has length $\ell-1+\ell(W)(k+1)$ and uses at most $\sum_{i\in[k]}\mu_i(Z)m+B$ vertices from any point cluster $Z$
where $\mu_i(Z)=\sum_{Z\in e,e\in H_{W_i}}\textbf{w}_i(e)$ and $B\ll\psi\mu mk$.
Note that $\ell\geq(1-\psi/4)\mu km-k$, it follows that
\begin{equation}
\nonumber
\begin{split}
&\sum_{Z\in V'}\sum_{i\in[k]}\mu_i(Z)m-\sum_{Z\in V'}|V(Q')\cap Z|\\
&\leq \mu km-(1-\frac{\psi}{4})\mu km+k+1-\ell(W)(k+1)\\
&\leq\frac{\psi}{4}\mu km+k+1\\
&\leq\frac{\psi}{4}(1-2\lambda)tm+k+1\leq\frac{\psi}{2} mt.
\end{split}
\end{equation}
Hence, there is a set $U\subseteq V(\mathcal{P})$ of size at most $\psi mt$ such that $U\cup V(Q')$ has $\lceil\sum_{i\in[k]}\mu_i(Z)m\rceil+B$ vertices from any point cluster $Z\in V'$.

For (1), we can choose a path $Q'$ in the same way.
The only difference is  that in this case $\textbf{w}$ is a single edge of weight 1 and $\ell=4k$.
Hence, $Q'$ is a path of length $4k-1+\ell(W)(k+1)$.

Finally, we use the above extensible paths to choose $c_1,\ldots,c_{k+1},w_1,\ldots,w_k,v_1$ and $f_1,\ldots,f_{k+1},$
$v_k',w_1',\ldots,w_k'$ in $T$ such that for
\[
Q=((c_1,\ldots,c_{k+1})C(Q')(f_1,\ldots,f_{k+1}),(w _1,\ldots,w_k,v_1)Q'(v_k',w_1',\ldots,w_k')),
\]
the concatenation $P_1QP_2$ is a sequentially path and $Q$ is disjoint from $S$, since $V(S)\cap T=\emptyset$ $T\cap V(Q')=\emptyset$.
It is obvious that the length of $Q$ in (1) is $4k-1+(\ell(W)+2)(k+1)$ and the length of $Q$ in (2) is $\ell(W)+1$ mod $k$.

\end{proof}

\begin{proposition}
\label{8.3}
Let $W$ be a sequentially walk in a $(1,k)$-graph $H$ on $[t]\cup V'$ which starts from $(1,k)$-tuple $X$ and ends with $(1,k)$-tuple $Y$ where $|V'|=t$.
There exists a sequentially walk $W'$ of length at most $kt^{k+1}$, which starts from $X$ and ends with $Y$.
Moreover, $\ell(W')=\ell(W)$ mod $k$.
\end{proposition}
\begin{proof}
Suppose that $\ell(W)=j$ mod $k$ for a $j\in[0,k-1]$.
Let $W'$ be a vertex-minimal sequentially tightly walk from $X$ to $Y$ of size $j$ mod $k$.
Our goal is to show that every $(1,k)$-tuple repeats at most $k$ times in $W'$.

Assume that $W'$ contains $k+1$ copies of the same $(1,k)$-tuple $Z$ and denote by $n_j$ the position in $W'$ where the $j$th repetition $Z$ begins.
It is obvious that $n_j-n_1\not\equiv 0$ mod $k$, otherwise it is contrary to the minimal of $W'$.
By the pigeonhole principle, there exist two indices $j, j'$ such that $n_j-n_1\equiv n_{j'}-n_1$ mod $k$ for $1\leq j<j'\leq k+1$.
That is, $n_j-n_{j'}\equiv0$ mod $k$.
We can also reduce the length of $W'$ by deleting the vertices between $n_j$ and $n_{j'}-1$, a contradiction.
\end{proof}

\begin{proposition}
\label{8.4}
Let $j,k,t\in \mathbb{N}$ with $j\in[k]$.
Let $W$ be a sequentially closed walk that is compatible with respect to an orientation $\overrightarrow{H}$ of a $(1,k)$-graph $H$ on $[t]\cup V'$ where $|V'|=t$.
Let $X_1$ and $X_2$ be consistent with $\overrightarrow{H}$.
There exists a sequentially walk $W'$ of length at most $kt^{k+1}$, which starts from $X_1$ and ends with $X_2$.
Moreover, if $W$ has length 1 ${\rm mod}\ k$, then $W'$ has length $j$ ${\rm mod}\ k$.
\end{proposition}
\begin{proof}
For the first part, by Proposition \ref{8.3}, it suffices to show that there is a sequentially walk starting from $X_1$ and ending with $X_2$.
Since $X_1$ is consistent with $\overrightarrow{H}$, there is a sequentially path $W_{X_1}$ of length at most $k-1$ from $X_1$ to $X_1'$ in $H$ where $X_1'$ is an oriented edge in $\overrightarrow{H}$ which is a cyclic shift of $X_1$.
Similarly, there is a sequentially path $W_{X_2}$ of length at most $k-1$ from $X_2$ to $X_2'$ in $H$ where $X_2'$ is an oriented edge in $\overrightarrow{H}$ which is a cyclic shift of $X_2$.
Since $W$ is compatible with respect to an orientation $\overrightarrow{H}$, there is a subwalk $W_{X_1'X_2'}\subseteq W$ starting from $X_1'$ and ending with $X_2'$, hence $(C(X_1)C(W_{X_1})C(W_{X_1'X_2'})C(W_{X_2})C(X_2),I(X_1)I(W_{X_1})I(W_{X_1'X_2'})I(W_{X_2})I(X_2))$ is the desired $W'$.

Note that we choose $W_{X_1'X_2'}$ such that $W'$ has length $j$ mod $k$ by extending $W_{X_1'X_2'}$ along the same $(1,k)$-tuple with copies of $W$, for an appropriate number of times.
This is possible since any number coprime to $k$ is a generator for the finite cyclic group $\mathbb{Z}/k\mathbb{Z}$.
\end{proof}
\begin{lemma}[Connecting Lemma]
\label{8.5}
Let $k,m,r,t\in \mathbb{N}$, $d_2,\ldots,d_{k+1},\varepsilon,\varepsilon_{k+1},p,\nu,\lambda,\zeta$ be such that
\[
1/m\ll1/r,\varepsilon\ll1/t,\zeta,\varepsilon_{k+1},d_2,\ldots,d_k,
\]
\[
\zeta\ll p\ll d_2,\ldots,d_k,
\]
\[
1/t\ll\varepsilon_{k+1}\ll d_{k+1},\nu\leq1/k,
\]
\[
\lambda\ll\nu\ll1/k.
\]
Let $\textbf{d}=(d_2,\ldots,d_{k+1})$ and let $(G,G_{\mathcal{J}},\mathcal{J},\mathcal{P},H)$ be a $(k,m,2t,\varepsilon,\varepsilon_{k+1},r,\textbf{d})$-regular setup with $H$ being sequentially tightly connected.
Let $\overrightarrow{H}$ be an orientation of $H$ with a compatible closed walk $W$.
Suppose that $\mathcal{C}$ is a collection of pairwise disjoint $(p,\nu)$-extensible paths consistent with $\overrightarrow{H}$ and with joint connection set $T$.
Assume that
\begin{itemize}
  \item[(1)] $|\mathcal{C}|\leq\zeta m$,
  \item[(2)] $V(\mathcal{C})$ is $\lambda$-sparse in $\mathcal{P}$,
  \item[(3)] $V(\mathcal{C})\cap T=\emptyset$.
\end{itemize}
Consider any two elements $P_1,P_2$ of $\mathcal{C}$, there is a sequentially path $P$ in $G$ such that
\begin{itemize}
  \item[(a)] $P$ connects every path of $\mathcal{C}$,
  \item[(b)] $P$ starts from $P_1$ and ends with $P_2$,
  \item[(c)] $V(P)\setminus V(\mathcal{C})\subseteq V(\mathcal{P})$,
  \item[(d)] $V(P)\setminus V(\mathcal{C})$ intersects in at most $10k^2\mathcal{C}_Z+t^{2t+3k+2}$ vertices with each cluster $Z\in \mathcal{P}$, where $\mathcal{C}_Z$ denotes the number of paths of $\mathcal{C}$ intersecting with $Z$.
\end{itemize}
\end{lemma}
\begin{proof}
Choose a set $T'$ from $V(G)$ by including each vertex of $V(\mathcal{P})$ independently at random with probability $p$.
By Proposition \ref{chernoff} and the union bound, we obtain that the set $T'$ is $(2p)$-sparse with probability $1-2t\exp(-\Omega(m))$.
By Proposition \ref{Mc}, we obtain that the set $T'$ is a connection set of a fixed $(p^{3k+2}/2,\nu)$-extensible path in $\mathcal{C}$ with probability $1-2m^{k-1}\exp(-\Omega(m))$.
Since $|\mathcal{C}|\leq\zeta m$, with positive probability, we get a set $T'$ satisfying all these properties.

Initiate $S=V(\mathcal{C})$.
While there are two paths $Q_1,Q_2\in \mathcal{A}$ such that the extension to the right of $Q_1$ equals to the left of $Q_2$, apply Lemma \ref{lem8.1} (1) with $\ell(W)=kp^{k+4}/2$ to obtain a path $Q$ of length $10k^2$ which avoids $S$ and has exactly $6k+2$ vertices in $T'$.
Add $V(Q)$ to $S$, replace $Q_1, Q_2$ with $Q$ in $\mathcal{C}$ and delete the $6k+2$ vertices used by $Q$ in $T'$.
Denote the set of paths after the procedure by $\mathcal{C}'$.

Note that the size of $S$ grows by at most $10k^2|\mathcal{C}|\leq10k^2\zeta m\leq\lambda m$, we delete at most $(6k+2)|\mathcal{C}|\leq(6k+2)\zeta m\leq p^{3k+2}m/4$ vertices from $T$ throughout this process since $\zeta\ll p$.
This implies that every path of $\mathcal{C}$ remains $(p^{3k+2}/4,\nu)$-extensible with connection set $T'$.
Hence the conditions of Lemma \ref{lem8.1} (1) are satisfied in every step and $\mathcal{C}'$ is well-defined.

Note that when the procedure ends, $\mathcal{C}'$ has size at most $t^{2t}$.
Moreover, the paths of $\mathcal{C}'$ inherit the property of being consistent with $\overrightarrow{H}$.
We continue by connecting up the paths of $\mathcal{C}'$ to the desired path $P$ along the orientation.
As the paths of $\mathcal{C}'$ are consistent with $\overrightarrow{H}$, the left and right extensions of each path in $\mathcal{C}'$ are contained in the walk $W$.
Since $W$ is compatible with $\overrightarrow{H}$, we can apply Proposition \ref{8.4} to obtain a sequentially walk in $H$ of length of at most $t^{2k+1}$ between the left and right end of each path in $\mathcal{C}'$.
Use Lemma \ref{restriction} and Lemma \ref{lem8.1} (1), we can connect up the paths of $\mathcal{C}'$ using at most $t^{2t+3k+2}$ further vertices of $V(\mathcal{P})$.

Thus, $P$ contains every path in $\mathcal{C}$ as a subpath and $V(P)\setminus V(\mathcal{C})\subseteq V(\mathcal{P})$.
Moreover, note that $V(\mathcal{C}')\setminus \mathcal{C}$ intersects in at most $10k^2\mathcal{C}_Z$ vertices for each $Z\in \mathcal{P}$, where $\mathcal{C}_Z$ denotes the number of paths of $\mathcal{C}$ that intersects with $Z$.
It is obvious that $P$ can start and end with any two paths of $\mathcal{C}$.
\end{proof}
\begin{proof}[Proof of Lemma \ref{almost cover}]
Let $P_1=P$.
Suppose that $P_1$ extends rightwards to $X$ and leftwards to $Y$, there exists a path $P_2$ of length $k-1$ which $(c,\nu)$-extends both leftwards and rightwards to $Y$ by Proposition \ref{8.1}.
Moreover, we can assume that $V(P_1)$ is disjoint from $V(P_2)$ and $T_2$, where $T_2$ is the connection set of $P_2$.
By Proposition \ref{chernoff} and Proposition \ref{Mc}, we can choose a $\lambda$-sparse vertex set $T'$ such that $P_1$, $P_2$ are $(c^{3k+2}/2,\nu)$-extensible paths with connection set $T'$.

Firstly, let $S_1=V(P_1)\cup V(P_2)$, and we choose $\kappa$ such that $\lambda\ll\kappa\ll\gamma$.
For each $Z\in \mathcal{P}$, we can select a subset $Z'$ of $Z$ of size $m'=\kappa m$ such that $Z\cap S_1\subseteq Z'$ since $S_1$ is $2\lambda$-sparse, $1/m\ll1/t\ll\alpha\ll\lambda$ and $2\lambda\ll\kappa$.
Let $\mathcal{P}'=\{Z'\}_{Z\in \mathcal{P}}$, $V(\mathcal{P}')=\bigcup_{Z\in \mathcal{P}}Z'$,
$G'=G[V(\mathcal{P}')]$, $G_{\mathcal{J}'}'=G_{\mathcal{J}}[V(\mathcal{P}')]$ be the corresponding induced subgraphs and $\mathcal{J}'=\mathcal{J}[V(\mathcal{P}')]$ be the induced subcomplex.
By Lemma \ref{restriction}, $\mathfrak{S}'=(G',G_{\mathcal{J}'}',\mathcal{J}',\mathcal{P}',H)$ is a $(k,m',2t,\sqrt{\varepsilon},\sqrt{\varepsilon_{k+1}},r,d_2,\ldots,d_{k+1})$-regular setup.

Now we define a fractional matching that complements the discrepancy of $S_1$ in the clusters of $\mathcal{P}$.
Consider $\textbf{b}_i\in \mathbb{R}^{V(H_{W_i})}$ by setting $\textbf{b}_i(Z')=|Z'\setminus S_1|/|Z'|$ for every $Z\in V(H_{W_i})$.
Recall that $|S_1\cap Z|\leq2\lambda m$, $|Z'|=\kappa m$ and $\lambda\ll\kappa,\gamma$.
It follows that
\[
1-\gamma\leq1-\frac{2\lambda}{\kappa}\leq1-\frac{|S_1|}{|Z'|}\leq \textbf{b}_i\leq1.
\]
Since $H_{W_i}$ is $\gamma$-robustly matchable, there is a fractional matching $\textbf{w}_i$ such that $\sum_{Z\in e,e\in H_{W_i}}\textbf{w}_i(e)=\textbf{b}_i(Z')/k$ for every cluster $Z'\in \mathcal{P}'$ of $H_{W_i}$ where $i\in[k]$.
Consider $\psi>0$ with $\varepsilon_{k+1}\ll\psi\ll\alpha$, there exists a sequentially path $Q_1$ in $G'$ such that $P_2Q_1P_1$ is a sequentially path in $G$ which contains no vertices of $S_1$ and $4k+2$ vertices of $T'$ by Lemma \ref{lem8.1}.
Moreover, there is a set $U\subseteq V(\mathcal{P})$ of size at most $\psi mt$ such that $U\cup V(Q_1)$ has
$\lceil\sum_{i\in[k]}\sum_{Z\in e,e\in H_{W_i}}\textbf{w}_i(e)\kappa m\rceil+B$  vertices in each point cluster $Z$.
In other words, $V(P_2Q_1P_1)\cup U$ has $\kappa m+B$ vertices in each point cluster of $V(H)$ and uses $(\kappa m+B)(1-\alpha)t$ vertices of $V$ since $|V(L_H(i))|\geq(1-\alpha)t$ for $i\in[t]$.

We now choose the second path $Q_2$.
Note that $P_2Q_1P_1$ has right extension $X$ and left extension $Y$, which are consistent with $\overrightarrow{H}$.
Since $W$ is compatible with $\overrightarrow{H}$, we can apply Proposition \ref{8.4} to obtain a sequentially walk $W'$ in $H$  of length $p\leq t^{2k+1}$ starting from $X$ and ending with $Y$.
Moreover, since $W$ has length coprime to $k$, we can choose $W'$ such that
\[
p+1=|V(G)\setminus V(P_2Q_1P_1)|\ {\rm mod}\ k.
\]
Let $S_2=V(P_2Q_1P_1)$ and $T''=T'\setminus S_2$.
Define $\textbf{c}_i\in \mathbb{R}^{V(H_{W_i})}$ by setting $\textbf{c}_i(Z)=(m-|Z\cap S_2|)/m$ for every $Z\in V(H_{W_i})$.
Note that $1-\gamma\leq1-\kappa-\psi\leq \textbf{c}_i\leq1$.
Since $H_{W_i}$ is robustly matchable, there is a fractional matching $\textbf{z}_i$ such that $\sum_{Z\in e,e\in H_{W_i}}\textbf{z}_i(e)=\textbf{c}_i(Z)/k$ for every $Z\in \mathcal{P}$ of $H_{W_i}$.
By Lemma \ref{lem8.1}, there exists a sequentially path $Q_2$ in $G$ of length $p+1$ mod $k$ which contains no vertices of $S_2$ and $4k+2$ vertices of $T''$ such that $P_2Q_1P_1Q_2$ is a sequentially cycle.
Besides, there is a set $U'\subseteq V(\mathcal{P})$ of size at most $\psi mt$ such that $U'\cup V(Q_2)$ has  
$\lceil\sum_{i\in[k]}\sum_{Z\in e,e\in H_{W_i}}\textbf{z}_i(e)m\rceil+B$  vertices in each point cluster $Z$.
Thus, $U'\cup V(Q_2)$ uses at least $\left((1-\kappa)m-B+B\right)(1-\alpha)t=(1-\kappa)m(1-\alpha)t$ vertices of $V$.
Denote the set of uncovered vertices in all clusters of $\mathcal{P}$ by $M$.

Note that $P_2Q_1P_1Q_2$ contains all vertices of $V(G)$ but $M$, $U$ and $U'$.
We know that $|M|\leq \alpha mt$, $|U|\leq\psi mt,|U'|\leq\psi mt$.
Thus $P_2Q_1P_1Q_2$ covers all but at most $\alpha mt+2\psi mt\leq3\alpha n\leq\eta n$ vertices.
Since the length of $Q_2$ is $p+1$ mod $k$, it follows that $|V\setminus V(P_2Q_1P_1Q_2)|$ is divisible by $k$.
\end{proof}

\section{Absorption}
We will give the proof of Lemma \ref{absorption lemma} in this section.
The method can be sketched as follows.
We define absorbing gadget to absorb a set $T$ of $k$ vertices and a set $O$ of $k$ colors.
For each $(T,O)$, the absorbing gadgets are numerous.
Based on the above properties, we can choose a small family of vertex-disjoint gadgets such that for every $(T,O)$, there are many absorbing gadgets.
Such a family is obtained by probabilistic method.
Connecting all these gadgets yields the desired absorbing path.

This section can be organised as follows.
In subsection~\ref{technical tools}, we attach vertices to regular complexes since the gadgets we need should be well-integrated in regular setups.
In section~\ref{absorbing gadget}, we count the number of absorbing gadgets for each $(T,O)$.
In section~\ref{absorbing lemma}, we select a well-behaved family of  absorbing gadgets, which is used to absorb a small number of arbitrary sets of $k$ vertices and $k$ colors.
\subsection{Technical Tools}
\label{technical tools}
In this part, we will obtain some results to help us attach vertices to regular complexes.
Let $H$ be a $(1,k)$-graph with vertex set $[n]\cup V$, $\mathcal{J}$ be a regular slice with cluster set $\mathcal{P}$.
Given a $(0,k-1)$-subset $X\subseteq \mathcal{P}$, $\mathcal{J}_X$ is an $|X|$-partite $|X|$-graph containing all edges of $|X|$-level of $\mathcal{J}$.
For any $v\in V$, $\delta>0$ and any color cluster $C$, let
\[
N_{\mathcal{J}}((v,C),\delta)=\{X\subseteq \mathcal{P}:|X|=k-1, {\rm for\ any\ }c\in C, |N_{H}((v,c);\mathcal{J}_X)|>\delta|\mathcal{J}_X|\},
\]
\begin{lemma}
\label{9.1}
Let $k,r,m,t\in \mathbb{N}$ and $d_2,\ldots,d_{k+1},\varepsilon,\varepsilon_{k+1},\mu,\delta$ be such that
\[
1/m\ll1/r,\varepsilon\ll\varepsilon_{k+1},d_2,\ldots,d_k,
\]
\[
\varepsilon_{k+1}\ll d_{k+1}\leq1/k,
\]
and
\[
\varepsilon_{k+1}\ll\mu\ll\delta.
\]
Let $\textbf{d}=(d_2,\ldots,d_{k+1})$ and let $(H,H_{\mathcal{J}},\mathcal{J},\mathcal{P},R)$ be a representative $(k,m,2t,\varepsilon,\varepsilon_{k+1},r,\textbf{d})$-regular setup.
Suppose that $H$ has minimum relative $(1,1)$-degree at least $\delta+\mu$ with vertex set $[n]\cup V$.
Then for any $v\in V$ and any color cluster $C$, we have
\[
|N_{\mathcal{J}}((v,C),\frac{\mu}{3})|\geq(\delta+\frac{\mu}{4})\binom{t}{k-1}.
\]
For any $c\in[n]$ and any point cluster $Z$, we have
\[
|N_{\mathcal{J}}((c,Z),\frac{\mu}{3})|\geq(\delta+\frac{\mu}{4})\binom{t}{k-1}.
\]
\begin{proof}
Let $v\in V$ and $c\in C$ be arbitrary.
The minimum relative degree condition implies that $\overline{\deg}_H(v,c)\geq\delta+\mu$.
Since the regular setup is representative and $\varepsilon_{k+1}\ll\mu$, we have $|\overline{\deg}_H(v,c)-\overline{\deg}_H((v,c);\mathcal{J})|<\varepsilon_{k+1}$ and
\[
\deg_H((v,c),\mathcal{J}^{(k-1)})\geq(\delta+\mu-\varepsilon_{k+1})|\mathcal{J}^{(k-1)}|\geq(\delta+\frac{2}{3}\mu)|\mathcal{J}^{(k-1)}|.
\]

For any $(0,k-1)$-subset $X$ of $\mathcal{P}$, $\mathcal{J}_X$ corresponds to the $(k-1)$-edges of $\mathcal{J}^{(k-1)}$ which are $X$-partite.
Define $d_X=\prod_{i=2}^{k-1}d_i^{\binom{k-1}{i}}$.
By Lemma \ref{7.7}, we have $|\mathcal{J}_X|=(1\pm\varepsilon_{k+1})d_Xm^{k-1}$.
By summing over all the $(0,k-1)$-subsets of $\mathcal{P}$, we have
\[
|\mathcal{J}^{(k-1)}|\geq(1-\varepsilon_{k+1})\binom{t}{k-1}d_Xm^{k-1}.
\]
Moreover, let $X$ range over all $(0,k-1)$-subsets of $\mathcal{P}$, we have
\[
\sum_X|N_H((v,c);\mathcal{J}_X)|=\deg_H((v,c);\mathcal{J}^{(k-1)})\geq(\delta+\frac{2}{3}\mu)|\mathcal{J}^{(k-1)}|.
\]
Finally, we obtain
\begin{equation}
\nonumber
\begin{split}
&(\delta+\frac{2}{3}\mu)|\mathcal{J}^{(k-1)}|\\
&\leq\sum_X|N_H((v,c);\mathcal{J}_X)|\leq\sum_{X\in N_{\mathcal{J}}((v,c),\mu/3)}|\mathcal{J}_X|+\sum_{X\notin N_{\mathcal{J}}((v,c),\mu/3)}\frac{\mu}{3}|\mathcal{J}_X|\\
&\leq\left(|N_{\mathcal{J}}((v,c),\mu/3)|+\frac{\mu}{3}\left(\binom{t}{k-1}-|N_{\mathcal{J}}((v,c),\mu/3)|\right)\right)(1+\varepsilon_{k+1})d_Xm^{k-1}\\
&\leq\left((1-\frac{\mu}{3})|N_{\mathcal{J}}((v,c),\mu/3)|+\frac{\mu}{3}\binom{t}{k-1}\right)\frac{1+\varepsilon_{k+1}}{1-\varepsilon_{k+1}}\frac{|\mathcal{J}^{(k-1)}|}{\binom{t}{k-1}}\\
&\leq\left(|N_{\mathcal{J}}((v,c),\mu/3)|+\frac{\mu}{3}\binom{t}{k-1}\right)(1+2\varepsilon_{k+1})\frac{|\mathcal{J}^{(k-1)}|}{\binom{t}{k-1}}.
\end{split}
\end{equation}
Thus, for any $v\in V$ and $c\in C$, we have
\[
|N_{\mathcal{J}}((v,c),\mu/3)|\geq(\delta+\frac{\mu}{4})\binom{t}{k-1},
\]
and by definition, the following holds for any $v\in V$ and color cluster $C$,
\[
|N_{\mathcal{J}}((v,C),\mu/3)|\geq(\delta+\frac{\mu}{4})\binom{t}{k-1}.
\]

Similarly, we can obtain the following result holds for any $c\in[n]$ and  point cluster $Z$,
\[
|N_{\mathcal{J}}((c,Z),\mu/3)|\geq(\delta+\frac{\mu}{4})\binom{t}{k-1}.
\]
\end{proof}
\end{lemma}

\begin{lemma}
\label{9.2}
Let $k,r,m,t\in \mathbb{N}$ and $d_2,\ldots,d_{k+1},\varepsilon,\varepsilon_{k+1},\mu,\lambda$ be such that
\[
1/m\ll1/r,\varepsilon\ll\varepsilon_{k+1},d_2,\ldots,d_k,
\]
\[
\varepsilon_{k+1}\ll d_{k+1}\leq1/k,
\]
and
\[
\varepsilon_{k+1}\ll\lambda\ll\mu.
\]
Let $\textbf{d}=(d_2,\ldots,d_{k+1})$ and let $(H,H_{\mathcal{J}},\mathcal{J},\mathcal{P},R)$ be a $(k,m,2t,\varepsilon,\varepsilon_{k+1},r,\textbf{d})$-regular setup.
Let $T\subseteq V(H)$ such that $|Z_1\cap T|=|Z_2\cap T|\leq\lambda m$ for every $Z_1,Z_2\in \mathcal{P}$.
Let $Z'=Z\setminus T$ for each $Z\in\mathcal{P}$, and let $\mathcal{J}'=\mathcal{J}[\bigcup Z']$ be the induced subcomplex.
For every $v\in V$ and color cluster $C$, we have
\[
|N_{\mathcal{J}}((v,C),2\mu)|\leq|N_{\mathcal{J}'}((v,C),\mu)|,
\]
and for every $c\in [n]$ and point cluster $Z$, we have
\[
|N_{\mathcal{J}}((c,Z),2\mu)|\leq|N_{\mathcal{J}'}((c,Z),\mu)|,
\]
\end{lemma}
\begin{proof}
For any $v\in V$, color cluster $C$ and a $(0,k-1)$-set $X\in N_{\mathcal{J}}((v,C),2\mu)$.
By the definition, we have $|N_{H}((v,c);\mathcal{J}_X)|>2\mu|\mathcal{J}_X|$ for any $c\in C$.
Let $X=\{X_1,\ldots,X_{k-1}\}$ and $X'=\{X_1',\ldots,X_{k-1}'\}$ be the corresponding clusters in the complex $\mathcal{J}'$.
Our goal is to prove that $X'\in N_{\mathcal{J}'}((v,C),\mu)$.

Let $\varepsilon\ll\beta\ll\varepsilon_{k+1}$ and $d_X=\prod_{i=2}^{k-1}d_i^{\binom{k-1}{i}}$.
By Lemma \ref{7.7}, we have
\[
|\mathcal{J}_X|=(1\pm\beta)d_Xm^{k-1}
\]
and
\[
|N_H((v,c);\mathcal{J}_X)|>2\mu|\mathcal{J}_X|\geq2\mu(1-\beta)d_Xm^{k-1}.
\]

Let $m'=|X_1\setminus T|$, we have $|Z'|=m'$ for each $Z\in \mathcal{P}$, note that $m'\geq(1-\lambda)m$.
By Lemma \ref{restriction}, $\mathcal{J}'$ is a $(\cdot,\cdot,\sqrt{\varepsilon},\sqrt{\varepsilon_{k+1}},r)$-regular slice.
By Lemma \ref{7.7}, we have
\[
(1+\beta)d_X(m')^{k-1}\geq|\mathcal{J}'_{X'}|\geq(1-\beta)d_X(m')^{k-1}\geq(1-\beta)(1-\lambda)^{k-1}d_Xm^{k-1}.
\]

Since $\beta\ll\varepsilon_{k+1}\ll\lambda\ll\mu$, we have
\begin{equation}
\nonumber
\begin{split}
|N_H((v,c);\mathcal{J}'_{X'})|&\geq|N_H((v,c);\mathcal{J}_X)|-(|\mathcal{J}_X|-|\mathcal{J}'_{X'}|)\\
&\geq(1-\beta)(2\mu-(1-(1-\lambda)^{k-1}))d_Xm^{k-1}\\
&\geq\mu(1+\beta)d_Xm^{k-1}\geq\mu|\mathcal{J}'_{X'}|.
\end{split}
\end{equation}
Thus, we obtain that $X\in N_{\mathcal{J}'}((v,C),\mu)$.

Similarly, for every $c\in [n]$ and point cluster $Z$, we have
\[
|N_{\mathcal{J}}((c,Z),2\mu)|\leq|N_{\mathcal{J}'}((c,Z),\mu)|.
\]
\end{proof}

In a $(k+1)$-uniform sequentially cycle, the link graph of a point corresponds to a $k$-uniform sequentially path.
Thus, we will look for sequentially paths in the neighbors of vertices inside a regular complex.
The following lemma states that by looking at a $\mu$-fraction of $ (1,k-1)$-edges of a regular complex, we will find lots of sequentially paths.
\begin{lemma}
\label{9.3}
Let $1/m\ll\varepsilon\ll d_2,\ldots,d_k,1/k,\mu$ and $k\geq3$.
Suppose that $\mathcal{J}$ is a $(\cdot,\cdot,\varepsilon)$-equitable complex with density vector $\textbf{d}=(d_2,\ldots,d_k)$ and ground partition $\mathcal{P}$, the size of each vertex class is $m$.
Let $W=\{W_0,W_1,\ldots,W_{k-1}\}\subseteq \mathcal{P}$.
Let $S\subseteq\mathcal{J}_W$ be with size at least $\mu|\mathcal{J}_W|$ and $Q$ be a $k$-uniform sequentially path $((c_1,\ldots,c_k),(v_1,\ldots,v_{2k-2}))$ with vertex classes $\{X_0,X_1,\ldots,X_{k-1}\}$ such that $v_i,v_{i+k-1}\in X_i$ for $i\in[k-1]$ and $c_j\in X_0$ for $j\in[k]$.
Let $\mathcal{Q}$ be the down-closed $k$-complex generated by $Q$ and $\mathcal{Q}_S\subseteq\mathcal{Q}_{\mathcal{J}}$ be the copies of $\mathcal{Q}$ whose edges in the $k$-th level are in $S$.
We have
\[
|\mathcal{Q}_S|\geq\frac{1}{2}\left(\frac{\mu}{8k}\right)^{k+1}|\mathcal{Q}_{\mathcal{J}}|.
\]
\end{lemma}
\begin{proof}
The proof consists of three steps.
Firstly, we use the dense version of the counting and extension lemma to count the number of various hypergraphs in $\mathcal{J}$.
Secondly, we remove some $(1,k-1)$-tuples without good properties.
Finally, we use an iterative procedure to return sequentially paths using good $(1,k-1)$-tuples, as desired.

Firstly, let $\beta$ be such that $\varepsilon\ll\beta\ll d_2,\ldots,d_k,1/k,\mu$.
Define
\[
d_a=\prod_{i=2}^{k-2}d_i^{\binom{k-2}{i}}, d_b=\prod_{i=2}^{k-2}d_i^{\binom{k}{i}-\binom{k-2}{i}}\cdot\prod_{i=k-1}^{k}d_i^{\binom{k}{i}}.
\]
Let $W'=W\setminus\{W_0,W_{k-1}\}$.
By Lemmas \ref{dense counting} and \ref{7.7}, we have
\begin{equation}
\label{in1}
|\mathcal{J}_W|=(1\pm\beta)d_ad_bm^k,
\end{equation}
\begin{equation}
\label{in2}
|\mathcal{J}_{W'}|=(1\pm\beta)d_am^{k-2},
\nonumber
\end{equation}
\begin{equation}
\label{in3}
|\mathcal{Q}_{\mathcal{J}}|=(1\pm\beta)d_ad_b^km^{3k-2}.
\nonumber
\end{equation}
Since $S\subseteq\mathcal{J}_W$ with $|S|\geq\mu|\mathcal{J}_W|$, with (\ref{in1}), we have
\begin{equation}
|S|\geq(1-\beta)\mu d_ad_bm^k.
\nonumber
\end{equation}
Let $B_{W'}\subseteq\mathcal{J}_{W'}$ be the $(k-2)$-edges which are not extensible to $(1\pm\beta)d_bm^2$ copies of a $k$-edge in $\mathcal{J}_W$.
By Lemma \ref{dense extension}, we have
\begin{equation}
|B_{W'}|\leq\beta|\mathcal{J}_{W'}|.
\nonumber
\end{equation}

Secondly, we delete from $S$ the edges which contain a $(k-2)$-set from $B_{W'}$ to obtain $S'$, the number of edges deleted is at most
\begin{equation}
|B_{W'}|m^2\leq\beta|\mathcal{J}_{W'}|m^2\leq\beta(1+\beta)d_am^k\leq|S|/3,
\nonumber
\end{equation}
since $\beta\ll\mu,d_2,\ldots,d_k$.
Thus, we have $|S'|\geq2|S|/3$.
Furthermore, if there is any partite $(k-2)$-set $T$ in $\mathcal{J}$ which lies in less than $\mu d_bm^2/(4k)$ edges of $S'$, then we delete all edges in $S'$ containing $T$ to obtain $S''$ and iterate this until no further deletions are possible.
Note that the number of partite $(k-2)$-sets supported in the clusters of  $W\setminus\{W_0\}$ is $(k-1)(1\pm\beta)d_am^{k-2}$.
Thus the number of edges deleted is at most
\[
(k-1)(1+\beta)d_am^{k-2}\frac{\mu d_bm^2}{4k}\leq(1+\beta)\frac{\mu d_ad_bm^k}{4}\leq\frac{|S|}{3}.
\]
Thus, $|S''|\geq |S|/3$.
Each partite $(k-2)$-set in $W_1,\ldots,W_{k-1}$ is either contained in zero edges of $S''$ or in at least $\mu d_bm^2/(4k)$ edges in $S''$.

Finally, we use the properties of $S''$ to construct many labelled partition-respecting paths in $\mathcal{Q}_S$.

\textbf{Step 1.} Select $T=\{x_1,\ldots,x_{k-2}\}\in \mathcal{J}_{W'}$ which is contained in at least $\mu d_bm^2/4$ edges in $S''$.

\textbf{Step 2.} Choose $(c_1,x_{k-1})$ such that $\{c_1,x_1,x_2,\ldots,x_{k-1}\}\in S''$ and $c_1,x_{k-1}$ are not in $T$.

\textbf{Step 3.} For $i\in[k,2k-2]$, choose $(c_{i-k+2},x_i)$ such that $\{c_{i-k+2},x_{i-k+2},\ldots,x_i\}\in S''$ and $c_{i-k+2},x_i$ are not used before.

This constructs a sequentially path $\mathcal{Q}_S$ on $3k-2$ vertices such that each edge in the $k$-th level is in $S''$, thus in $S$.
Next, we count the size of $\mathcal{Q}_S$.

In Step 1, let $G\subseteq\mathcal{J}_{W'}$ be the set of $(k-2)$-sets which are contained in less than $\mu d_bm^2/4$ edges in $S''$, we have
\[
\frac{|S|}{3}\leq|S''|=\sum_{T\in\mathcal{J}_{W'}}\deg_{S''}(T)\leq |G|\frac{\mu}{4}d_bm^2+(|\mathcal{J}_{W'}|-|G|)d_bm^2(1+\beta),
\]
it gives that $|G|\leq(1-\beta)(1-\mu/12)d_am^{k-2}$,
thus, the choices for $T$ is at least $|\mathcal{J}_{W'}|-|G|\geq\mu/13d_am^{k-2}$.
In Step 2, we have at least $\mu d_bm^2/4$ choices for $(c_1,x_{k-1})$.
In Step 3, $\{x_{i-k+2},\ldots,x_{i-1}\}$ is a $(k-2)$-set contained in $S''$, by the construction of $S''$, there are at least $\mu d_bm^2/(4k)$ choices for $(c_{i-k+2},x_i)$, furthermore, at least $\mu d_bm^2/(8k)$ are different from the previous choices.

Thus, the number of paths in $\mathcal{Q}_S$ is at least
\[
\left(\frac{\mu}{13}d_am^{k-2}\right)\left(\frac{\mu}{4}d_bm^2\right)\left(\frac{\mu}{8k}d_bm^2\right)^{k-1}\geq(\frac{\mu}{8k})^{k+1}d_ad_b^km^{3k-2}\geq\frac{1}{2}(\frac{\mu}{8k})^{k+1}|\mathcal{Q}_{\mathcal{J}}|,
\]
since $\beta\ll\mu,1/k$.
\end{proof}
\begin{lemma}
\label{color}
Let $1/m\ll\varepsilon\ll d_2,\ldots,d_k,1/k,\mu$ and $k\geq3$.
Suppose that $\mathcal{J}$ is a $(\cdot,\cdot,\varepsilon)$-equitable complex with density vector $\textbf{d}=(d_2,\ldots,d_k)$ and ground partition $\mathcal{P}$, the size of each vertex class is $m$.
Let $W=\{W_1,\ldots,W_{k-1},W_k\}\subseteq \mathcal{P}$.
Let $S\subseteq\mathcal{J}_W$ be with size at least $\mu|\mathcal{J}_W|$ and $Q$ be a $k$-uniform tight path $v_1,\ldots,v_{k-1},b,v_k,\ldots,v_{2k-2}$ with vertex classes $\{X_1,\ldots,X_{k-1},X_k\}$ such that $v_i,v_{i+k-1}\in X_i$ for $i\in[k-1]$ and $b\in X_k$.
Let $\mathcal{Q}$ be the down-closed $k$-complex generated by $Q$ and $\mathcal{Q}_S\subseteq\mathcal{Q}_{\mathcal{J}}$ be the copies of $\mathcal{Q}$ whose edges in the $k$-th level are in $S$.
We have
\[
|\mathcal{Q}_S|\geq\frac{1}{2}\left(\frac{\mu}{8k}\right)^{k+1}|\mathcal{Q}_{\mathcal{J}}|.
\]
\end{lemma}
\begin{proof}
The proof consists of three steps.
Firstly, we use the dense version of the counting and extension lemma to count the number of various hypergraphs in $\mathcal{J}$.
Secondly, we remove some $k$-tuples without good properties.
Finally, we use an iterative procedure to return a tight path using good $k$-tuples, as desired.

Firstly, let $\beta$ be such that $\varepsilon\ll\beta\ll d_2,\ldots,d_k,1/k,\mu$.
Define
\[
d_a=\prod_{i=2}^{k-1}d_i^{\binom{k-1}{i}}, d_b=\prod_{i=2}^{k}d_i^{\binom{k-1}{i-1}}.
\]
Let $W'=W\setminus\{W_k\}$.
By Lemma \ref{dense counting} and \ref{7.7}, we have
\begin{equation}
\label{in1}
|\mathcal{J}_W|=(1\pm\beta)d_ad_bm^k,
\end{equation}
\begin{equation}
\label{in2}
|\mathcal{J}_{W'}|=(1\pm\beta)d_am^{k-1},
\nonumber
\end{equation}
\begin{equation}
\label{in3}
|\mathcal{Q}_{\mathcal{J}}|=(1\pm\beta)d_ad_b^km^{2k-1}.
\nonumber
\end{equation}
Since $S\subseteq\mathcal{J}_W$ with $|S|\geq\mu|\mathcal{J}_W|$, with (\ref{in1}), we have
\begin{equation}
|S|\geq(1-\beta)\mu d_ad_bm^k.
\nonumber
\end{equation}
Let $B_{W'}\subseteq\mathcal{J}_{W'}$ be the $(k-1)$-edges which are not extensible to $(1\pm\beta)d_bm$ copies of a $k$-edge in $\mathcal{J}_W$.
By Lemma \ref{dense extension}, we have
\begin{equation}
|B_{W'}|\leq\beta|\mathcal{J}_{W'}|.
\nonumber
\end{equation}

Secondly, we delete from $S$ the edges which contain a $(k-1)$-set from $B_{W'}$ to obtain $S'$, the number of edges deleted is at most
\begin{equation}
|B_{W'}|m\leq\beta|\mathcal{J}_{W'}|m\leq\beta(1+\beta)d_am^k\leq|S|/3,
\nonumber
\end{equation}
since $\beta\ll\mu,d_2,\ldots,d_k$.
Thus, we have $|S'|\geq2|S|/3$.
Furthermore, if there is any partite $(k-1)$-set $T$ in $\mathcal{J}$ which lies in less than $\mu d_bm/(4k)$ edges of $S'$, then we delete all edges in $S'$ containing $T$ to obtain $S''$ and iterate this until no further deletions are possible.
Note that the number of partite $(k-1)$-sets supported in the clusters of  $W$ is $k(1\pm\beta)d_am^{k-1}$.
Thus the number of edges deleted is at most
\[
k(1+\beta)d_am^{k-1}\frac{\mu d_bm}{4k}\leq(1+\beta)\frac{\mu d_ad_bm^k}{4}\leq\frac{|S|}{3}.
\]
Thus, $|S''|\geq |S|/3$.
Each partite $(k-1)$-set in $W_1,\ldots,W_k$ is either contained in zero edges of $S''$ or in at least $\mu d_bm/(4k)$ edges in $S''$.

Finally, we use the properties of $S''$ to construct many labelled partition-respecting paths in $\mathcal{Q}_S$.

\textbf{Step 1.} Select $T=\{x_1,\ldots,x_{k-1}\}\in \mathcal{J}_{W'}$ which is contained in at least $\mu d_bm/4$ edges in $S''$.

\textbf{Step 2.} Choose $b$ such that $\{x_1,x_2,\ldots,x_{k-1},b\}\in S''$ and $b\notin T$.

\textbf{Step 3.} For $i\in[k,2k-2]$, choose $x_i$ such that $\{x_{i-k+2},\ldots,x_{k-1},b,x_k,\ldots,x_i\}\in S''$ and $x_i$ is not used before.

This constructs a sequentially path $\mathcal{Q}_S$ on $2k-1$ vertices such that each edge in the $k$-th level is in $S''$, thus in $S$.
Next, we count the size of $\mathcal{Q}_S$.

In Step 1, let $G\subseteq\mathcal{J}_{W'}$ be the set of $(k-1)$-sets which are contained in less than $\mu d_bm/4$ edges in $S''$, we have
\[
\frac{|S|}{3}\leq|S''|=\sum_{T\in\mathcal{J}_{W'}}\deg_{S''}(T)\leq |G|\frac{\mu}{4}d_bm+(|\mathcal{J}_{W'}|-|G|)d_bm(1+\beta),
\]
it gives that $|G|\leq(1-\beta)(1-\mu/12)d_am^{k-1}$,
thus, the choices for $T$ is at least $|\mathcal{J}_{W'}|-|G|\geq\mu/13d_am^{k-1}$.
In Step 2, we have at least $\mu d_bm/4$ choices for $b$.
In Step 3, $\{x_{i-k+2},\ldots,x_{k-1},b,x_k,\ldots,$ $x_{i-1}\}$ is a $(k-1)$-set contained in $S''$, by the construction of $S''$, there are at least $\mu d_bm/(4k)$ choices for $x_i$, furthermore, at least $\mu d_bm/(8k)$ are different from the previous choices.

Thus, the number of paths in $\mathcal{Q}_S$ is at least
\[
\left(\frac{\mu}{13}d_am^{k-1}\right)\left(\frac{\mu}{4}d_bm\right)\left(\frac{\mu}{8k}d_bm\right)^{k-1}\geq(\frac{\mu}{8k})^{k+1}d_ad_b^km^{2k-1}\geq\frac{1}{2}(\frac{\mu}{8k})^{k+1}|\mathcal{Q}_{\mathcal{J}}|,
\]
since $\beta\ll\mu,1/k$.
\end{proof}
\subsection{Absorbing Gadget}
\label{absorbing gadget}
Before we build the absorbing path, we need to define absorbing gadget, which is useful to absorb a particular set $T$ of $k$ vertices and a particular set $O$ of $k$ colors.
Next, we will show that for every $(T,O)$, there are numerous absorbing gadgets to absorb $(T,O)$.
\begin{definition}[Absorbing gadget]
\label{10.1absorbing}
Let $T=\{t_1,\ldots,t_k\}$ be a $k$-set of points of $G$ and $O=\{o_1,\ldots,o_k\}$ be a $k$-set of colors of $G$.
We say that $F\subseteq G$ is an absorbing gadget for $(T,O)$ if $F=F_1\cup F_2$ where $F_1=A\cup B\cup E\cup \bigcup_{i=1}^k(P_i\cup Q_i)\cup C\cup \bigcup_{i=1}^kC_k$ and $F_2=A'\cup B'\cup E'\cup \bigcup_{i=1}^k(P_i'\cup Q_i')\cup C'\cup \bigcup_{i=1}^kC_k'$ such that
\begin{itemize}
  \item[(1)] $A,B,E$,$P_1,Q_1,\ldots,P_k,Q_k$,$A',B',E'$, $P_1',Q_1',\ldots,P_k',Q_k'$ are pairwise disjoint and also disjoint from $T$. $C,C_1,\ldots, C_k,C',C_1',\ldots,C_k'$ are pairwise disjoint and also disjoint from $O$,
  \item[(2)] $C_i=(c_{i,1},\ldots,c_{i,k-1})$ and $C_i'=(c_{i,1}',\ldots,c_{i,k-1}')$ for $i\in[k]$,
  \item[(3)] $A,B,E,A',B',E'$ are $k$-tuples of points of $G$, $C$ and $C'$ are $(k+1)$-tuples of colors of $G$, $(C,AE)$, $(C',A'E')$ and $(C'(c_{1,1},\ldots,c_{k,1}),A'B'E')$ are sequentially paths,
  \item[(4)] for $B=(b_1,\ldots,b_k)$, each of $P_i,Q_i$ has $k-1$ vertices for $i\in[k]$, both $(C_i,P_ib_iQ_i)$ and $(\{o_i\}\cup C_i\setminus\{c_{i,1}\},P_ib_iQ_i)$ are sequentially paths of length $2k-1$ for $i\in[k]$.
   \item[(5)] for $B'=(b_1',\ldots,b_k')$, each of $P_i',Q_i'$ has $k-1$ vertices for $i\in[k]$, both $(C_i',P_i'b_i'Q_i')$ and $(C_i',P_i't_iQ_i')$ are sequentially paths of length $2k-1$ for $i\in[k]$.
\end{itemize}
\end{definition}

Note that an absorbing gadget $F$ spans $4k^2+2k$ points together with $2k^2+2k+2$ colors.
\begin{figure}[htb]
\begin{center}
\begin{tikzpicture}{center}
\filldraw [red] (-3,2) circle (1.2pt);
\filldraw [red] (-1.5,2) circle (1.2pt);
\filldraw [red] (0,2) circle (1.2pt);
\filldraw [red] (1.5,2) circle (1.2pt);
\node [red] at (-1,2.6){$C$};
\filldraw [black] (-4.5,0) circle (1.2pt);
\filldraw [black] (-3,0) circle (1.2pt);
\filldraw [black] (-1.5,0) circle (1.2pt);
\filldraw [black] (0,0) circle (1.2pt);
\filldraw [black] (1.5,0) circle (1.2pt);
\filldraw [black] (3,0) circle (1.2pt);
\node at (-4.5,-0.4){$a_1$};
\node at (-3,-0.4){$a_2$};
\node at (-1.5,-0.4){$a_3$};
\node at (0,-0.4){$e_1$};
\node at (1.5,-0.4){$e_2$};
\node at (3,-0.4){$e_3$};
\draw[rounded corners=0.3cm,line width =1pt] (-4.8,-0.15)--(-1.2,-0.15)--(-3,2.3)--cycle;
\draw[rounded corners=0.3cm,line width =1pt] (-3.3,-0.15)--(0.3,-0.15)--(-1.5,2.3)--cycle;
\draw[rounded corners=0.3cm,line width =1pt] (-1.8,-0.15)--(1.8,-0.15)--(0,2.3)--cycle;
\draw[rounded corners=0.3cm,line width =1pt] (-0.3,-0.15)--(3.3,-0.15)--(1.5,2.3)--cycle;
\filldraw [black] (-3.23,-1.96) circle (1pt);
\filldraw [black] (-0.97,-1.96) circle (1pt);
\filldraw [black] (1.2,-1.96) circle (1pt);
\node at (-3.2,-1.65){$b_1$};
\node at (-1,-1.65){$b_2$};
\node at (1.2,-1.65){$b_3$};
\filldraw [black] (-4,-3) circle (1.2pt);
\filldraw [black] (-3.5,-3.1) circle (1.2pt);
\filldraw [black] (-3,-3.1) circle (1.2pt);
\filldraw [black] (-2.5,-3) circle (1.2pt);
\filldraw [blue] (-3.75,-3.05) circle (1.2pt);
\filldraw [blue] (-3.25,-3.1) circle (1.2pt);
\filldraw [blue] (-2.75,-3.05) circle (1.2pt);
\node[blue] at (-3.25,-3.6){$C_1$};
\draw[rounded corners=0.2cm,line width =1pt] (-4.2,-3.1)--(-3.35,-3.25)--(-3.15,-1.75)--cycle;
\draw[rounded corners=0.2cm,line width =1pt] (-3.7,-3.2)--(-2.8,-3.2)--(-3.22,-1.75)--cycle;
\draw[rounded corners=0.2cm,line width =1pt] (-3.15,-3.25)--(-2.25,-3.05)--(-3.27,-1.75)--cycle;
\filldraw [black] (-0.7,-3.1) circle (1.2pt);
\filldraw [black] (-0.2,-3) circle (1.2pt);
\filldraw [black] (-1.2,-3.1) circle (1.2pt);
\filldraw [black] (-1.7,-3) circle (1.2pt);
\filldraw [orange] (-1.45,-3.05) circle (1.2pt);
\filldraw [orange] (-0.95,-3.1) circle (1.2pt);
\filldraw [orange] (-0.45,-3.05) circle (1.2pt);
\node[orange] at (-0.95,-3.6){$C_2$};
\draw[rounded corners=0.2cm,line width =1pt] (-1.9,-3.1)--(-1.05,-3.25)--(-0.9,-1.75)--cycle;
\draw[rounded corners=0.2cm,line width =1pt] (-1.4,-3.2)--(-0.5,-3.2)--(-0.95,-1.75)--cycle;
\draw[rounded corners=0.2cm,line width =1pt] (-0.85,-3.25)--(0.05,-3.05)--(-1.03,-1.765)--cycle;
\filldraw [black] (2,-3) circle (1.2pt);
\filldraw [black] (1.5,-3.1) circle (1.2pt);
\filldraw [black] (1,-3.1) circle (1.2pt);
\filldraw [black] (0.5,-3) circle (1.2pt);
\filldraw [gray] (0.75,-3.05) circle (1.2pt);
\filldraw [gray] (1.25,-3.1) circle (1.2pt);
\filldraw [gray] (1.75,-3.05) circle (1.2pt);
\node[gray] at (1.25,-3.6){$C_3$};
\node[black] at (-4.25,-2.5){$P_1$};
\node[black] at (-2.2,-2.5){$Q_1$};
\node[black] at (-1.8,-3.5){$P_2$};
\node[black] at (0,-3.5){$Q_2$};
\node[black] at (0.3,-2.5){$P_3$};
\node[black] at (2.2,-2.5){$Q_3$};
\node[black] at (4.75,-2.5){$P_1'$};
\node[black] at (6.75,-2.5){$Q_1'$};
\node[black] at (7.15,-3.5){$P_2'$};
\node[black] at (8.95,-3.5){$Q_2'$};
\node[black] at (9.25,-2.5){$P_3'$};
\node[black] at (11.15,-2.5){$Q_3'$};
\draw[rounded corners=0.2cm,line width =1pt] (0.3,-3.1)--(1.15,-3.25)--(1.25,-1.75)--cycle;
\draw[rounded corners=0.2cm,line width =1pt] (0.82,-3.2)--(1.7,-3.2)--(1.2,-1.75)--cycle;
\draw[rounded corners=0.2cm,line width =1pt] (1.35,-3.25)--(2.25,-3.05)--(1.15,-1.765)--cycle;

\filldraw [green] (6,2) circle (1.2pt);
\filldraw [green] (7.5,2) circle (1.2pt);
\filldraw [green] (9,2) circle (1.2pt);
\filldraw [green] (10.5,2) circle (1.2pt);
\node [green] at (8,2.6){$C'$};
\filldraw [black] (4.5,0) circle (1.2pt);
\filldraw [black] (6,0) circle (1.2pt);
\filldraw [black] (7.5,0) circle (1.2pt);
\filldraw [black] (9,0) circle (1.2pt);
\filldraw [black] (10.5,0) circle (1.2pt);
\filldraw [black] (12,0) circle (1.2pt);
\node at (4.5,-0.4){$a_1'$};
\node at (6,-0.4){$a_2'$};
\node at (7.5,-0.4){$a_3'$};
\node at (9,-0.4){$e_1'$};
\node at (10.5,-0.4){$e_2'$};
\node at (12,-0.4){$e_3'$};
\draw[rounded corners=0.3cm,line width =1pt] (4.2,-0.15)--(7.8,-0.15)--(6,2.3)--cycle;
\draw[rounded corners=0.3cm,line width =1pt] (5.7,-0.15)--(9.3,-0.15)--(7.5,2.3)--cycle;
\draw[rounded corners=0.3cm,line width =1pt] (7.2,-0.15)--(10.8,-0.15)--(9,2.3)--cycle;
\draw[rounded corners=0.3cm,line width =1pt] (8.7,-0.15)--(12.3,-0.15)--(10.5,2.3)--cycle;
\filldraw [black] (5.77,-1.96) circle (1pt);
\filldraw [black] (8.03,-1.96) circle (1pt);
\filldraw [black] (10.2,-1.96) circle (1pt);
\node at (5.8,-1.65){$b_1'$};
\node at (8,-1.65){$b_2'$};
\node at (10.2,-1.65){$b_3'$};
\filldraw [black] (5,-3) circle (1.2pt);
\filldraw [black] (5.5,-3.1) circle (1.2pt);
\filldraw [black] (6,-3.1) circle (1.2pt);
\filldraw [black] (6.5,-3) circle (1.2pt);
\filldraw [yellow] (5.25,-3.05) circle (1.2pt);
\filldraw [yellow] (5.75,-3.1) circle (1.2pt);
\filldraw [yellow] (6.25,-3.05) circle (1.2pt);
\node[yellow] at (5.75,-3.6){$C_1'$};
\draw[rounded corners=0.2cm,line width =1pt] (4.8,-3.1)--(5.65,-3.25)--(5.85,-1.75)--cycle;
\draw[rounded corners=0.2cm,line width =1pt] (5.3,-3.2)--(6.2,-3.2)--(5.78,-1.75)--cycle;
\draw[rounded corners=0.2cm,line width =1pt] (5.85,-3.25)--(6.75,-3.05)--(5.73,-1.75)--cycle;
\filldraw [black] (8.3,-3.1) circle (1.2pt);
\filldraw [black] (8.8,-3) circle (1.2pt);
\filldraw [black] (7.8,-3.1) circle (1.2pt);
\filldraw [black] (7.3,-3) circle (1.2pt);
\filldraw [purple] (7.55,-3.05) circle (1.2pt);
\filldraw [purple] (8.05,-3.1) circle (1.2pt);
\filldraw [purple] (8.55,-3.05) circle (1.2pt);
\node[purple] at (8.05,-3.6){$C_2'$};
\draw[rounded corners=0.2cm,line width =1pt] (7.1,-3.1)--(7.95,-3.25)--(8.1,-1.75)--cycle;
\draw[rounded corners=0.2cm,line width =1pt] (7.6,-3.2)--(8.5,-3.2)--(8.05,-1.75)--cycle;
\draw[rounded corners=0.2cm,line width =1pt] (8.15,-3.25)--(9.05,-3.05)--(7.97,-1.765)--cycle;
\filldraw [black] (11,-3) circle (1.2pt);
\filldraw [black] (10.5,-3.1) circle (1.2pt);
\filldraw [black] (10,-3.1) circle (1.2pt);
\filldraw [black] (9.5,-3) circle (1.2pt);
\filldraw [pink] (9.75,-3.05) circle (1.2pt);
\filldraw [pink] (10.25,-3.1) circle (1.2pt);
\filldraw [pink] (10.75,-3.05) circle (1.2pt);
\node[pink] at (10.25,-3.6){$C_3'$};
\draw[rounded corners=0.2cm,line width =1pt] (9.3,-3.1)--(10.15,-3.25)--(10.25,-1.75)--cycle;
\draw[rounded corners=0.2cm,line width =1pt] (9.82,-3.2)--(10.7,-3.2)--(10.2,-1.75)--cycle;
\draw[rounded corners=0.2cm,line width =1pt] (10.35,-3.25)--(11.25,-3.05)--(10.15,-1.765)--cycle;
\end{tikzpicture}
\end{center}
\caption{Before the Absorption ($k=3$), the paths $(C,AE),(C',A'E')$,
$(C_1,P_1b_1Q_1),(C_1',P_1'b_1'Q_1'),(C_2,P_2b_2Q_2),(C_2',P_2'b_2'Q_2'),(C_3,P_3b_3Q_3),(C_3',P_3'b_3'Q_3')$ are sequentially paths, the black dots represent the points while other dots with the same color come from the same color set.}
\end{figure}
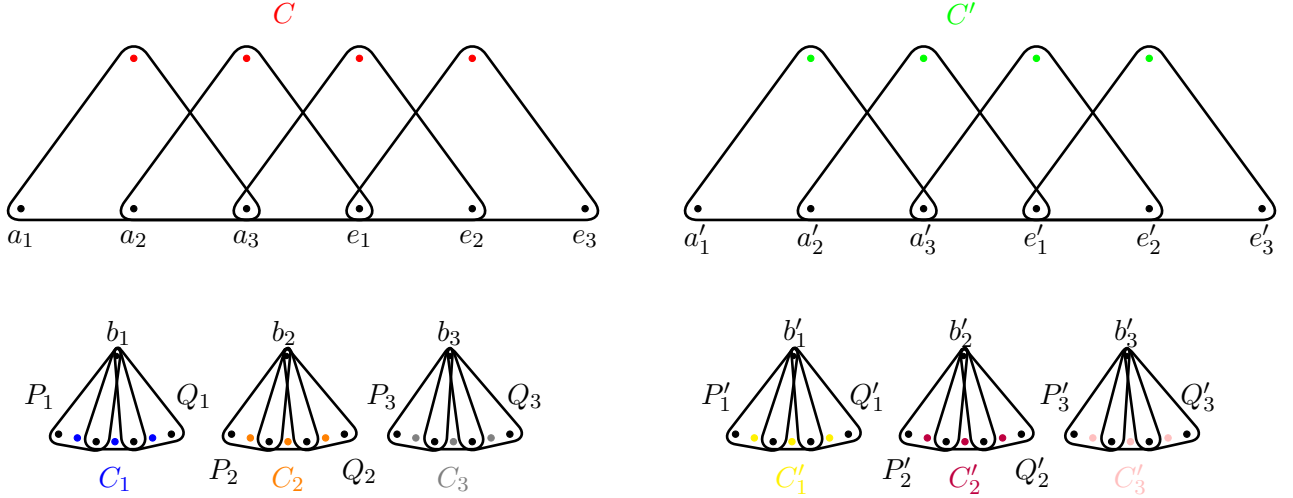

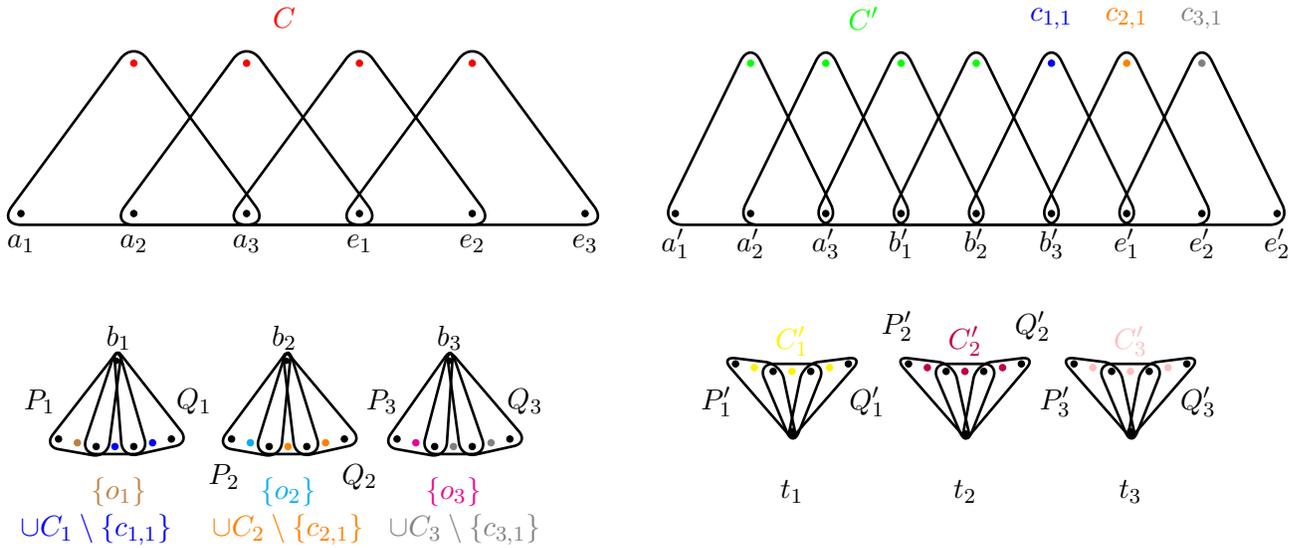
\begin{figure}[htb]
\begin{center}
\begin{tikzpicture}{center}
\filldraw [red] (-3,2) circle (1.2pt);
\filldraw [red] (-1.5,2) circle (1.2pt);
\filldraw [red] (0,2) circle (1.2pt);
\filldraw [red] (1.5,2) circle (1.2pt);
\node [red] at (-1,2.6){$C$};
\filldraw [black] (-4.5,0) circle (1.2pt);
\filldraw [black] (-3,0) circle (1.2pt);
\filldraw [black] (-1.5,0) circle (1.2pt);
\filldraw [black] (0,0) circle (1.2pt);
\filldraw [black] (1.5,0) circle (1.2pt);
\filldraw [black] (3,0) circle (1.2pt);
\node at (-4.5,-0.4){$a_1$};
\node at (-3,-0.4){$a_2$};
\node at (-1.5,-0.4){$a_3$};
\node at (0,-0.4){$e_1$};
\node at (1.5,-0.4){$e_2$};
\node at (3,-0.4){$e_3$};
\draw[rounded corners=0.3cm,line width =1pt] (-4.8,-0.15)--(-1.2,-0.15)--(-3,2.3)--cycle;
\draw[rounded corners=0.3cm,line width =1pt] (-3.3,-0.15)--(0.3,-0.15)--(-1.5,2.3)--cycle;
\draw[rounded corners=0.3cm,line width =1pt] (-1.8,-0.15)--(1.8,-0.15)--(0,2.3)--cycle;
\draw[rounded corners=0.3cm,line width =1pt] (-0.3,-0.15)--(3.3,-0.15)--(1.5,2.3)--cycle;
\filldraw [black] (-3.23,-1.96) circle (1pt);
\filldraw [black] (-0.97,-1.96) circle (1pt);
\filldraw [black] (1.2,-1.96) circle (1pt);
\node at (-3.2,-1.65){$b_1$};
\node at (-1,-1.65){$b_2$};
\node at (1.2,-1.65){$b_3$};
\filldraw [black] (-4,-3) circle (1.2pt);
\filldraw [black] (-3.5,-3.1) circle (1.2pt);
\filldraw [black] (-3,-3.1) circle (1.2pt);
\filldraw [black] (-2.5,-3) circle (1.2pt);
\filldraw [brown] (-3.75,-3.05) circle (1.2pt);
\filldraw [blue] (-3.25,-3.1) circle (1.2pt);
\filldraw [blue] (-2.75,-3.05) circle (1.2pt);
\node[brown] at (-3.2,-3.7){$\{o_1\}$};
\node[blue] at (-3.5,-4.2){$\cup C_1\setminus\{c_{1,1}\}$};

\node[black] at (-4.25,-2.5){$P_1$};
\node[black] at (-2.2,-2.5){$Q_1$};
\node[black] at (-1.8,-3.5){$P_2$};
\node[black] at (0,-3.5){$Q_2$};
\node[black] at (0.3,-2.5){$P_3$};
\node[black] at (2.2,-2.5){$Q_3$};
\node[black] at (4.75,-2.5){$P_1'$};
\node[black] at (6.75,-2.5){$Q_1'$};
\node[black] at (7.15,-1.5){$P_2'$};
\node[black] at (8.95,-1.5){$Q_2'$};
\node[black] at (9.25,-2.5){$P_3'$};
\node[black] at (11.15,-2.5){$Q_3'$};
\draw[rounded corners=0.2cm,line width =1pt] (-4.2,-3.1)--(-3.35,-3.25)--(-3.15,-1.75)--cycle;
\draw[rounded corners=0.2cm,line width =1pt] (-3.7,-3.2)--(-2.8,-3.2)--(-3.22,-1.75)--cycle;
\draw[rounded corners=0.2cm,line width =1pt] (-3.15,-3.25)--(-2.25,-3.05)--(-3.27,-1.75)--cycle;
\filldraw [black] (-0.7,-3.1) circle (1.2pt);
\filldraw [black] (-0.2,-3) circle (1.2pt);
\filldraw [black] (-1.2,-3.1) circle (1.2pt);
\filldraw [black] (-1.7,-3) circle (1.2pt);
\filldraw [cyan] (-1.45,-3.05) circle (1.2pt);
\filldraw [orange] (-0.95,-3.1) circle (1.2pt);
\filldraw [orange] (-0.45,-3.05) circle (1.2pt);
\node[cyan]at (-0.95,-3.7){$\{o_2\}$};
\node[orange] at (-0.95,-4.2){$\cup C_2\setminus\{c_{2,1}\}$};
\draw[rounded corners=0.2cm,line width =1pt] (-1.9,-3.1)--(-1.05,-3.25)--(-0.9,-1.75)--cycle;
\draw[rounded corners=0.2cm,line width =1pt] (-1.4,-3.2)--(-0.5,-3.2)--(-0.95,-1.75)--cycle;
\draw[rounded corners=0.2cm,line width =1pt] (-0.85,-3.25)--(0.05,-3.05)--(-1.03,-1.765)--cycle;
\filldraw [black] (2,-3) circle (1.2pt);
\filldraw [black] (1.5,-3.1) circle (1.2pt);
\filldraw [black] (1,-3.1) circle (1.2pt);
\filldraw [black] (0.5,-3) circle (1.2pt);
\filldraw [magenta] (0.75,-3.05) circle (1.2pt);
\filldraw [gray] (1.25,-3.1) circle (1.2pt);
\filldraw [gray] (1.75,-3.05) circle (1.2pt);
\node [magenta] at (1.25,-3.7){$\{o_3\}$};
\node [gray] at (1.4,-4.2){$\cup C_3\setminus\{c_{3,1}\}$};
\draw[rounded corners=0.2cm,line width =1pt] (0.3,-3.1)--(1.15,-3.25)--(1.25,-1.75)--cycle;
\draw[rounded corners=0.2cm,line width =1pt] (0.82,-3.2)--(1.7,-3.2)--(1.2,-1.75)--cycle;
\draw[rounded corners=0.2cm,line width =1pt] (1.35,-3.25)--(2.25,-3.05)--(1.15,-1.765)--cycle;

\filldraw [green] (5.2,2) circle (1.2pt);
\filldraw [green] (6.2,2) circle (1.2pt);
\filldraw [green] (7.2,2) circle (1.2pt);
\filldraw [green] (8.2,2) circle (1.2pt);
\filldraw [blue] (9.2,2) circle (1.2pt);
\filldraw [orange] (10.2,2) circle (1.2pt);
\filldraw [gray] (11.2,2) circle (1.2pt);
\node [green] at (6.7,2.6){$C'$};
\node [blue] at (9.2,2.6){$c_{1,1}$};
\node [orange] at (10.2,2.6){$c_{2,1}$};
\node [gray] at (11.2,2.6){$c_{3,1}$};
\filldraw [black] (4.2,0) circle (1.2pt);
\filldraw [black] (5.2,0) circle (1.2pt);
\filldraw [black] (6.2,0) circle (1.2pt);
\filldraw [black] (7.2,0) circle (1.2pt);
\filldraw [black] (8.2,0) circle (1.2pt);
\filldraw [black] (9.2,0) circle (1.2pt);
\filldraw [black] (10.2,0) circle (1.2pt);
\filldraw [black] (11.2,0) circle (1.2pt);
\filldraw [black] (12.2,0) circle (1.2pt);
\node at (4.2,-0.4){$a_1'$};
\node at (5.2,-0.4){$a_2'$};
\node at (6.2,-0.4){$a_3'$};
\node at (7.2,-0.4){$b_1'$};
\node at (8.2,-0.4){$b_2'$};
\node at (9.2,-0.4){$b_3'$};
\node at (10.2,-0.4){$e_1'$};
\node at (11.2,-0.4){$e_2'$};
\node at (12.2,-0.4){$e_2'$};
\draw[rounded corners=0.3cm,line width =1pt] (4,-0.15)--(6.4,-0.15)--(5.2,2.3)--cycle;
\draw[rounded corners=0.3cm,line width =1pt] (5,-0.15)--(7.4,-0.15)--(6.2,2.3)--cycle;
\draw[rounded corners=0.3cm,line width =1pt] (6,-0.15)--(8.4,-0.15)--(7.2,2.3)--cycle;
\draw[rounded corners=0.3cm,line width =1pt] (7,-0.15)--(9.4,-0.15)--(8.2,2.3)--cycle;
\draw[rounded corners=0.3cm,line width =1pt] (8,-0.15)--(10.4,-0.15)--(9.2,2.3)--cycle;
\draw[rounded corners=0.3cm,line width =1pt] (9,-0.15)--(11.4,-0.15)--(10.2,2.3)--cycle;
\draw[rounded corners=0.3cm,line width =1pt] (10,-0.15)--(12.4,-0.15)--(11.2,2.3)--cycle;

\filldraw [black] (5,-2) circle (1.2pt);
\filldraw [black] (5.5,-2.1) circle (1.2pt);
\filldraw [black] (6,-2.1) circle (1.2pt);
\filldraw [black] (6.5,-2) circle (1.2pt);
\filldraw [yellow] (5.25,-2.05) circle (1.2pt);
\filldraw [yellow] (5.75,-2.1) circle (1.2pt);
\filldraw [yellow] (6.25,-2.05) circle (1.2pt);
\filldraw [black] (5.75,-2.9) circle (1.2pt);
\node[yellow] at (5.75,-1.7){$C_1'$};
\draw[rounded corners=0.2cm,line width =1pt] (4.8,-1.9)--(5.6,-2)--(5.83,-3.085)--cycle;
\draw[rounded corners=0.2cm,line width =1pt] (5.3,-2)--(6.2,-2)--(5.765,-3.1)--cycle;
\draw[rounded corners=0.2cm,line width =1pt] (5.9,-2)--(6.7,-1.9)--(5.68,-3.07)--cycle;

\filldraw [black] (8.3,-2.1) circle (1.2pt);
\filldraw [black] (8.8,-2) circle (1.2pt);
\filldraw [black] (7.8,-2.1) circle (1.2pt);
\filldraw [black] (7.3,-2) circle (1.2pt);
\filldraw [purple] (7.55,-2.05) circle (1.2pt);
\filldraw [purple] (8.05,-2.1) circle (1.2pt);
\filldraw [purple] (8.55,-2.05) circle (1.2pt);
\filldraw [black] (8.05,-2.93) circle (1.2pt);
\node[purple] at (8.05,-1.7){$C_2'$};
\draw[rounded corners=0.2cm,line width =1pt] (7.1,-1.9)--(7.9,-2)--(8.14,-3.085)--cycle;
\draw[rounded corners=0.2cm,line width =1pt] (7.6,-2)--(8.5,-2)--(8.065,-3.1)--cycle;
\draw[rounded corners=0.2cm,line width =1pt] (8.2,-2)--(9,-1.9)--(7.995,-3.09)--cycle;

\filldraw [black] (11,-2) circle (1.2pt);
\filldraw [black] (10.5,-2.1) circle (1.2pt);
\filldraw [black] (10,-2.1) circle (1.2pt);
\filldraw [black] (9.5,-2) circle (1.2pt);
\filldraw [pink] (9.75,-2.05) circle (1.2pt);
\filldraw [pink] (10.25,-2.1) circle (1.2pt);
\filldraw [pink] (10.75,-2.05) circle (1.2pt);
\filldraw [black] (10.25,-2.9) circle (1.2pt);
\node[pink] at (10.25,-1.7){$C_3'$};
\draw[rounded corners=0.2cm,line width =1pt] (9.3,-1.9)--(10.1,-2)--(10.34,-3.085)--cycle;
\draw[rounded corners=0.2cm,line width =1pt] (9.8,-2)--(10.7,-2)--(10.265,-3.1)--cycle;
\draw[rounded corners=0.2cm,line width =1pt] (10.4,-2)--(11.2,-1.9)--(10.2,-3.05)--cycle;
\node[black] at (5.75,-3.7){$t_1$};
\node[black] at (8.05,-3.7){$t_2$};
\node[black] at (10.25,-3.7){$t_3$};
\end{tikzpicture}
\end{center}
\caption{After the Absorption ($k=3$), the paths
$(\{o_1\}\cup C_1\setminus\{c_{1,1}\},P_1b_1Q_1)$,
$(\{o_2\}\cup C_2\setminus\{c_{2,1}\},P_2b_2Q_2)$,
$(\{o_3\}\cup C_3\setminus\{c_{3,1}\},P_3b_3Q_3)$,
$(C_1',P_1't_1Q_1')$, $(C_2',P_2't_2Q_2')$,$(C_3',P_3't_3Q_3')$,
$(C,AE)$,
$(C'(c_{1,1},c_{2,1},c_{3,1}),A'B'E')$ are sequentially paths, the black dots represent the points  while other dots with the same color come from the same color set.}
\end{figure}
\begin{definition}[$\mathfrak{S}$-gadget]
\label{10.2}
\label{gadget}
Suppose $F=F_1\cup F_2$ is an absorbing gadget where $F_1=A\cup B\cup E\cup \bigcup_{i=1}^k(P_i\cup Q_i)\cup C\cup \bigcup_{i=1}^kC_k$ and $F_2=A'\cup B'\cup E'\cup \bigcup_{i=1}^k(P_i'\cup Q_i')\cup C'\cup \bigcup_{i=1}^kC_k'$   with $A=(a_1,\ldots,a_k)$, $B=(b_1,\ldots,b_k)$, $E=(e_1,\ldots,e_k)$, $C=(c_1,\ldots,c_{k+1})$, $C_i=(c_{i,1},\ldots,c_{i,k})$, $P_i=(p_{i,1},\ldots,p_{i,k-1})$ and $Q_i=(q_{i,1},\ldots,q_{i,k-1})$ for $i\in[k]$,
$A'=(a_1',\ldots,a_k')$, $B'=(b_1',\ldots,b_k')$, $E'=(e_1',\ldots,e_k')$, $C'=(c_1',\ldots,c_{k+1}')$, $C_i'=(c_{i,1}',\ldots,c_{i,k}')$, $P_i'=(p_{i,1}',\ldots,p_{i,k-1}')$ and $Q_i'=(q_{i,1}',\ldots,q_{i,k-1}')$ for $i\in[k]$.
Suppose that $\varepsilon,\varepsilon_{k+1},d_2,\ldots,d_{k+1},c,\nu>0$.
Let $\textbf{d}=(d_2,\ldots,d_{k+1})$ and suppose that $\mathfrak{S}=(G,G_{\mathcal{J}},\mathcal{J},\mathcal{P},\overrightarrow{H})$ is an oriented $(k+1,m,2t,\varepsilon,\varepsilon_{k+1},r,\textbf{d})$-regular setup.
We say that $F$ is an $\mathfrak{S}$-gadget if
\begin{itemize}
  \item[(G1)] there exists an oriented edge $Y'=(Y_0,Z_1,\ldots,Z_k)\in \overrightarrow{H}$ and a color cluster $Z_0$, such that $C\cup C'\cup\bigcup_{i\in[k]}C_i\subseteq Y_0$, $\bigcup_{i\in[k]}C_i'\subseteq Z_0$, $a_i, b_i, e_i\in Z_i$ for $i\in[k]$,
  \item[(G2)] there exists an oriented edge $Y=(Y_0,Y_1,\ldots,Y_k)\in \overrightarrow{H}$, such that $a_i', b_i', e_i'\in Y_i$ for $i\in[k]$,
  \item[(G3)] there exists an ordered $k$-tuple of clusters $W_i=(W_{i,1},\ldots, W_{i,k-1})$ such that $W_i\cup \{Y_0,Z_i\}$ is an edge in $H$ and $(Y_0, W_{i,1},\ldots, W_{i,k-1},Z_i)$ is consistent with $\overrightarrow{H}$, $p_{i,j},q_{i,j}\in W_{i,j}$ for $i\in[k], j\in[k-1]$,
  \item[(G4)] there exists an ordered $k$-tuple of clusters $W_i'=(W_{i,1}',\ldots, W_{i,k-1}')$ such that $W_i'\cup \{Z_0,Y_i\}$ is an edge in $H$ and $(Z_0, W_{i,1}',\ldots, W_{i,k-1}',Y_i)$ is consistent with $\overrightarrow{H}$, $p_{i,j}',q_{i,j}'\in W_{i,j}'$ for $i\in[k], j\in[k-1]$,
  \item[(G5)] $F\subseteq G_{\mathcal{J}}$,
\end{itemize}

We will further say that $F$ is $(c,\nu)$-extensible if the following also holds:
\begin{itemize}
  \item[(G6)] The path $(C,AE)$ is $(c,\nu)$-extensible both left- and rightwards to the ordered tuple $Y'=(Y_0,Z_1,\ldots,Z_k)$ and the path $(C_i,P_ib_iQ_i)$ is $(c,\nu)$-extensible leftwards to $(Y_0,W_{i,1},\ldots, W_{i,k-1},Z_i)$ and rightwards to $(Y_0,Z_i,W_{i,1},\ldots, W_{i,k-1})$ for $i\in[k]$.
  \item[(G7)]  The path $(C',A'E')$ is $(c,\nu)$-extensible both left- and rightwards to the ordered tuple $Y=(Y_0,Y_1,\ldots,Y_k)$ and the path $(C_i',P_i'b_i'Q_i')$ is $(c,\nu)$-extensible leftwards to $(Z_0,W_{i,1}',\ldots, W_{i,k-1}',Y_i)$ and rightwards to $(Z_0,Y_i,W_{i,1}',\ldots, W_{i,k-1}')$ for $i\in[k]$.
\end{itemize}

\end{definition}
\begin{definition}[Reduced gadget]
\label{reduced gadget}
A reduced gadget is a $(1,k)$-graph $L$ consisting of $Y\cup W_1\cup \cdots\cup W_k\cup Z_0\cup Z_1\cup\ldots\cup Z_k\cup W_1'\cup \cdots\cup W_k'$ where $Y=\{Y_0,Y_1,\ldots,Y_k\}$, $W_i=\{W_{i,1},\ldots, W_{i,k-1}\}$ for $i\in[k]$, $W_i'=\{W_{i,1}',\ldots, W_{i,k-1}'\}$ for $i\in[k]$ and $2(k+1)$ edges given by $Y, Y'=\{Y_0, Z_1,\ldots,Z_k\}$, $W_i\cup \{Y_0,Z_i\}$ for $i\in[k]$ and $W_i'\cup \{Z_0,Y_i\}$ for $i\in[k]$.
We refer to $Y$ and $Y'$ as the core edges of $L$ and $W_i, W_i', i\in[k]$ as the peripheral sets of $L$.
\end{definition}
\vspace{12pt}
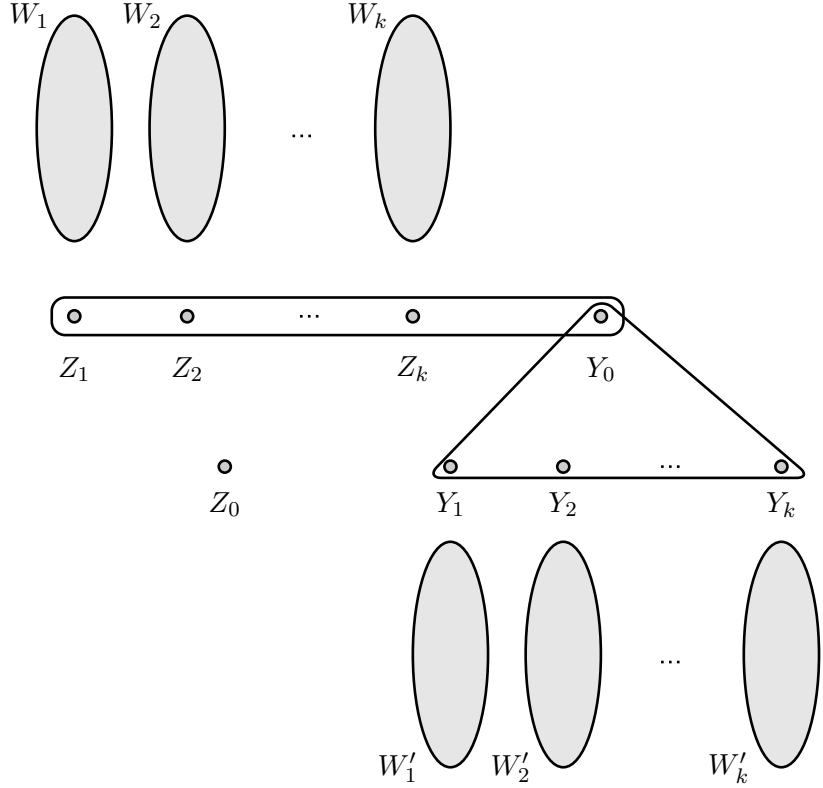
\begin{figure}[htb]
\begin{center}
\begin{tikzpicture}{center}
\draw [black,fill=gray!40,line width =1pt] (-3,2) circle(0.08);
\draw [black,fill=gray!40,line width =1pt] (0,2) circle (0.08);
\draw [black,fill=gray!40,line width =1pt] (1.5,2) circle (0.08);
\draw [black,fill=gray!40,line width =1pt] (4.4,2) circle (0.08);
\draw [black,fill=gray!40,line width =1pt] (2,4) circle (0.08);
\node [black] at (-3,1.5){$Z_0$};
\node [black] at (0,1.5){$Y_1$};
\node [black] at (1.5,1.5){$Y_2$};
\draw[dotted,line width =1pt] (2.8,2)--(3.05,2);
\draw[dotted,line width =1pt] (2.8,-0.6)--(3.05,-0.6);
\node [black] at (4.4,1.5){$Y_k$};
\node [black] at (2,3.3){$Y_0$};
\node [black] at (-5,3.3){$Z_1$};
\node [black] at (-3.5,3.3){$Z_2$};
\node [black] at (-0.5,3.3){$Z_k$};
\node [black] at (-0.7,-2){$W_1'$};
\node [black] at (0.8,-2){$W_2'$};
\node [black] at (3.7,-2){$W_k'$};
\draw[rounded corners=0.2cm,line width =1pt] (-0.32,1.85)--(4.8,1.85)--(2,4.25)--cycle;
\draw[rounded corners=0.18cm,line width =1pt] (2.3,3.75)rectangle(-5.3,4.25);
\filldraw[draw=black,fill=gray!20,line width =1pt](0,-0.5) ellipse (0.5 and 1.5);
\filldraw[draw=black,fill=gray!20,line width =1pt](1.5,-0.5) ellipse (0.5 and 1.5);
\filldraw[draw=black,fill=gray!20,line width =1pt](4.4,-0.5) ellipse (0.5 and 1.5);

\draw [black,fill=gray!40,line width =1pt] (-5,4) circle (0.08);
\draw [black,fill=gray!40,line width =1pt] (-3.5,4) circle (0.08);
\draw[dotted,line width =1pt] (-2,4)--(-1.75,4);
\draw [black,fill=gray!40,line width =1pt] (-0.5,4) circle (0.08);
\filldraw[draw=black,fill=gray!20,line width =1pt](-5,6.5) ellipse (0.5 and 1.5);
\filldraw[draw=black,fill=gray!20,line width =1pt](-3.5,6.5) ellipse (0.5 and 1.5);
\filldraw[draw=black,fill=gray!20,line width =1pt](-0.5,6.5) ellipse (0.5 and 1.5);
\node [black] at (-5.6,8){$W_1$};
\node [black] at (-4.1,8){$W_2$};
\node [black] at (-1.1,8){$W_k$};
\draw[dotted,line width =1pt] (-2.1,6.4)--(-1.8,6.4);
\end{tikzpicture}
\end{center}
\caption{Reduced Gadget}
\end{figure}
Given an oriented $(1,k)$-graph $\overrightarrow{H}$, a reduced gadget in $\overrightarrow{H}$ is a copy of $L$ such that $Y$ coincides with the orientation of that edge in $\overrightarrow{H}$ and such that $(Z_0, W_{i,1},\ldots, W_{i,k-1},Y_i)$ is consistent with that edge in $\overrightarrow{H}$.

Let $\mathfrak{S}=(G,G_{\mathcal{J}},\mathcal{J},\mathcal{P},\overrightarrow{H})$ be an oriented regular setup.
Let $c,\nu>0$, $T=\{t_1,\ldots,t_k\}$ be a $k$-set of $V$ and $O=\{o_1,\ldots,o_k\}$ be a $k$-set of $[n]$, and $L$ be a reduced gadget in $\overrightarrow{H}$. We define the following sets:
\begin{enumerate}
  \item Denote the set of all reduced gadgets in $\overrightarrow{H}$ by $\mathfrak{L}_{\overrightarrow{H}}$,
  \item Denote the set of $\mathfrak{S}$-gadgets which use precisely the clusters of $L$ as in Definition \ref{reduced gadget} by $\mathfrak{F}_{L}$,
  \item Denote the set of $\mathfrak{S}$-gadgets in $\mathfrak{F}_{L}$ which are $(c,\nu,V(G))$-extensible by $\mathfrak{F}_{L}^{\rm ext}$,
  \item Denote the set of all $\mathfrak{S}$-gadgets by $\mathfrak{F}$,
  \item Denote the set of all $(c,\nu,V(G))$-extensible $\mathfrak{S}$-gadgets by $\mathfrak{F}^{\rm ext}\subseteq\mathfrak{F}$,
  \item For any $k$-subset $T$ of $V$ and any $k$-subset $O$ of $[n]$, let $\mathfrak{F}_{(T,O)}\subseteq\mathfrak{F}$ be the set of absorbing $\mathfrak{S}$-gadgets for $(T,O)$,
  \item Denote the set of $\mathfrak{S}$-gadgets absorbing $(T,O)$ which are $(c,\nu)$-extensible by $\mathfrak{F}_{(T,O)}^{\rm ext}=\mathfrak{F}_{(T,O)}\cap\mathfrak{F}^{\rm ext}$.
\end{enumerate}
\begin{lemma}
\label{10.5}
Let $k,r,m,t\in \mathbb{N}$ and $d_2,\ldots,d_{k+1},\varepsilon,\varepsilon_{k+1},c,\nu,\beta$ be such that
\begin{equation}
\begin{split}
1/m&\ll1/r,\varepsilon\ll1/t,c,\varepsilon_{k+1},d_2,\ldots,d_k,\\
c&\ll d_2,\ldots,d_k,\\
1/t&\ll\varepsilon_{k+1}\ll \beta,d_{k+1}\leq1/k,\\
\varepsilon_{k+1}&\ll\nu.
\nonumber
\end{split}
\end{equation}
Let $\textbf{d}=(d_2,\ldots,d_{k+1})$ and let $\mathfrak{S}=(G,G_{\mathcal{J}},\mathcal{J},\mathcal{P},\overrightarrow{H})$ be an oriented $(k,m,2t,\varepsilon,\varepsilon_{k+1},r,\textbf{d})$-regular setup and $L\in\mathcal{L}_{\overrightarrow{H}}$ be a reduced gadget in $\overrightarrow{H}$.
Let $\mathcal{F}$ be the $(k+1)$-complex corresponding to the down-closure of $(1,k)$-graph $F$ as in Definition \ref{gadget}.
Then
\begin{equation}\label{7.6}
|\mathfrak{F}_L|=(1\pm\beta)\left(\prod_{i=2}^{k+1}d_i^{e_i(\mathcal{F})}\right)m^{6k^2+4k+2},
\end{equation}
\begin{equation}
|\mathfrak{F}_L\setminus\mathfrak{F}_L^{\rm ext}|\leq\beta|\mathfrak{F}_L|.
\nonumber
\end{equation}
\end{lemma}
\begin{proof}
Let $Y=(Y_0,Y_1,\ldots,Y_k), Y'=(Y_0,Z_1,\ldots,Z_k)\in \overrightarrow{H}$ be the ordered core edge of $L$ and $W_i=\{W_{i,1},\ldots,W_{i,k-1}\}$, $W_i'=\{W_{i,1}',\ldots,W_{i,k-1}'\}$ for $i\in[k]$, be the peripheral sets, ordered such that $(Y_0,W_{i,1},\ldots,W_{i,k-1},Z_i)$ and $(Z_0,W_{i,1}',\ldots,W_{i,k-1}',Y_i)$ are consistent with $\overrightarrow{H}$.
Note that $|V(F)|=6k^2+4k+2$.
The bounds on $|\mathfrak{F}_L|$ are given by Lemma \ref{counting lemma} directly.

Let $Y^*=(Y_1,\ldots,Y_{k-1})$ and denote the ordered tuples in the $(k-1)$-th level of $\mathcal{J}$ in the clusters $\{Y_1,\ldots,Y_{k-1}\}$ by $\mathcal{J}_{Y^*}$.
Let $d_{Y^*}=\prod_{i=2}^{k-1}d_i^{\binom{k-1}{i}}$.
By Lemma \ref{7.7} we have
\begin{equation}
|\mathcal{J}_{Y^*}|=(1\pm\beta)d_{Y^*}m^{k-1}.
\nonumber
\end{equation}

Let $\beta_1$ be such that $\varepsilon_{k+1}\ll\beta_1\ll\beta,d_k,d_{k+1},1/k$.
Let $B_1\subseteq\mathcal{J}_{Y^*}$ be the set of $(k-1)$-tuples which are not $(c,\nu)$-extensible leftwards to $(Y_0,Y_1,\ldots,Y_k)$.
By Proposition \ref{8.1} with $\beta_1$ playing the role of $\beta$, we deduce that
\begin{equation}
|B_1|\leq\beta_1|\mathcal{J}_{Y^*}|.
\nonumber
\end{equation}

Let $\beta_2$ be such that $\varepsilon\ll\beta_2\ll\varepsilon_{k+1},d_2,\ldots,d_{k-1}$.
Let $\phi: V(F)\rightarrow L$ be the homomorphism and $Z\subseteq V(F)$ corresponds to the first $k-1$ vertices $\{a_1,\ldots,a_{k-1}\}$ of path $AE$.
Let $\mathcal{F}^-$ be the $(k-1)$-complex generated by removing the $(k+1)$-st and $k$-th layer from the down-closure $\mathcal{F}$ of $F$.
Let $\mathcal{Z}=\mathcal{F}^-[Z]$ be the induced subcomplex of $\mathcal{F}^-$ in $Z$.
Note that $\phi(a_i)=Y_i$ for $i\in[k-1]$.
Thus the labelled partition-respecting copies of $\mathcal{Z}$ in $\mathcal{J}$ correspond exactly to $\mathcal{J}_{Y^*}$.
Define
\[
d_{\mathcal{F}^-\setminus \mathcal{Z}}=\prod_{i=2}^{k-1}d_i^{e_i(\mathcal{F}^-)-e_i(\mathcal{Z})}.
\]
Let $B_2\subseteq\mathcal{J}_{Y^*}$ be the set of $(k-1)$-tuples which are not extensible to $(1\pm\beta_2)d_{\mathcal{F}^-\setminus \mathcal{Z}}m^{6k^2+3k+3}$
labelled partition-respecting copies of $\mathcal{F}^-$ in $\mathcal{J}$.
By Lemma \ref{dense extension} with $\beta_2$ playing the role of $\beta$, we have
\begin{equation}
|B_2|\leq\beta_2|\mathcal{J}_{Y^*}|.
\nonumber
\end{equation}
By (\ref{7.6}), we have
\begin{equation}
|\mathfrak{F}_L|=(1\pm\beta)d_{k+1}^{e_{k+1}(\mathcal{F})}d_k^{e_k(\mathcal{F})}d_{\mathcal{F}^-\setminus \mathcal{Z}}d_{Y^*}m^{6k^2+4k+2}.
\nonumber
\end{equation}

Let $\mathcal{G}=\mathcal{J}\cup G_{\mathcal{J}}$.
Say that a labelled partition-respecting copy of $\mathcal{F}$ in $\mathcal{G}$ is \emph{nice} if the vertices of $\{a_1,\ldots,a_{k-1}\}$ are not in $B_1\cup B_2$.
For every $Z\in\mathcal{J}_{Y^*}$, let $N^*(Z)$ be the number of labelled partition-respecting copies of $\mathcal{F}$ in $\mathcal{G}$ which extend $Z$.
Note that $0\leq N^*(Z)\leq m^{6k^2+3k+3}$ and we have
\begin{equation}
\begin{split}
\sum_{Z\in B_1\cup B_2}N^*(Z)&=\sum_{Z\in B_1\setminus B_2}N^*(Z)+\sum_{Z\in B_2}N^*(Z)\\
&\leq[|B_1|(1+\beta_2)d_{\mathcal{F}^-\setminus \mathcal{Z}}+|B_2|]m^{6k^2+3k+3}\\
&\leq[\beta_1(1+\beta_2)d_{\mathcal{F}^-\setminus \mathcal{Z}}+\beta_2]|\mathcal{J}_{Y'}|m^{6k^2+3k+3}\\
&\leq3\beta_1d_{\mathcal{F}^-\setminus \mathcal{Z}}|\mathcal{J}_{Y^*}|m^{6k^2+3k+3}\\
&\leq3\beta_1(1+\beta)d_{\mathcal{F}^-\setminus \mathcal{Z}}d_{Y^*}m^{6k^2+4k+2}\\
&\leq\frac{3\beta_1(1+\beta)}{(1-\beta)d_{k+1}^{e_{k+1}(\mathcal{F})}d_k^{e_k(\mathcal{F})}}|\mathcal{F}_L|\\
&\leq\frac{\beta}{4k+4}|\mathcal{F}_L|,
\nonumber
\end{split}
\end{equation}
since $\beta_1\ll\beta,d_k,d_{k+1},1/k$ and $\beta_2\ll d_2,\ldots,d_{k-1},\varepsilon_{k+1}$.

The same analysis shows that we define nice tuples for any $(k-1)$-set of vertices of $F$, the number of copies of $F$ which are not nice with respect to that $(k-1)$-set is at most $\beta|\mathcal{F}_L|/(4k+4)$.
Note that $F\in\mathfrak{F}_L$ is extensible if and only if paths $(C,AE)$, $(C',A'E')$, $(C_i,P_ib_iQ_i)$ and $(C_i',P_i'b_i'Q_i')$ for $i\in[k]$ contained in $F$ are extensible with certain edges of the reduced graph.
This means that $4(k+1)$ many $(k-1)$-tuples are extensible with certain edges of the reduced graph.
Thus, $F\in \mathfrak{F}_L\setminus\mathfrak{F}_L^{\rm ext}$ implies that $F$ is not nice with one of $4k+4$ many $(k-1)$-sets.
Thus,
\[
|\mathfrak{F}_L\setminus\mathfrak{F}_L^{\rm ext}|\leq (4k+4)\frac{\beta}{4k+4}|\mathcal{F}_L|=\beta|\mathcal{F}_L|.
\]
\end{proof}
\begin{lemma}
\label{10.6}
Let $k,r,m,t\in \mathbb{N}$ and $d_2,\ldots,d_{k+1},\varepsilon,\varepsilon_{k+1},c,\nu,\beta,\mu$ be such that
\begin{equation}
\begin{split}
1/m&\ll1/r,\varepsilon\ll1/t,c,\varepsilon_{k+1},d_2,\ldots,d_k,\\
c&\ll d_2,\ldots,d_k,\\
1/t&\ll\varepsilon_{k+1}\ll \beta,d_{k+1}\leq1/k,\\
\varepsilon_{k+1}&\ll\nu,\mu,\\
\alpha&\ll\mu.
\nonumber
\end{split}
\end{equation}
Let $\textbf{d}=(d_2,\ldots,d_{k+1})$ and let $\mathfrak{S}=(G,G_{\mathcal{J}},\mathcal{J},\mathcal{P},\overrightarrow{H})$ be an oriented $(k,m,2t,\varepsilon,\varepsilon_{k+1},r,\textbf{d})$-regular setup.
Suppose that for each color cluster $C$, there are at least $(1-\alpha)t$ point clusters $Z$ such that $\{C,Z\}$ has relative  $(1,1)$-degree at least $\mu$ in $H$, then
\[
\frac{\mu^{2k+2}}{8}\binom{t}{k}^2\binom{t}{k-1}^{2k} t(t-1)\leq|\mathfrak{L}_{\overrightarrow{H}}|\leq\binom{t}{k}^2\binom{t}{k-1}^{2k} t(t-1).
\]
Let $\mathcal{F}$ be the $(k+1)$-complex corresponding to the down-closure of the $(1,k)$-graph $F$.
For each reduced gadget $L\in\mathfrak{L}_{\overrightarrow{H}}$ in $\overrightarrow{H}$, we have
\[
|\mathfrak{F}_L^{ext}|=(1\pm\beta)\left(\prod_{i=2}^{k+1}d_i^{e_i(\mathcal{F})}\right)m^{6k^2+4k+2}
\]
and
\[
|\mathfrak{F}^{ext}|=(1\pm\beta)\left(\prod_{i=2}^{k+1}d_i^{e_i(\mathcal{F})}\right)m^{6k^2+4k+2}|\mathfrak{L}_{\overrightarrow{H}}|.
\]
\end{lemma}
\begin{proof}
The lower bound of $\mathfrak{L}_{\overrightarrow{H}}$ can be done as follows.
Let $Y=(Y_0,Y_1,\ldots,Y_k)$,
$Y'=(Y_0,Z_1,\ldots,Z_k)$ $\in \overrightarrow{H}$ be the ordered core edge of $L$ and $W_i=\{W_{i,1},\ldots,W_{i,k-1}\}$, $W_i'=\{W_{i,1}',\ldots,W_{i,k-1}'\}$ for $i\in[k]$, be the peripheral sets, ordered such that $(Z_0,W_{i,1}',\ldots,W_{i,k-1}',Y_i)$ and $(Y_0,W_{i,1},\ldots,W_{i,k-1},Z_i)$ are consistent with $\overrightarrow{H}$.
We first choose $Y_0,Z_0$ arbitrarily, there are at least $t(t-1)$ choices.
For $(Y_1,\ldots,Y_k)$, there are at least $\mu\binom{t}{k}-\alpha t\binom{t}{k-1}\geq\mu\binom{t}{k}/2$ choices.
Similarly, for $(Z_1,\ldots,Z_k)$, there are at least $\mu\binom{t}{k}/2$ choices.
Furthermore, $W_i'$ and $W_i$ for $i\in[k]$ can be chosen in at least $\mu\binom{t}{k-1}$ ways for $i\in[k]$, but we need to delete the possible choices of intersecting reduced gadgets, whose number is at most $t(t-1)(2k^2)^2t^{2k^2-2}\leq(2k^2)^2t^{2k^2}$.
We have
\[
|\mathfrak{L}_{\overrightarrow{H}}|\geq\frac{\mu^{2k+2}}{4}\binom{t}{k}^2\binom{t}{k-1}^{2k} t(t-1)-(2k^2)^2t^{2k^2}\geq\frac{\mu^{2k+2}}{8}\binom{t}{k}^2\binom{t}{k-1}^{2k} t(t-1),
\]
since $1/t\ll\mu,1/k$.

While the upper bound is obvious.

We choose $\beta'$ such that $\varepsilon_{k+1}\ll\beta'\ll\beta,d_k,d_{k+1},1/k$.
By Lemma \ref{10.5} (with $\beta'$ in place of $\beta$), we obtain that
\[
(1-\beta)\left(\prod_{i=2}^{k+1}d_i^{e_i(\mathcal{F})}\right)m^{6k^2+4k+2}\leq(1-\beta')^2\left(\prod_{i=2}^{k+1}d_i^{e_i(\mathcal{F})}\right)m^{6k^2+4k+2}\leq(1-\beta')|\mathfrak{F}_L|\leq|\mathfrak{F}_L^{\rm ext}|,
\]
\[
|\mathfrak{F}_L^{\rm ext}|\leq|\mathfrak{F}_L|\leq(1+\beta')\left(\prod_{i=2}^{k+1}d_i^{e_i(\mathcal{F})}\right)m^{6k^2+4k+2}\leq(1+\beta)\left(\prod_{i=2}^{k+1}d_i^{e_i(\mathcal{F})}\right)m^{6k^2+4k+2}.
\]

Note that
\[
\mathfrak{F}^{\rm ext}=\bigcup_{L\in \mathfrak{L}_{\overrightarrow{H}}}\mathfrak{F}_L^{\rm ext},
\]
and the union is disjoint, the bounds of $|\mathfrak{F}^{\rm ext}|$ are easy to see.
\end{proof}

\begin{lemma}
\label{10.7}
Let $k,r,m,t\in \mathbb{N}$ and $d_2,\ldots,d_{k+1},\varepsilon,\varepsilon_{k+1},c,\nu,\theta,\mu$ be such that
\begin{equation}
\begin{split}
1/m&\ll1/r,\varepsilon\ll1/t,c,\varepsilon_{k+1},d_2,\ldots,d_k,\\
c&\ll d_2,\ldots,d_k,\\
1/t&\ll\varepsilon_{k+1}\ll d_{k+1}\leq1/k,\\
\varepsilon_{k+1}&\ll\nu\ll\theta\ll\mu\ll1/k.
\nonumber
\end{split}
\end{equation}
Let $\textbf{d}=(d_2,\ldots,d_{k+1})$ and let $\mathfrak{S}=(G,G_{\mathcal{J}},\mathcal{J},\mathcal{P},\overrightarrow{H})$ be an oriented $(k,m,2t,\varepsilon,\varepsilon_{k+1},r,\textbf{d})$-regular setup.
Suppose that for each color cluster $C$, there are at least $(1-\alpha)t$ point clusters $Z$ such that $\{C,Z\}$ has relative  $(1,1)$-degree at least $\mu$ in $H$.
For any point $v$ of $G$, color cluster $C$, there are at least $(1-\alpha)t$ point clusters $Z\in \mathcal{P}$ such that $|N_{\mathcal{J}}((v,C),\mu)\cap N_H(Z,C)|\geq\mu \binom{t}{k-1}$.
And for every $c\in [n]$, color cluster $C$, there are at least  $(1-\alpha)t$ point clusters $Z\in \mathcal{P}$ such that $|N_{\mathcal{J}}((c,Z),\mu)\cap N_H(C,Z)|\geq\mu \binom{t}{k-1}$.
Let $T\subseteq V$ be a $k$-set and $O\subseteq[n]$ be a $k$-set, we have
\[
|\mathfrak{F}_{(T,O)}^{\rm ext}|\geq\theta|\mathfrak{F}^{\rm ext}|.
\]
\end{lemma}
Given a $k$-subset $T=\{t_1,\ldots,t_k\}$ of $V$ and a $k$-subset $O=(o_1,\ldots,o_k)$, the family $\mathfrak{L}_{\overrightarrow{H}}$ and $\mu>0$, we define $\mathfrak{L}_{\overrightarrow{H},(T,O),\mu}$ of \emph{reduced $((T,O),\mu)$-absorbers} as the set of $(T,O)$-absorbers $Y\cup W_1\cup \cdots\cup W_k\cup Z_0\cup Z_1\cup\ldots\cup Z_k\cup W_1'\cup \cdots\cup W_k'$, where $W_i\subseteq N_{\mathcal{J}}((c_i,Z_i),\mu)$ and $W_i'\subseteq N_{\mathcal{J}}((t_i,Z_0),\mu)$ for $i\in[k]$.

\begin{claim}
\label{10.8}
Let $k,r,m,t\in \mathbb{N}$ and $d_2,\ldots,d_{k+1},\varepsilon,\varepsilon_{k+1},c,\nu,\theta,\mu$ be such that
\begin{equation}
\begin{split}
1/m&\ll1/r,\varepsilon\ll1/t,c,\varepsilon_{k+1},d_2,\ldots,d_k,\\
c&\ll d_2,\ldots,d_{k+1},\\
1/t&\ll\varepsilon_{k+1}\ll d_{k+1}\leq1/k,\\
\varepsilon_{k+1}&\ll\nu\ll\theta\ll\mu\ll1/k,\\
\alpha&\ll\mu.
\nonumber
\end{split}
\end{equation}
Let $\textbf{d}=(d_2,\ldots,d_{k+1})$ and let $\mathfrak{S}=(G,G_{\mathcal{J}},\mathcal{J},\mathcal{P},\overrightarrow{H})$ be an oriented $(k,m,2t,\varepsilon,\varepsilon_{k+1},r,\textbf{d})$-regular setup.
Suppose that for each color cluster $C$, there are at least $(1-\alpha)t$ point clusters $Z$ such that $\{C,Z\}$ has relative  $(1,1)$-degree at least $\mu$ in $H$.
For any point $v$ of $G$, color cluster $C$, there are at least $(1-\alpha)t$ point clusters $Z\in \mathcal{P}$ such that $|N_{\mathcal{J}}((v,C),\mu)\cap N_H(Z,C)|\geq\mu \binom{t}{k-1}$.
And for every $c\in [n]$, color cluster $C$, there are at least  $(1-\alpha)t$ point clusters $Z\in \mathcal{P}$ such that $|N_{\mathcal{J}}((c,Z),\mu)\cap N_H(C,Z)|\geq\mu \binom{t}{k-1}$.
Let $T\subseteq V$ be a $k$-set and $O\subseteq[n]$ be a $k$-set, we have
\[
|\mathfrak{L}_{\overrightarrow{H},(T,O),\mu}|\geq\theta|\mathfrak{L}_{\overrightarrow{H}}|.
\]
\end{claim}
\begin{proof}
Let $T=\{t_1,\ldots,t_k\}$ and $O=(o_1,\ldots,o_k)$.
Since $H$ has minimum relative $(1,1)$-degree at least $\mu$, there are at least $\mu t\binom{t}{k}-t\alpha t\binom{t}{k-1}\geq\mu t\binom{t}{k}/2$ choices for $Y$.
Besides, there are at least $t-1$ choices for $Z_0$.
For $(Z_1,\ldots,Z_k)$, there are at least $\mu\binom{t}{k}/2-k^2\binom{t}{k-1}\geq\mu\binom{t}{k}/3$ choices.
Each $W_i$ is chosen from $N_{\mathcal{J}}((o_i,Z_i),\mu)\cap N_H(Y_0,Z_i)$ for $i\in[k]$, thus, $W_i$ can be chosen in at least $\mu\binom{t}{k-1}-(k-1)((i-1)(k-1)+2k)\binom{t}{k-2}\geq\mu\binom{t}{k-1}/2$ ways for $i\in[k]$, since there are at most $(k-1)((i-1)(k-1)+2k)\binom{t}{k-2}$ choices for $W_i'$ which intersects with $Y\setminus\{Y_0\}, Z_1,\ldots,Z_k, W_1,\ldots,W_{i-1}$.

And each $W_i'$ is chosen from $N_{\mathcal{J}}((t_i,Z_0),\mu)\cap N_H(Y_i,Z_0)$ for $i\in[k]$.
Similarly, there are at least $(\mu/2)\binom{t}{k-1}$ possible choices for each $W_i'$ for $i\in[k]$.
Thus, the number of reduced $((T,O),\mu)$-absorbers is at least
\[
\frac{\mu t}{2}\binom{t}{k}(t-1)\frac{\mu }{3} \binom{t}{k}\left(\frac{\mu}{2}\binom{t}{k-1}\right)^{2k}\geq\theta\binom{t}{k}^2\binom{t}{k-1}^{2k} t(t-1)\geq\theta|\mathfrak{L}_{\overrightarrow{H}}|
\]
since $\theta\ll\mu$.
\end{proof}
\begin{claim}
\label{10.9}
Let $k,r,m,t\in \mathbb{N}$ and $d_2,\ldots,d_{k+1},\varepsilon,\varepsilon_{k+1},c,\nu,\theta,\mu$ be such that
\begin{equation}
\begin{split}
1/m&\ll1/r,\varepsilon\ll1/t,c,\varepsilon_{k+1},d_2,\ldots,d_k,\\
c&\ll d_2,\ldots,d_k,\\
1/t&\ll\varepsilon_{k+1}\ll d_{k+1}\leq1/k,\\
\varepsilon_{k+1}&\ll\nu\ll\theta\ll\mu\ll1/k.
\nonumber
\end{split}
\end{equation}
Let $\textbf{d}=(d_2,\ldots,d_{k+1})$ and let $\mathfrak{S}=(G,G_{\mathcal{J}},\mathcal{J},\mathcal{P},\overrightarrow{H})$ be an oriented $(k,m,2t,\varepsilon,\varepsilon_{k+1},r,\textbf{d})$-regular setup.
Let $T\subseteq V$ and $O\subseteq[n]$ be $k$-sets and let $L\in\mathfrak{L}_{\overrightarrow{H}}$ be a reduced $((T,O),\mu)$-gadget in $\overrightarrow{H}$.
We have
\[
|\mathfrak{F}_L\cap\mathfrak{F}_{(T,O)}|\geq\theta|\mathfrak{F}_L|.
\]
\end{claim}
\begin{proof}
Let $T=\{t_1,\ldots,t_k\}$ and $O=\{o_1,\ldots,o_k\}$, $L=Y\cup W_1\cup \cdots\cup W_k\cup Z_0\cup Z_1\cup\ldots\cup Z_k\cup W_1'\cup \cdots\cup W_k'$ where $W_i=\{W_{i,1},\ldots,W_{i,k-1}\}$ and $W_i'=\{W_{i,1}',\ldots,W_{i,k-1}'\}$.
Choose $P_i, Q_i$ in $W_i$ and $P_i', Q_i'$ in $W_i'$, let $\mathcal{Q}_{Z_i,W_i}$ be the set of $k$-uniform tight paths $(b_i,v_1,\ldots,v_{2k-2})$ such that $b_i\in Z_i$, $v_{\ell},v_{\ell+k-1}\in W_{i,\ell}$ for $i,j\in[k]$, $\ell\in[k-1]$ and its down-closure is in $\mathcal{J}$.
Let $\mathcal{Q}_{o_i,(Z_i,W_i)}\subseteq\mathcal{Q}_{Z_i,W_i}$ be the set of those paths whose edges in the $k$-th level are in $N_G(o_i)$.
Note that $F$ is the absorbing gadget for $(T,O)$.
Let $\mathcal{F}$ be the down-closure of $F$.
Since $L$ is a reduced $(T,\mu)$-gadget, we have $W_i\in N_H(Y_0,Z_i)\cap N_{\mathcal{J}}((o_i,Z_i),\mu)$, thus $|N_G((o_i,Z_i),\mathcal{J}_{W_i})|\geq\mu|\mathcal{J}_{W_i}|$.
By Lemma \ref{color} with $S$ being the set of $k$-sets where each $k$-set consists of $k-1$ points from $N_G((o_i,Z_i),\mathcal{J}_{W_i})$ and one point from $Z_i$, we have
\[
|\mathcal{Q}_{o_i,(Z_i,W_i)}|\geq\frac{1}{2}\left(\frac{\mu}{8k}\right)^{k+1}|\mathcal{Q}_{Z_i,W_i}|.
\]

Let $\mathcal{Q}_{Z_0,W_i'}$ be the set of $k$-uniform sequentially paths $(c_1',\ldots,c_k',v_1',\ldots,v_{2k-2}')$ such that $c_j'\in Z_0$, $v_{\ell}',v_{\ell+k-1}'\in W_{i,\ell}'$ for $i,j\in[k]$, $\ell\in[k-1]$ and its down-closure is in $\mathcal{J}$.
Let $\mathcal{Q}_{t_i,(Z_0,W_i')}\subseteq\mathcal{Q}_{Z_0,W_i'}$ be the set of those paths whose edges in the $k$-th level are in $N_G(t_i)$.
Since $L$ is a reduced $((T,O),\mu)$-gadget, we have $W_i'\in N_H(Z_0,Y_i)\cap N_{\mathcal{J}}((t_i,Z_0),\mu)$, thus $|N_G((t_i,Z_0),\mathcal{J}_{W_i'})|\geq\mu|\mathcal{J}_{W_i'}|$.
By Lemma \ref{9.3} with $S$ being the set of $k$-sets where each $k$-set consists $k-1$ points from $N_G((t_i,Z_0),\mathcal{J}_{W_i'})$ and one color from $Z_0$, we have
\begin{equation}
|\mathcal{Q}_{t_i,(Z_0,W_i')}|\geq\frac{1}{2}\left(\frac{\mu}{8k}\right)^{k+1}|\mathcal{Q}_{Z_0,W_i'}|.
\nonumber
\end{equation}

Let $\phi:V(F)\rightarrow V(L)$ be the homomorphism which labels the copies of $F$ in $\mathfrak{F}_L$.
Set $Z=\{b_1,\ldots,b_k\}\cup \bigcup_{i=1}^k(V(P_i)\cup V(Q_i))\cup\bigcup_{i=1}^k(C_i'\cup V(P_i')\cup V(Q_i'))$.
Thus, $|Z|=5k^2-3k$.
Let $\mathcal{Z}=\mathcal{F}[Z]$ be the induced subcomplex of $\mathcal{F}$ in $Z$.
Note that $\mathcal{Z}$ consists of  $k$ vertex-disjoint $k$-uniform tight paths of length $2k-1$ where the $i$-th path lies in $\mathcal{Q}_{o_i,(Z_i,W_i)}$ and $k$ vertex-disjoint $k$-uniform sequentially paths of length $2k-2$ where the $i$-th path lies in $\mathcal{Q}_{t_i,(Z_0,W_i')}$.
Let $\mathcal{G}=\mathcal{J}\cup G_{\mathcal{J}}$ and $\mathcal{Z}_{\mathcal{G}}$ be the set of labelled partition-respecting copies of $\mathcal{Z}$ in $\mathcal{G}$.
Let $\beta_1$ be such that $\varepsilon\ll\beta_1\ll d_2,\ldots,d_k,\varepsilon_{k+1}$ and define $d_{\mathcal{Z}}=\prod_{i=2}^kd_i^{e_i(\mathcal{Z})}$.
By Lemma \ref{dense counting}, we have
\begin{equation}
|\mathcal{Z}_{\mathcal{G}}|=\prod_{i=1}^k|\mathcal{Q}_{Z_i,W_i}||\mathcal{Q}_{Z_0,W_i'}|=(1\pm\beta_1)d_{\mathcal{Z}}m^{5k^2-3k}.
\nonumber
\end{equation}

Let $\mathcal{Z}_{(T,O),\mathcal{G}}\subseteq\mathcal{Z}_{\mathcal{G}}$ be the labelled partition-respecting copies of $\mathcal{Z}$ absorbing $(T,O)$, thus we have
\begin{equation}
|\mathcal{Z}_{(T,O),\mathcal{G}}|\geq \prod_{i=1}^k|\mathcal{Q}_{o_i,(Z_i,W_i)}||\mathcal{Q}_{t_i,(Z_0,W_i')}|\geq\left(\frac{1}{2}\left(\frac{\mu}{8k}\right)^{k+1}\right)^{2k}\prod_{i=1}^k|\mathcal{Q}_{Z_i,W_i}||\mathcal{Q}_{Z_0,W_i'}|\geq3\theta|\mathcal{Z}_{\mathcal{G}}|,
\nonumber
\end{equation}
since $\theta\ll\mu,1/k$.

Let $\beta_2$ be such that $\varepsilon_{k+1}\ll\beta_2\ll \theta,d_{k+1},1/k$ and define $d_{\mathcal{F}-\mathcal{Z}}=\prod_{i=2}^{k+1}d_i^{e_i(\mathcal{F})-e_i(\mathcal{Z})}$.
Let $I\subseteq\mathcal{Z}_{\mathcal{G}}$ be the set of labelled partition-respecting copies of $\mathcal{Z}$ which are not extensible to $(1\pm\beta_2)d_{\mathcal{F}-\mathcal{Z}}m^{k^2+7k+2}$ labelled partition-respecting copies of $\mathcal{F}$ in $\mathcal{G}$.
By Lemma \ref{extension}, we have
\[
|I|\leq\beta_2|\mathcal{Z}_{\mathcal{G}}|\leq\theta|\mathcal{Z}_{\mathcal{G}}|,
\]
since $\beta_2\ll \theta$.
By Lemma \ref{10.5}, we have
\[
|\mathfrak{F}_L|=(1\pm\beta_2)d_{\mathcal{F}-\mathcal{Z}}d_{\mathcal{Z}}m^{6k^2+4k+2},
\]
since $\varepsilon_{k+1}\ll\beta_2\ll \theta,d_{k+1},1/k$.

Note that a labelled partition-respecting copy of $\mathcal{F}$ in $\mathcal{G}$ containing a $Z\in\mathcal{Z}_{(T,O),\mathcal{G}}$ yields exactly one gadget in $\mathfrak{F}_L\cap\mathfrak{F}_{(T,O)}$, we have
\begin{equation}
\begin{split}
|\mathfrak{F}_L\cap\mathfrak{F}_{(T,O)}|&\geq|\mathcal{Z}_{(T,O),\mathcal{G}}\setminus I|(1-\beta_2)d_{\mathcal{F}-\mathcal{Z}}m^{k^2+7k+2}\\
&\geq(|\mathcal{Z}_{(T,O),\mathcal{G}}|-|I|)(1-\beta_2)d_{\mathcal{F}-\mathcal{Z}}m^{k^2+7k+2}\\
&\geq2\theta|\mathcal{Z}_{\mathcal{G}}|(1-\beta_2)d_{\mathcal{F}-\mathcal{Z}}m^{k^2+7k+2}\\
&\geq2\theta(1-\beta_2)(1-\beta_1)d_{\mathcal{Z}}m^{5k^2-3k}d_{\mathcal{F}-\mathcal{Z}}m^{k^2+7k+2}\\
&\geq2\theta(1-2\beta_2)d_{\mathcal{Z}}d_{\mathcal{F}-\mathcal{Z}}m^{6k^2+4k+2}\\
&\geq2\theta\frac{1-2\beta_2}{1+\beta_2}|\mathfrak{F}_L|\\
&\geq\theta|\mathfrak{F}_L|,
\nonumber
\end{split}
\end{equation}
since $\beta_2\ll\theta$.
\end{proof}
\begin{proof}[Proof of Lemma \ref{10.7}]
Let $\theta\ll\theta'\ll\mu$.
By Claim \ref{10.9} with $\theta'$, we have for every reduced $((T,O),\mu)$-gadget $L\in\mathfrak{L}_{\overrightarrow{H}}$,
\[
|\mathfrak{F}_L\cap\mathfrak{F}_{(T,O)}|\geq\theta'|\mathfrak{F}_L|.
\]
Let $\beta$ be such that $\varepsilon_{k+1}\ll\beta\ll d_{k+1},\theta'$, by Lemma \ref{10.5} with $\theta'$, we have $|\mathfrak{F}_L\setminus\mathfrak{F}_L^{\rm ext}|\leq\beta|\mathfrak{F}_L|\leq\theta'|\mathfrak{F}_L|/2$.
Thus,
\[
|\mathfrak{F}_{(T,O)}^{\rm ext}\cap\mathfrak{F}_L|\geq|\mathfrak{F}_L\cap\mathfrak{F}_{(T,O)}|-|\mathfrak{F}_L\setminus\mathfrak{F}_L^{\rm ext}|\geq\frac{\theta'}{2}|\mathfrak{F}_L|.
\]
By Claim \ref{10.8} with $\theta'$ and Lemma \ref{10.6}, we have $|\mathfrak{L}_{\overrightarrow{H},(T,O),\mu}|\geq\theta'|\mathfrak{L}_{\overrightarrow{H}}|$ and
\[
|\mathfrak{F}_{(T,O)}^{\rm ext}|\geq\sum_{L\in\mathfrak{L}_{\overrightarrow{H},(T,O),\mu}}|\mathfrak{F}_{(T,O)}^{\rm ext}\cap\mathfrak{F}_L|\geq\frac{\theta'}{2}\sum_{L\in\mathfrak{L}_{\overrightarrow{H},(T,O),\mu}}|\mathfrak{F}_L|\geq\theta|\mathfrak{F}^{\rm ext}|.
\]
\end{proof}

\subsection{Absorbing Lemma}
\label{absorbing lemma}
\begin{lemma}
\label{10.10}
Let $k,r,m,t\in \mathbb{N}$ and $d_2,\ldots,d_{k+1},\varepsilon,\varepsilon_{k+1},c,\nu,\theta,\mu,\alpha,\zeta$ be such that
\begin{equation}
\begin{split}
1/m&\ll1/r,\varepsilon\ll1/t,\zeta,\varepsilon_{k+1},d_2,\ldots,d_k,\\
\zeta&\ll c\ll d_2,\ldots,d_k,\\
1/t&\ll\varepsilon_{k+1}\ll d_{k+1},\nu\leq1/k,\\
c&\ll\varepsilon_{k+1}\ll\alpha\ll\theta\ll\mu\ll1/k.
\nonumber
\end{split}
\end{equation}
Let $\textbf{d}=(d_2,\ldots,d_{k+1})$ and let $\mathfrak{S}=(G,G_{\mathcal{J}},\mathcal{J},\mathcal{P},\overrightarrow{H})$ be an oriented $(k,m,2t,\varepsilon,\varepsilon_{k+1},r,\textbf{d})$-regular setup.
Suppose that $V(G)=[n]\cup V$ where $|V|=n\leq(1+\alpha)mt$ and $V(H)=[t]\cup V'$ where $|V'|=t$.
Suppose that for each color cluster $C$, there are at least $(1-\alpha)t$ point clusters $Z$ such that $\{C,Z\}$ has relative  $(1,1)$-degree at least $\mu$ in $H$.
For any point $v$ of $G$, color cluster $C$, there are at least $(1-\alpha)t$ point clusters $Z\in \mathcal{P}$ such that $|N_{\mathcal{J}}((v,C),\mu)\cap N_H(Z,C)|\geq\mu \binom{t}{k-1}$.
And for every $c\in [n]$, color cluster $C$, there are at least  $(1-\alpha)t$ point clusters $Z\in \mathcal{P}$ such that $|N_{\mathcal{J}}((c,Z),\mu)\cap N_H(C,Z)|\geq\mu \binom{t}{k-1}$. Then there exists a family $\mathfrak{F}''$ of pairwise disjoint $\mathfrak{S}$-gadgets which are $(c,\nu)$-extensible with the following properties.
\begin{itemize}
  \item[(1)] $|\mathfrak{F}''|\leq\zeta m,$
  \item[(2)] $|\mathfrak{F}''\cap\mathfrak{F}_{(T,O)}^{\rm ext}|\geq\zeta\theta m$ for any $k$-subset $T$ of $V$ and $k$-subset $O$ of $[n]$,
  \item[(3)] $V(\mathfrak{F}'')$ is $(2(k+1)\zeta/t)$-sparse in $\mathcal{P}$.
\end{itemize}
\end{lemma}
\begin{proof}
Let $\beta>0$ be such that $\varepsilon_{k+1}\ll\beta\ll d_{k+1}$.
Let $F$ be the $(1,k)$-graph as in Definition \ref{10.2} and let $\mathcal{F}$ be the $(k+1)$-complex generated by its down-closure.
Let $d_F=\prod_{i=2}^{k+1}d_i^{e_i(\mathcal{F})}$.
By Lemma \ref{10.6}, we have
\begin{equation}
|\mathfrak{F}^{\rm ext}|\leq(1+\beta)d_Fm^{6k^2+4k+2}\binom{t}{k}^2\binom{t}{k-1}^{2k}t (t-1)\leq d_Fm^{6k^2+4k+2}t^{2k^2+2},
\nonumber
\end{equation}
\begin{equation}
\begin{aligned}
|\mathfrak{F}^{\rm ext}|&\geq\frac{\mu^{k+1}}{2}(1-\beta)d_Fm^{6k^2+4k+2}\binom{t}{k}^2\binom{t}{k-1}^{2k}t(t-1)\\
&\geq \frac{\mu^{k+1}}{2^{k+2}k^{2k}(k-1)^{2k^2}}d_Fm^{6k^2+4k+2}t^{2k^2+2}\\
&\geq6\theta^{1/2}d_Fm^{6k^2+4k+2}t^{2k^2+2},
\nonumber
\end{aligned}
\end{equation}
since $1/t\ll\varepsilon_{k+1}\ll\beta\ll d_{k+1}\ll1/k$ and $\theta\ll\mu,1/k$.
By Lemma \ref{10.6}, for each reduced gadget  $L\in\mathfrak{L}_{\overrightarrow{H}}$ in $\overrightarrow{H}$, we have
\[
|\mathfrak{F}_L^{ext}|\leq2d_Fm^{6k^2+4k+2}.
\]
By Lemma \ref{10.7} with $\theta^{1/2}$, for any $k$-set $T\subseteq V$ and any $k$-set $O\subseteq[n]$, we have
\[
|\mathfrak{F}_{(T,O)}^{\rm ext}|\geq\theta^{1/2}|\mathfrak{F}^{\rm ext}|\geq6\theta d_Fm^{6k^2+4k+2}t^{2k^2+2}.
\]

Choose a family $\mathfrak{F}'$ from $\mathfrak{F}^{\rm ext}$ by including each $\mathfrak{S}$-gadget independently at random with probability
\[
p=\frac{\zeta m}{2d_Fm^{6k^2+4k+2}t^{2k^2+2}}.
\]
Note that $|\mathfrak{F}'|$, $|\mathfrak{F}'\cap\mathfrak{F}_{(T,O)}^{\rm ext}|$ are binomial random variables, for any $k$-set $T\subseteq V$ and any $k$-set $O\subseteq[n]$, we have
\[
\mathbb{E}|\mathfrak{F}'|=p|\mathfrak{F}^{\rm ext}|\leq\frac{\zeta m}{2},
\]
\[
\mathbb{E}|\mathfrak{F}'\cap\mathfrak{F}_{(T,O)}^{\rm ext}|=p|\mathfrak{F}_{(T,O)}^{\rm ext}|\geq3\theta\zeta m.
\]

For each $Z\in \mathcal{P}$,
note that $Z$ exists in at most $t^{2k^2+1}$ reduced gadgets, thus, there are at most $2d_Fm^{6k^2+4k+2}t^{2k^2+1}$ $\mathfrak{S}$-gadgets with vertices in $Z$.
Note that each $\mathfrak{S}$-gadget contains at most $k^2+2k+2$ vertices in a cluster.
Hence, for each cluster $Z\in \mathcal{P}$, we have
\[
\mathbb{E}|V(\mathfrak{F}')\cap Z|\leq2(k^2+2k+2)d_Fm^{6k^2+4k+2}t^{2k^2+1}p=\frac{(k^2+2k+2)\zeta m}{t}.
\]
By Proposition \ref{chernoff}, with probability $1-o(1)$, the family $\mathfrak{F}'$ satisfies the following properties.
\[
|\mathfrak{F}'|\leq2\mathbb{E}|\mathfrak{F}'|\leq \zeta m,
\]
\[
|\mathfrak{F}'\cap\mathfrak{F}_{(T,O)}^{\rm ext}|\geq2\theta\zeta m,
\]
\[
|V(\mathfrak{F}')\cap Z|\leq\frac{2(k^2+k+1)\zeta m}{t}
\]
for any $k$-set $T\subseteq V$, $k$-set $O\subseteq [n]$ and cluster $Z\in \mathcal{P}$.
We say that two $\mathfrak{S}$-gadgets are \emph{intersecting} if they share at least one vertex.
Note that there at most $(2k^2+2)^2t^{4k^2+3}$ pairs of intersecting reduced gadgets.
Hence, there are at most $(6k^2+4k+2)^2m^{12k^2+8k+1}(2k^2+2)^2t^{4k^2+3}$ pairs of intersecting $\mathfrak{S}$-gadgets.
We can bound the expected number of pairs of intersecting $\mathfrak{S}$-gadgets by
\begin{equation}
\begin{split}
&(6k^2+4k+2)^2m^{12k^2+8k+3}(2k^2+2)^2t^{4k^2+3}p^2\\
&=\frac{\zeta^2(6k^2+4k+2)^2(2k^2+2)^2m}{4d_F^2t}\leq\frac{\zeta\theta m}{2},
\end{split}
\nonumber
\end{equation}
since $\zeta\ll d_2,\ldots,d_{k+1},\theta,1/k$.
Using Markov's inequality, we derive that with probability at least $1/2$, $\mathfrak{F}'$ contains at most $\zeta\theta m$ pairs intersecting $\mathfrak{S}$-gadgets.
Remove one gadget from each intersecting pair in such a family and remove gadgets that are not absorbing for any $(T,O)$ where $T\subseteq V$, $O\subseteq[n]$ and $|T|=|O|$.
We obtain a subfamily $\mathfrak{F}''$, satisfying the following properties.
\begin{itemize}
  \item[(1)] $|\mathfrak{F}''|\leq \zeta m$,
  \item[(2)] $|\mathfrak{F}''\cap\mathfrak{F}_{(T,O)}^{\rm ext}|\geq\theta\zeta m$,
  \item[(3)] $V(\mathfrak{F}'')$ is $(2(k^2+k+1)\zeta/t)$-sparse in $\mathcal{P}$,
\end{itemize}
as desired.
\end{proof}
\begin{proof}[The proof of Lemma \ref{absorption lemma}]
Since $G$ has minimum relative $(1,1)$-degree at least $\delta+\mu$ and $\mathfrak{S}$ is a representative setup.
For any $v\in V$ and any color cluster $C$, we have
\[
|N_{\mathcal{J}}((v,C),\frac{\mu}{3})|\geq(\delta+\frac{\mu}{4})\binom{t}{k-1}.
\]
For any $c\in[n]$ and any point cluster $Z$, we have
\[
|N_{\mathcal{J}}((c,Z),\frac{\mu}{3})|\geq(\delta+\frac{\mu}{4})\binom{t}{k-1}.
\]
by Lemma \ref{9.1}.
Let $\zeta>0$ with $1/r,\varepsilon\ll\zeta\ll c$ and let $\theta>0$ with $\eta\ll \theta\ll\mu,1/k$ and $M:=\lceil\eta t/(\theta\zeta)\rceil$.
Firstly, we need the following claim.
\begin{claim}
For each $j\in [0,M]$, and any $S\subseteq V$ of size at most $j\theta\zeta n/t$ divisible by $k$ and any $O\subseteq[n]$ of size $|S|$, there is a sequentially path $P_j\subseteq G$ such that the following holds.
\begin{itemize}
  \item[(i)] $P_j$ is $(S,O)$-absorbing in $G$,
  \item[(ii)] $P_j$ is $(c,\nu)$-extensible and consistent with $\overrightarrow{H}$,
  \item[(iii)] $V(P_j)$ is $(100k^3j\zeta/t)$-sparse in $\mathcal{P}$ and $V(P_j)\cap T_j=\emptyset$, where $T_j$ denotes the connection set of $P_j$.
\end{itemize}
\end{claim}
\begin{proof}[Proof of the claim]
Take $P_0$ to be the empty path and $P_j$ satisfy the above conditions for $j\in[0,M)$.

Select a subset $Z'\subseteq Z\setminus V(P_j)$ of size $m'=(1-\lambda)m$ since $100k^3j\zeta/t\leq(2\eta t/(\zeta\theta))(100k^3\zeta/t)\leq\lambda$ which follows from $\zeta\ll c\ll\eta\ll\lambda,\theta$.
Also, since $n\leq(1+\alpha)mt$, we have $m'\geq n/(2t)$.
Let $\mathcal{P}'=\{Z'\}_{Z\in \mathcal{P}}$,  $\mathcal{J}'=\mathcal{J}[V(\mathcal{P}')]$ and
$G_{\mathcal{J}'}'=G_{\mathcal{J}}[V(\mathcal{P}')]$ .
By lemma \ref{restriction}, $\mathfrak{S}':=(G',G_{\mathcal{J}'}',\mathcal{J}',\mathcal{P}',H)$ is a $(k,m',2t,\sqrt{\varepsilon},\sqrt{\varepsilon_{k+1}},r,\textbf{d})$-regular setup.

By Lemma \ref{9.2}, for every $v\in V$ and color cluster $C$, we have
\[
|N_{\mathcal{J}'}((v,C),\mu/6)|\geq|N_{\mathcal{J}}((v,C),\mu/3)|\geq (\delta+\mu/4)\binom{t}{k-1},
\]
and for every $o\in [n]$ and point cluster $Z$, we have
\[
|N_{\mathcal{J}'}((o,Z),\mu/6)|\geq|N_{\mathcal{J}}((o,Z),\mu/3)|\geq (\delta+\mu/4)\binom{t}{k-1},
\]
Thus, we obtain that for every $v\in V$, $o\in [n]$, color cluster $C$, there are at least $(1-\alpha)t$ point clusters $Z\in \mathcal{P}$, we have \[|N_{\mathcal{J}}((v,C),\mu/6)\cap N_H(Z,C)|\geq\frac{\mu}{5} \binom{t}{k-1},\]
and
\[|N_{\mathcal{J}}((o,Z),\mu/6)\cap N_H(C,Z)|\geq\frac{\mu}{5} \binom{t}{k-1}.
\]
By Lemma \ref{10.10} with $4c$ instead of $c$, $2\zeta$ instead of $\zeta$, we obtain a set $\mathcal{A}'$ of pairwise-disjoint $\mathfrak{S}'$-gadgets which are $(4c,\nu)$-extensible and such that
\begin{itemize}
  \item[(1)] $|\mathcal{A}'|\leq2\zeta m',$
  \item[(2)] $|\mathcal{A}'\cap\mathfrak{F}_{(T,O)}|\geq2\zeta\theta m'$ for any $k$-subset of $V$,
  \item[(3)] $V(\mathcal{A}')$ is $(4(k^2+k+1)\zeta/t)$-sparse in $\mathcal{P}'$.
\end{itemize}
Next, we would connect all paths of absorbing gadgets in $\mathcal{A}'$ and $P_j$ to obtain $P_{j+1}$.
By Definition \ref{10.2}, there are $2(k+1)$ pairwise disjoint sequentially paths  in each $\mathfrak{S}'$-gadget in $\mathcal{A}'$ which are $(4c,\nu)$-extensible in $\mathfrak{S}'$.
Let $\mathcal{A}$ be the union of all such sequentially paths of all gadgets of $\mathcal{A}'$ and $P_j$.
Set $T_{j+1}=V(G)\setminus V(\mathcal{A})$, it is obvious that $\mathcal{A}$ is a set of pairwise disjoint sequentially paths in $G$ such that
\begin{itemize}
  \item[(1')] $|\mathcal{A}|\leq4(k+1)\zeta m'+1,$
  \item[(2')] $V(\mathcal{A})$ is $(100k^3j\zeta/t+4(k^2+k+1)\zeta/t)$-sparse in $\mathcal{P}$ and $V(\mathcal{A})\cap T_{j+1}=\emptyset$,
  \item[(3')] every path in $\mathcal{A}\setminus \{P_j\}$ is $(2c,\nu,T_{j+1})$-extensible in $\mathfrak{S}$ and consistent with $\overrightarrow{H}$.
      $P_j$ is $(c,\nu,T_{j+1})$-extensible in $\mathfrak{S}$ and consistent with $\overrightarrow{H}$.
\end{itemize}
Note that (1') follows from (1) and the addition of $P_j$.
(2') follows from (iii), (3) and the definition of $T_{j+1}$.
(3') follows from (ii) and (3) since $4(k^2+k+1)\zeta m/t\leq2c m$.
In particular, $P_j$ is $(c,\nu)$-extensible by (ii) while all other paths go from $(4c,\nu)$-extensible in $\mathfrak{S}'$ to $(2c,\nu)$-extensible in $\mathfrak{S}$.
The consistency with $\overrightarrow{H}$ is given by the consistency of $P_j$ and the definition of $\mathfrak{S}'$-gadgets.

By Lemma \ref{8.5}, we obtain a sequentially path $P_{j+1}$ with the following properties.
\begin{itemize}
  \item[(A)] $P_{j+1}$ contains every path of $\mathcal{A}$,
  \item[(B)] $P_{j+1}$ starts and ends with two paths different from $P_j$,
  \item[(C)] $V(P_{j+1})\setminus V(\mathcal{A})\subseteq V(\mathcal{P}')$,
  \item[(D)] $V(P_{j+1})\setminus V(\mathcal{A})$ intersects in at most $10k^2\mathcal{A}_Z+t^{2t+3k+2}$ vertices with each cluster $Z\in \mathcal{P}$, where $\mathcal{A}_Z$ denotes the number of paths of $\mathcal{A}$ that intersect with $Z$.
\end{itemize}

We claim that $P_{j+1}$ satisfies (i)-(iii).
First, we prove (iii).
Note that for every cluster $Z\in \mathcal{P}$, the number of paths of $\mathcal{A}$ that intersect with $Z$ is bounded by $4(k+1)\zeta m/t+1$.
(D) implies that $V(P_{j+1})\setminus V(\mathcal{A})$ intersects in at most $100k^3\zeta m/t$ vertices with each cluster $Z\in \mathcal{P}$.
Together with (iii), it follows that $\mathcal{A}$ is $(100k^3(j+1)\zeta/t)$-sparse in $\mathcal{P}$.

Next, we want to prove (ii), $V(P_{j+1})\setminus V(\mathcal{A})$ intersects in at most $100k^3\zeta m/t\leq cm/4$ vertices with each cluster $Z\in \mathcal{P}$, since $\zeta\ll c$.
Also, we have $V(\mathcal{A})\cap T_{j+1}=\emptyset$.
Hence, we obtain (ii) after deleting the vertices of $P_{j+1}$ from $T_{j+1}$.
After the deletion, we go from $(2c,\nu)$-extensible in (3') to $(c,\nu)$-extensible.
It is crucial that $P_{j+1}$ starts and ends with two paths different from $P_j$ by (B).

Finally, we claim that $P_{j+1}$ is $(S,O)$-absorbing in $G$ for any $S\subseteq V$ of size divisible by $k$ and at most $(j+1)\zeta\theta n/t$ and any $O\subseteq[n]$ of size $|S|$.
Partition $S$ into two sets $S_1$ and $S_2$ such that both $|S_1|,|S_2|$ are divisible by $k$ and $S_1$ is maximal such that $|S_1|\leq j\zeta\theta n/t$.
Partition $O$ into two sets $O_1$ and $O_2$ such that $|O_1|=|S_1|$ and $|O_2|=|S_2|$.
Since $P_j$ is $(S',O')$-absorbing in $G$ for any set $S'\subseteq V$ of size at most $(j\zeta\theta n/t)$ and $|O'|=|S'|$, there exists a path $P_j'$ with the same endpoints as $P_j$ such that $I(P_j')=S_1\cup I(P_j)$ and $C(P_j')=O_1\cup C(P_j)$, besides, $P_j$ is a subpath of $P_{j+1}$.
So it remains to absorb $S_2$.
By the choice of $S_1$, we have $|S_2|\leq\zeta\theta n/t+k\leq2\zeta^3n/t\leq2(1+\alpha)\zeta^3m\leq5\zeta^3m/2$.
Therefore, we can partition $S_2$ and $O_2$ into $\ell\leq 5\zeta^3m/(2k)\leq 2\zeta\theta m'$ sets of size $k$ each, let $D_1,\ldots,D_{\ell}$ and $R_1,\ldots,R_{\ell}$ be those sets.
By (2), we have $|\mathfrak{F}_{(D_i,R_i)}\cap \mathcal{A}'|\geq\ell$.
Thus, we can associate each $(D_i, R_i)$ with a different gadget $F_i\in \mathcal{A}'$ for each $i\in[\ell]$.
Each $F_i$ yields a collection of $2(k+1)$ sequentially paths $P_{i,1},\ldots,P_{i,2(k+1)}$ and we can replace those paths with a collection of different paths with the same endpoints.
Since $P_j$ and each $P_{i,u}$, $i\in[\ell], u\in[2(k+1)]$, are subpaths of $P_{j+1}$, the sequentially path $P_{j+1}'$ has the same endpoints with $P_{j+1}$.
Also, $P_{j+1}'$ is exactly $(C(P_{j+1})\cup O,I(P_{j+1})\cup S)$.
\end{proof}
To finish, note that $P_M$ and $C_M$ has the desired properties.
By the choice of $M=\lceil\eta t/(\zeta\theta)\rceil$, we have $M\zeta\theta/t\geq\eta$, so $P_M$ with $C_M$ is $\eta$- absorbing in $G$.
Moreover, since $M(100k^3\zeta/t)\leq200k^2\eta/\theta\leq\lambda$ and $\eta\ll\lambda$, $V(P_M)$ is $\lambda$-sparse in $\mathcal{P}$.
\end{proof}

\section{Concluding Remarks}

Inspired by a series of very recent successes on rainbow matchings~\cite{2020A,LU,Hongliang2018ON,2020C}, rainbow Hamilton cycles~\cite{CHWWY,2019Rainbow,MR4171383} and rainbow factors~\cite{2020B,MR4055023,factor}, we suspect the threshold for a rainbow spanning subgraph in (hyper)graph system is
asymptotically same
with the threshold for a spanning subgraph in a (hyper)graph.

Let $1\leq d,\ell\leq k-1$. For $n\in(k-\ell)\mathbb{N}$, define $h_d^{\ell}(k,n)$ to be the smallest integer $h$ such that every $n$-vertex $k$-graph $H$ satisfying $\delta_d(H)\geq h$ contains a Hamilton $\ell$-cycle.
Han and Zhao~\cite{2015Forbidding} gave the result that
\begin{equation}\label{Han}
h_d^{k-1}(k,n)\geq\left(1-\binom{t}{\lfloor t/2\rfloor}\frac{\lceil t/2\rceil^{\lceil t/2\rceil}(\lfloor t/2\rfloor+1)^{\lfloor t/2\rfloor}}{(t+1)^t}+o(1)\right)\binom{n}{t}
\end{equation}
where $d\in[k-1]$ and $t=k-d$.
In particular, $h_{d}^{k-1}(k,n)\geq(5/9+o(1))\binom{n}{2},(5/8+o(1))\binom{n}{3}$ for $k-d=2,3$.
Lang and Sanhueza-Matamala~\cite{Lang} conjectured that the minimum $d$-degree threshold for $k$-uniform tight Hamilton cycles coincides with the lower bounds in (\ref{Han}).
This leads to the following conjecture.
\begin{conjecture}
For every $k\geq4, \mu>0$, there exists $n_0$ such that the following holds for $n\geq n_0$.
Given a $k$-graph system $\textbf{G}=\{G_i\}_{i\in[n]}$, if $\delta_{k-3}(G_i)\geq(5/8+\mu)\binom{n}{3}$ for $i\in[n]$, then $\textbf{G}$ admits a rainbow Hamilton cycle.
\end{conjecture}
Furthermore, we believe the following holds.
\begin{conjecture}
For every $k,d, \mu>0$, there exists $n_0$ such that the following holds for $n\geq n_0$.
Given a $k$-graph system $\textbf{G}=\{G_i\}_{i\in[n]}$, if $\delta_{d}(G_i)\geq h_{d}^{k-1}(k,n)+\mu\binom{n}{d}$ for $i\in[n]$, then $\textbf{G}$ admits a rainbow Hamilton cycle.
\end{conjecture}

In fact, due to the whole proof of this paper, we believe that it is interesting to study rainbow Hamilton vicinities or rainbow Hamilton frameworks to determine the thresholds of Hamilton cycles.

\section{Acknowledgement}

This work was supported by the Natural Science Foundation of China (12231018,11871311,11901292) and Youth Interdisciplinary Innovation Group of Shandong University.


\end{document}